\documentclass[11pt]{article}

\usepackage{epsfig,psfrag,mathrsfs,amsmath,amsfonts,amssymb,amsthm}
\usepackage{color}
\usepackage{array}
\usepackage{lscape}
\usepackage{hyperref}
\usepackage{algorithm}
\usepackage{rotating}
\usepackage{algpseudocode}
\usepackage{cleveref}
\usepackage{multirow}
\usepackage{siunitx}
\usepackage{booktabs}
\usepackage{placeins}
\usepackage{seqsplit}
\usepackage{picinpar}
\usepackage{graphicx}
\usepackage{xcolor}
\usepackage{tikz}
\usetikzlibrary{calc,arrows.meta}
\usepackage[normalem]{ulem}

\newtheorem{remark}{Remark}[section]

\newtheorem{lemma}{Lemma}[section]
\newtheorem{example}{Example}[section]

\def\D{\mathrm{d}}

\newcommand{\ignore}[1]{}

\definecolor{inkblue}{RGB}{75,0,130}
\definecolor{afcwpurple}{RGB}{102,0,153}

\newcommand{\inputwithoutblue}[1]{%
  \begingroup
  \let\AFCWOriginalTextColor\textcolor
  \def\AFCWBlue{blue}%
  \renewcommand{\textcolor}[2]{%
    \def\AFCWCurrentColor{##1}%
    \ifx\AFCWCurrentColor\AFCWBlue
    \else
      \AFCWOriginalTextColor{##1}{##2}%
    \fi}%
  \input{#1}%
  \endgroup}

\title{Fully Discrete Active Flux Methods for the Euler Equations:
  Exploring Different Reconstructions and New Multi-Dimensional
  Evolution Operators}

\author{Christiane Helzel\thanks{Heinrich-Heine-University D\"usseldorf, Germany (christiane.helzel@hhu.de)},
M\'aria Luk\'a\v{c}ov\'a-Medvid'ov\'a\thanks{Johannes Gutenberg University Mainz, Germany (lukacova@uni-mainz.de)},
Amelie Porfetye\thanks{Heinrich-Heine-University D\"usseldorf, Germany (amelie.porfetye@hhu.de)},
Zhuyan Tang\thanks{Johannes Gutenberg University Mainz, Germany (zhtang@uni-mainz.de)}}

\date{}

\begin{document}
\maketitle

\begin{abstract}
We present new fully discrete Cartesian grid Active Flux methods for
the two-dimensional Euler equations of gas dynamics.
A key component of these methods is the evolution operator used
to evolve point-value degrees of freedom.
Therefore, we  derive several new, third-order accurate, multidimensional
evolution operators for  point-value updates of the linearized Euler
equations using the method of bicharacteristics.
We also investigate how the
reconstruction and choice of approximate evolution operator affect
the accuracy and stability of the resulting methods.
\end{abstract}

\section{Introduction}

The Active Flux (AF) method, originally introduced by Roe and his
former student Eymann
\cite{eymann2011active,eymann2013multidimensional} and subsequently
developed by Roe et al.\
\cite{FR2015,article:Roe2017,article:Roe2018,article:Roe2020,article:Roe2021,article:SR2023,proc:Roe2025}
is an innovative variant of a finite volume method for hyperbolic
conservation laws that is receiving increasing attention.
Its defining characteristic is the inclusion of point values
located along cell boundaries as independent, active degrees of
freedom (DoFs), which are evolved simultaneously alongside the classical cell
averages. Earlier contributions that focused on third-order accurate
methods for advection and acoustics revealed desirable properties.
For example, Cartesian grid AF methods for
acoustics with an exact evolution operator for the point
values have been shown to preserve vorticity \cite{article:BHKR2019}. Earlier,
Morton and Roe \cite{article:MR2001} speculated that
numerical methods for the Euler equations capable of controlling vorticity to
some extent may be less susceptible to anomalous solutions. Thus, AF methods emerge as promising candidates for investigating this
hypothesis.

In recent years, several new finite volume methods for the Euler equations
that use  point and cell average values as DoFs have emerged from Roe's original AF approach.
A crucial aspect of these methods is the evolution of the point
values. Motivated by the work of Abgrall and Barsukow \cite{article:AB2023}, who
proposed different concepts for constructing methods of arbitrary order of accuracy using the DoFs of AF
in combination with a method of lines approach,
rapid progress has been made in the development of so-called
generalized AF \cite{article:ABK2025,article:DBK2025,article:BKLNOR2026} and
Point-Average-Moment PolynomiAl-interpreted (PAMPA) \cite{article:AJLW2025,article:AL2026,article:ALOe2026}
schemes.
Using point and cell average values results in compact stencils for
the discretizations of spatial operators. Furthermore, the method of lines approach
circumvents the need for evolution operators to evolve point
values. However, each stage of the time-stepping method
increases its stencil, and stability requirements seem to
lead to severe time-step restrictions.

To construct fully discrete AF methods
when exact evolution operators for point-value updates are
unavailable, Chudzik et al.\ \cite{article:CHL2024} combined
the multidimensional evolution operators developed by
Luk\'a\v{c}ov\'a et al.\
\cite{article:LMW2000,article:LSW2002,article:LMW2004}
using the method of bicharacteristics with the reconstruction
and flux computation of the Cartesian grid
AF method.

For acoustics, linear stability analysis
\cite{article:CHK2021,article:CHL2024}
showed that Cartesian grid AF methods
using the third-order EG2 evolution operator from
\cite{article:LSW2002} have a
CFL limit of $0.279$,
compared with $0.5$ when the exact evolution operator is used.
To improve this stability limit, we recently revisited
the bicharacteristic formulation and derived new third-order
evolution operators for acoustics \cite{preprint:PTCHL2025}.

For the Euler equations, the construction in
\cite{article:CHL2024} yielded a fully discrete,
third-order accurate AF method
for the linearized system and a fully discrete extension
to sufficiently smooth solutions of the nonlinear system.
Building on this approach, \cite{article:CHP2026} developed
a third-order method for the two-dimensional nonlinear
Euler equations using appropriate local linearizations
and a corrected point-value update.
That work also introduced a limiting procedure for point
values and cell averages that preserves the positivity
of density and pressure and suppresses unphysical
oscillations near shocks.
In the present work, we extend the new acoustic evolution
operators to the linearized Euler equations
and incorporate them into this fully discrete
AF method for the Euler equations.

In addition to the evolution operators, we investigate
an alternative to the globally continuous, piecewise
quadratic reconstruction used in classical
AF methods.
In our acoustic study \cite{preprint:PTCHL2025},
we replaced this reconstruction with a third-order
central weighted essentially
non-oscillatory (CWENO) reconstruction based on neighboring
cell averages.
The resulting method has a wider stencil and exhibits
improved linear stability properties.
For the nonlinear Euler equations, the nonlinear weights
of CWENO provide a further motivation
for this choice, as they help suppress oscillations near
discontinuities.
However, the reconstruction is no longer globally continuous.
Comparing it with the classical AF reconstruction therefore provides
a setting in which to examine Roe's question
\cite{article:Roe2017}:
\emph{Is discontinuous reconstruction really a good idea?}
We refer to the AF method employing
CWENO reconstruction as
AFCW.
Through linear stability analysis and numerical experiments,
we examine the effects of the evolution operator and
reconstruction on stability, accuracy, and robustness.

The AFCW formulation also connects
this work to finite volume evolution
Galerkin (FVEG) methods.
Unlike classical AF methods,
it reconstructs from cell averages alone and uses evolved
point values for flux computation.
It can thus be interpreted as a third-order extension
of the second-order fully discrete
FVEG methods in
\cite{article:LMW2004,article:LSW2002}.

A related approach, recently proposed by Barsukow \cite{barsukow2026active},
approximates the point values in
fully discrete AF methods for the Euler equations
by additively splitting the Euler equations in primitive variables
into an acoustic part, which can be solved exactly, and a nonlinear transport part.
While this method 
introduces a splitting error, which reduces the 
accuracy to second order, numerical results show almost third-order
accuracy on grids of practical interest. 
In \cite{preprint:HP2026}, we  proposed approximating
the locally linearized
Euler equations in primitive variables
using a moving grid approach, once again making use of
the exact acoustic solver.  This method provides an exact solution of the
linearized Euler equations.  When this solution is used to
approximate the nonlinear Euler equations, a linearization error is
introduced, limiting the resulting method to second-order accuracy.
However, this linearization error was already analyzed in \cite{article:CHP2026}
and can be eliminated to achieve full third-order accuracy.
A related moving grid approach was simultaneously introduced by
Duraisamy \cite{duraisamy2026fully}.   Note that the fully discrete methods
proposed in
\cite{barsukow2026active,duraisamy2026fully,preprint:HP2026} make use
of the precise structure of the Euler equations, while the method of
bicharacteristics used in our paper is a more general approach for linear hyperbolic
systems.

The paper is organized as follows:  \Cref{sec:governing_equations} presents the
governing equations in the various forms employed by our methods, and
introduces an exact as well as an approximate evolution operator for the linearized Euler
equations.  
\Cref{sec:3} describes the fully discrete AF and AFCW
formulations, including the two reconstructions,
the point-value updates based on local linearization
and linear evolution operators, and the correction
required to retain third-order accuracy. \Cref{sec:4} is
devoted to the derivation of several new third-order accurate evolution
operators and investigates the stability of the resulting
AF methods for the linearized Euler equations.  \Cref{sec:efix} discusses difficulties in approximating transonic rarefaction waves and introduces an entropy fix to address them.
	Finally, \Cref{sec:numerical-results} presents and compares numerical results of the two
methods applied to the nonlinear Euler equations for a variety of test problems.

\section{Governing Equations and Evolution Operators}
\label{sec:governing_equations}

In this section, we review the two-dimensional Euler
equations in conservative and quasilinear forms.
The conservative form is used to update cell averages,
whereas point-value updates rely on a local linearization
of the quasilinear form.
We then introduce the resulting linearized system and
present its exact evolution operator, derived using
the method of bicharacteristics [26].
This representation provides the basis for constructing
the approximate evolution operators developed in this paper.

\subsection{Euler Equations and their Local Linearization}
The two-dimensional  Euler equations in conservative form are given by
\begin{equation}\label{eqn:conservative}
  \partial_t\mathbf{q} + \partial_x \mathbf{f(q)} + \partial_y \mathbf{g(q)}= \mathbf{0},
\end{equation}
with the conserved quantities $\mathbf{q}:\mathbb{R}^2\times \mathbb{R}^+ \rightarrow \mathbb{R}^4$ defined as
$\mathbf{q} = \left( \rho, \rho u, \rho v, E\right)^T$ and the flux functions
$\mathbf{f}(\mathbf{q}) = \left( \rho u, \rho u^2 + p, \rho u v, u(E+p)\right)^T$ and  $\mathbf{g}(\mathbf{q}) = \left( \rho v, \rho u v, \rho v^2 + p, v(E+p)\right)^T$.
The system is closed by the ideal gas equation of state $E = p/(\gamma - 1) + \frac{1}{2}\rho(u^2 + v^2)$. Here, $\rho$ denotes density, $u$ and $v$ are the velocity
components in the $x$- and $y$-directions, and $p$ is pressure.
The ratio of specific heats is denoted by $\gamma$; we use $\gamma=1.4$ unless stated otherwise. We will also make use of the local speed of sound
	$c:= \sqrt{\gamma \frac{p}{\rho}}.$

To facilitate the characteristic analysis and the application of the method of bicharacteristics, we rewrite the system in a quasi-linear form with respect to the primitive variables $\mathbf{u} = (\rho, u, v, p)^{T}$:
\begin{equation}\label{eq:euler_primitive}
    \partial_t \mathbf{u} + A_1(\mathbf{u}) \partial_x \mathbf{u} + A_2(\mathbf{u}) \partial_y \mathbf{u} = \mathbf{0},
\end{equation}
where $A_1(\mathbf{u})$ and $A_2(\mathbf{u})$ are 
\begin{equation}\label{eon:quasilinAB}
		A_1(\mathbf{u}) := \left( \begin{array}{cccc}u & \rho & 0 & 0 \\ 0 & u &
			0 & 1/\rho \\ 0 & 0 & u & 0 \\ 0 & \gamma p & 0 & u\end{array}\right), \quad  A_2(\mathbf{u}) :=  \left( \begin{array}{cccc}v & 0 & \rho & 0 \\ 0 & v & 0 & 0 \\ 0 & 0 & v & 1/\rho \\ 0 & 0 & \gamma p & v \end{array}\right).
	\end{equation}

Following the approach of \cite{article:LMW2004}, to derive the exact integral equations, we linearize the system \eqref{eq:euler_primitive} by freezing the Jacobian matrices at a specific evaluation point $\bar{P} = (\bar{x}, \bar{y}, \bar{t})$. {In our fully discrete AF method, these points naturally correspond to the independent DoFs at the corners of the cells and the midpoints of the cell interfaces, which are then used for flux integration. }

We denote the local constant reference state at this point $\bar{P}$ by ${\mathbf{u}}' = ({\rho'}, {u'}, {v'}, {p'})^{T}$ and the local speed of sound by ${c'} := \sqrt{\gamma {p'} / {\rho'}}$. By freezing the Jacobian matrices at ${\mathbf{u}'}$, we obtain the locally linearized system
\begin{equation}\label{eq:euler_linearized}
    \partial_t \mathbf{v} + A_1 \partial_x \mathbf{v} + A_2 \partial_y \mathbf{v} = 0,
  \end{equation}
with $A_1 := A_1({\mathbf{u}}')$ and $A_2 := A_2({\mathbf{u}}')$. Here, $\mathbf v$ denotes the solution of the locally
linearized system. In the evolution formulas below, we denote its components by $\rho$, $u$, $v$, and $p$. 

The eigenvalues of the resulting local matrix pencil $A = A_1\cos\theta + A_2\sin\theta$ are
\begin{equation}
    \lambda_1 = {u'}\cos\theta + {v'}\sin\theta - {c'}, \quad 
    \lambda_2 = \lambda_3 = {u'}\cos\theta + {v'}\sin\theta, \quad 
    \lambda_4 = {u'}\cos\theta + {v'}\sin\theta + {c'}.
\end{equation}
The eigenvalues $\lambda_1$ and $\lambda_4$ correspond physically to the acoustic (pressure) waves and $\lambda_2$ and $\lambda_3$ to the linearly degenerate entropy and shear waves  that are purely advected by the background flow.

\subsection{Exact Integral Equations via Bicharacteristics for the Linearized Problem}
Let $R({\mathbf{u}'})$ be the matrix whose columns are the right eigenvectors of $A({\mathbf{u}'})$, and let $L({\mathbf{u}'}) = R^{-1}({\mathbf{u}'})$ be the corresponding matrix of left eigenvectors. By multiplying the locally linearized system \eqref{eq:euler_linearized} from the left by $L({\mathbf{u}'})$, we project the primitive variables into the local characteristic space, thereby obtaining a quasi-diagonalized characteristic system. 

\begin{figure}[htbp]
    \centering
\includegraphics[width=0.5\textwidth]{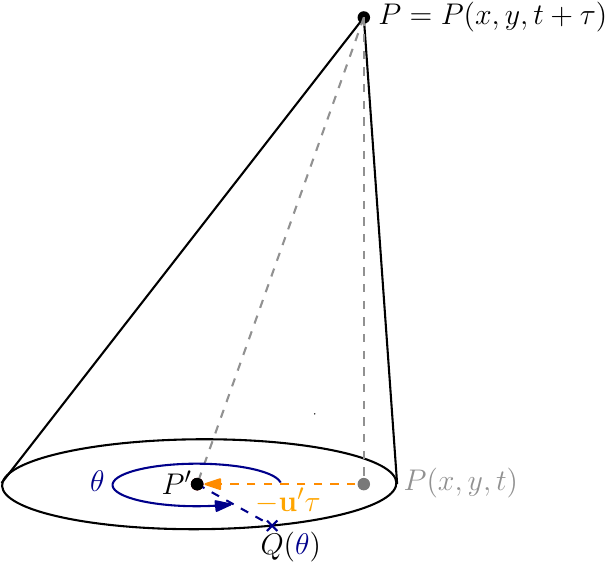}
    \caption{Bicharacteristics along the Mach cone through $P$ and $Q(\theta)$.}
    \label{fig:mach_cone}
\end{figure}

Consider the evaluation point $P=(x,y,t+\tau)$, with
$\tau>0$. The bicharacteristics through $P$ form an
advected Mach cone, as illustrated in Figure \ref{fig:mach_cone}.
At time $t$, their foot points lie on the circular
base of the cone. The space-time coordinates of the
base perimeter and its center are given by
\begin{align}
  Q(\theta) &= \left(x - ({u'} - {c'}\cos\theta)\tau, \, y - ({v'}
              - {c'}\sin\theta)\tau, \, t\right), \\
    P' &= \left(x - {u'}{\color{afcwpurple}\tau}, \, y - {v'}\tau, \, t\right).
\end{align}
Here, $Q(\theta)$ parameterizes the perimeter of the circular base of the advected acoustic cone, while $P'$ represents the purely advected center.

We write $\mathbf x=(x,y)$ for the spatial position,
$\mathbf a'=(u',v')$ for the constant background velocity,
and $\mathbf n(\theta)=(\cos\theta,\sin\theta)$ for the
unit direction vector.
For $P=(x,y,t+\tau)$ with $\tau>0$, integration over the
Mach cone mantle yields the following exact evolution
equations:
\begin{align}
\rho(P)
&= \rho(P')-\frac{p(P')}{c'^2}
 +\frac{1}{2\pi}\int_0^{2\pi}
 \left[
 \frac{p(Q)}{c'^2}
 -\frac{\rho'}{c'}
 \bigl(u(Q)\cos\theta+v(Q)\sin\theta\bigr)
 \right]\,\mathrm d\theta
 \notag\\
&\quad
 -\frac{\rho'}{2\pi c'}
 \int_0^{2\pi}\int_t^{t+\tau}
 S\!\left(
 \mathbf x-\bigl(\mathbf a'-c'\mathbf n(\theta)\bigr)
 (t+\tau-\tilde t),
 \theta,\tilde t
 \right)
 \,\mathrm d\tilde t\,\mathrm d\theta,
 \label{eq:exact_rho}
\\
u(P)
&= \frac{1}{2\pi}\int_0^{2\pi}
 \left[
 -\frac{p(Q)}{\rho'c'}\cos\theta
 +\bigl(u(Q)\cos\theta+v(Q)\sin\theta\bigr)\cos\theta
 \right]\,\mathrm d\theta
 \notag\\
&\quad
 +\frac{1}{2\pi}
 \int_0^{2\pi}\int_t^{t+\tau}
 \cos\theta\,
 S\!\left(
 \mathbf x-\bigl(\mathbf a'-c'\mathbf n(\theta)\bigr)
 (t+\tau-\tilde t),
 \theta,\tilde t
 \right)
 \,\mathrm d\tilde t\,\mathrm d\theta
 \notag\\
&\quad
 +\frac12 u(P')
 -\frac{1}{2\rho'}
 \int_t^{t+\tau}
 p_x\bigl(P'(\tilde t)\bigr)\,\mathrm d\widetilde t,
 \label{eq:exact_u}
\\
v(P)
&= \frac{1}{2\pi}\int_0^{2\pi}
 \left[
 -\frac{p(Q)}{\rho'c'}\sin\theta
 +\bigl(u(Q)\cos\theta+v(Q)\sin\theta\bigr)\sin\theta
 \right]\,\mathrm d\theta
 \notag\\
&\quad
 +\frac{1}{2\pi}
 \int_0^{2\pi}\int_t^{t+\tau}
 \sin\theta\,
 S\!\left(
 \mathbf x-\bigl(\mathbf a'-c'\mathbf n(\theta)\bigr)
 (t+\tau-\tilde t),
 \theta,\tilde t
 \right)
 \,\mathrm d\tilde t\,\mathrm d\theta
 \notag\\
&\quad
 +\frac12 v(P')
 -\frac{1}{2\rho'}
 \int_t^{t+\tau}
 p_y\bigl(P'(\tilde t)\bigr)\,\mathrm d\tilde t,
 \label{eq:exact_v}
\\
p(P)
&= \frac{1}{2\pi}\int_0^{2\pi}
 \left[
 p(Q)-\rho'c'
 \bigl(u(Q)\cos\theta+v(Q)\sin\theta\bigr)
 \right]\,\mathrm d\theta
 \notag\\
&\quad
 -\frac{\rho'c'}{2\pi}
 \int_0^{2\pi}\int_t^{t+\tau}
 S\!\left(
 \mathbf x-\bigl(\mathbf a'-c'\mathbf n(\theta)\bigr)
 (t+\tau-\tilde t),
 \theta,\tilde t
 \right)
 \,\mathrm d\tilde t\,\mathrm d\theta,
 \label{eq:exact_p}
\end{align}
where $Q=Q(\theta)$ denotes the foot point defined above and
\[
P'(\tilde t)
=
\bigl(
x-u'(t+\tau-\tilde t),
y-v'(t+\tau-\tilde t),
\tilde t
\bigr),
\qquad \tilde t\in[t,t+\tau],
\]
parameterizes the advected center line of the cone,
with $P'(t)=P'$ and $P'(t+\tau)=P$.
The source term $S$ is defined by
\begin{equation}
S(\mathbf x,\theta,\tilde t)
:=
c'\left[
u_x\sin^2\theta
-(u_y+v_x)\sin\theta\cos\theta
+v_y\cos^2\theta
\right],
\end{equation}
with all spatial derivatives evaluated at
$(\mathbf x,\tilde t)$. These exact integral representations elegantly separate the incoming states from the characteristic cone base and the continuous multi-dimensional wave interactions across the cone mantle.

Note that the exact expressions for $\rho(P)$ and $p(P)$ additionally satisfy
\begin{equation} \label{eq:rho-exact-prelation}
\rho(P) = \rho(P') + \frac{1}{{c'}^2} \left( p(P) - p(P')\right),
\end{equation}
which allows us to express the exact evolution operator in the form \eqref{eq:exact_u}-\eqref{eq:exact_p} and \eqref{eq:rho-exact-prelation}. This relation also holds for all approximate evolution operators considered in this paper. 

\subsection{The EG2 Approximate Evolution Operator}
While the exact integral representations \eqref{eq:exact_u}-\eqref{eq:exact_p} and \eqref{eq:rho-exact-prelation} rigorously capture the complete multi-dimensional wave physics, the integrals involving the cross-derivative source term $S$ over the Mach cone mantle are computationally intractable for a fully discrete scheme.

In the foundational FVEG methods \cite{article:LMW2004} and in our previous work on linear acoustics \cite{article:CHL2024,preprint:PTCHL2025}, we approximated these mantle integrals using the rectangle rule in time. This approximation yields the so-called EG2 evolution operator. For the linearized Euler equations, the EG2 operator updates the primitive variables using only the spatial data along the advected base circle $Q(\theta)$ and the center $P'$, leading to equations of the form
\begin{align}
    u(P) &= \frac{1}{\pi} \int_0^{2\pi} \left[ u(Q(\theta))\left(2\cos^2\theta - \frac{1}{2}\right) + 2v(Q(\theta))\sin\theta\cos\theta \right] {\rm d}\theta \notag \\
    &\quad - \frac{1}{\pi{\rho'}{c'}} \int_0^{2\pi} p(Q(\theta))\cos\theta \,{\rm d}\theta + \mathcal{O}(\tau ^3), \label{eq:eg2_u} \\[1.5ex]
    v(P) &= \frac{1}{\pi} \int_0^{2\pi} \left[ v(Q(\theta))\left(2\sin^2\theta - \frac{1}{2}\right) + 2u(Q(\theta))\sin\theta\cos\theta \right] {\rm d}\theta \notag \\
    &\quad - \frac{1}{\pi{\rho'}{c'}} \int_0^{2\pi} p(Q(\theta))\sin\theta \,{\rm d}\theta + \mathcal{O}(\tau ^3), \label{eq:eg2_v} \\[1.5ex]
    p(P) &= -p(P') + \frac{1}{\pi} \int_0^{2\pi} p(Q(\theta)) {\rm d}\theta \notag \\
    &\quad - \frac{{\rho'}{c'}}{\pi} \int_0^{2\pi} \left[ u(Q(\theta))\cos\theta + v(Q(\theta))\sin\theta \right] {\rm d}\theta + \mathcal{O}(\tau^3) \label{eq:eg2_p}
\end{align}
together with (\ref{eq:rho-exact-prelation}).

Although the EG2 operator is formally third-order accurate, recent stability analyses show that it yields restrictive CFL conditions when used with AF reconstruction (see \cite{article:CHL2024}). Furthermore, it fails to exactly preserve one-dimensional planar waves for piecewise quadratic initial data. These limitations motivate the novel approximate evolution operators presented in \Cref{sec:4}.

\subsection{Point Value Approximation of the Nonlinear Equations}\label{sec:correction}
The evolution operators discussed above approximate
the solution of a locally linearized Euler system.
When they are used to update point values of the nonlinear
Euler equations, the additional error introduced by
local linearization must also be taken into account.
In general, local linearization alone does not provide
a third-order accurate approximation of the nonlinear
solution.
In \cite[Theorem 2.1]{article:CHP2026}, a local linearization was described whose leading second-order term in the local linearization error has a particularly simple structure. We also use this linearization in this paper. 
Third-order accuracy can then be obtained by adding a correction term which eliminates the leading term of the linearization error as described in \cite[Corollary 2.2]{article:CHP2026}.

\section{Fully Discrete AF and AFCW Methods}\label{sec:3}


In this section, we describe the fully discrete AF and
AFCW methods for the two-dimensional Euler equations, with particular attention to their different
reconstructions.

Both formulations combine a piecewise quadratic
reconstruction, space-time flux quadrature, and point value
updates with the accuracy required for a third-order scheme.
We first present the cell-average updates and the flux
quadrature rules used in the two methods.
We then describe their reconstructions and local
linearization procedures for the point-value updates.
Finally, we discuss the simplifications used for
the linear stability analysis.

We consider a Cartesian grid with grid cell $(i,j)$ given by $[x_{i-\frac{1}{2}},x_{i+\frac{1}{2}}]\times [y_{j-\frac{1}{2}},y_{j+\frac{1}{2}}]$.
The cell average values of the conserved quantities $\mathbf{q}$ in cell $(i,j)$ at time $t_n$ are denoted by $\bar{\mathbf{Q}}_{i,j}^n$. The AF method also uses point values along the grid cell boundary as DoFs. These point values are located at the corners of the grid cell, i.e. at $(x_{i\pm \frac{1}{2}},y_{j\pm \frac{1}{2}})$ and at the edge midpoints, i.e. $(x_{i\pm \frac{1}{2}},y_j)$ and $(x_{i},y_{j\pm \frac{1}{2}})$. Point values of the conserved quantities  at time $t_n$ are denoted by $\mathbf{Q}_{i\pm\frac{1}{2},j\pm\frac{1}{2}}^n$,$\mathbf{Q}_{i,j\pm \frac{1}{2}}^n$, $\mathbf{Q}_{i\pm \frac{1}{2},j}^n$.  
Point values and cell averages in primitive variables
are denoted by $\mathbf{U}$ and $\overline{\mathbf{U}}$,
respectively, with the corresponding indices.

\subsection{Approximation of Cell Average Values}\label{sec:fvm-fluxes}
On a Cartesian mesh, the cell average values of the conserved quantities are updated by a finite volume scheme for (\ref{eqn:conservative}),
	\begin{equation}\label{eqn:fvm}
		\bar{\mathbf{Q}}_{i,j}^{n+1}
		=
		\bar{\mathbf{Q}}_{i,j}^{n}
		-
		\frac{\Delta t}{\Delta x}
		\left(
		F_{i+\frac12,j}
		-
		F_{i-\frac12,j}
		\right)
		-
		\frac{\Delta t}{\Delta y}
		\left(
		G_{i,j+\frac12}
		-
		G_{i,j-\frac12}
		\right),
	\end{equation}
	where the numerical fluxes are obtained by approximating the space-time flux integrals with a suitable quadrature rule, i.e.\
        \begin{equation}\label{eqn:fluxes}
          \begin{split}
            F_{i-\frac{1}{2},j} & \approx \frac{1}{\Delta t \Delta y}
			\int_{t_n}^{t_{n+1}}
			\int_{y_{j-\frac{1}{2}}}^{y_{j+\frac{1}{2}}}
			\mathbf{f}(\mathbf{q}(x_{i-\frac{1}{2}},y,t))
			\,{\rm d}y\,{\rm d}t, \\
                        G_{i,j-\frac{1}{2}}  & \approx \frac{1}{\Delta t \Delta x}
                        \int_{t_n}^{t_{n+1}}
			\int_{x_{i-\frac{1}{2}}}^{x_{i+\frac{1}{2}}} 	            \mathbf{g}(\mathbf{q}(x,y_{j-\frac{1}{2}},t))
			\,{\rm d}x\,{\rm d}t .\\
            \end{split}
          \end{equation}

          In both FVEG and AF methods, the numerical flux is computed using a quadrature rule. In \cite{article:LSW2002}, Simpson's rule in space was combined with the midpoint rule in time, resulting in 
          \begin{equation}
            F_{i-\frac{1}{2},j}^{FVEG} = \frac{1}{6}\left( \mathbf{f}(\mathbf{Q}_{i-\frac{1}{2},j-\frac{1}{2}}^{n+\frac{1}{2}}) + 4 \mathbf{f}(\mathbf{Q}_{i-\frac{1}{2},j}^{n+\frac{1}{2}}) +
              \mathbf{f}(\mathbf{Q}_{i-\frac{1}{2},j+\frac{1}{2}}^{n+\frac{1}{2}} ) \right).
          \end{equation}
          The point values $\mathbf{Q}_{i-\frac{1}{2},j-\frac{1}{2}}^{n+\frac{1}{2}}\approx \mathbf{q}(x_{i-\frac{1}{2}},y_{j-\frac{1}{2}},t_n+\frac{\Delta t}{2})$,  $\mathbf{Q}_{i-\frac{1}{2},j}^{n+\frac{1}{2}}$ and $\mathbf{Q}_{i-\frac{1}{2},j+\frac{1}{2}}^{n+\frac{1}{2}}$ have been computed
using an evolution operator derived by the method of bicharacteristics as explained in Section \ref{sec:governing_equations}.
The numerical flux  $G_{i,j-\frac{1}{2}}$ was computed analogously. In order to obtain a second-order accurate FVEG method, a second-order accurate approximation of the locally linearized Euler equations is appropriate, as this also provides a second-order accurate approximation of the nonlinear problem. 

In AF methods, the fluxes (\ref{eqn:fluxes}) are typically approximated using Simpson's rule in space and time, i.e.\
\begin{equation}\label{eqn:2dSimpson}
  \begin{split}
    F_{i-\frac{1}{2},j}^{AF} = &
\frac{1}{36}
			\Bigl(
			\mathbf{f}(\mathbf{Q}_{i-\frac{1}{2},j-\frac{1}{2}}^n)
			+4\mathbf{f}(\mathbf{Q}_{i-\frac{1}{2},j}^n)
			+\mathbf{f}(\mathbf{Q}_{i-\frac{1}{2},j+\frac{1}{2}}^n)
			\\
			& \qquad
			+4\mathbf{f}(\mathbf{Q}_{i-\frac{1}{2},j-\frac{1}{2}}^{n+\frac{1}{2}})
			+16\mathbf{f}(\mathbf{Q}_{i-\frac{1}{2},j}^{n+\frac{1}{2}})
			+4\mathbf{f}(\mathbf{Q}_{i-\frac{1}{2},j+\frac{1}{2}}^{n+\frac{1}{2}})
			\\
			& \qquad
			+\mathbf{f}(\mathbf{Q}_{i-\frac{1}{2},j-\frac{1}{2}}^{n+1})
			+4\mathbf{f}(\mathbf{Q}_{i-\frac{1}{2},j}^{n+1})
			+\mathbf{f}(\mathbf{Q}_{i-\frac{1}{2},j+\frac{1}{2}}^{n+1})
			\Bigr)
                      \end{split}
                    \end{equation}
 and analogously for $G_{i,j-\frac{1}{2}}$.                   
 In \cite{article:CHP2026}, third-order accurate AF methods for the Euler equations were developed, in which the point values in (\ref{eqn:2dSimpson}) at times $t_{n+\frac{1}{2}}$ and $t_{n+1}$ were approximated using the EG2 evolution operator, followed by correction terms to achieve third-order accuracy.  The point values at time $t_n$ are known DoFs. This makes Simpson's rule in time particularly attractive for AF methods involving globally continuous reconstruction,
since it uses the point values for both reconstruction and flux approximation.

If point values are only used for the flux computation and not for the reconstruction,
an alternative quadrature for third-order accurate flux computations is to use Simpson's rule in space and Gauss quadrature in time, with nodes \ignore{ and weights}
\begin{equation} \label{eq:Gauss-nodes}
  c_{1,2}^{\mathrm G}=\frac12\mp\frac{\sqrt3}{6}\ignore{,\qquad b_1=b_2=\frac12}.   
\end{equation}
This results in
\begin{equation}\label{eqn:gaussSimpson}
  \begin{split}
    F_{i-\frac{1}{2},j}^{SiGa} & = \frac{1}{12} \left( \mathbf{f}(\mathbf{Q}_{i-\frac{1}{2},j-\frac{1}{2}}^{n+\frac{1}{2}-\frac{\sqrt3}{6}}) + 4 \mathbf{f}(\mathbf{Q} _{i-\frac{1}{2},j}^{n+\frac{1}{2}-\frac{\sqrt3}{6}})+
      \mathbf{f}(\mathbf{Q} _{i-\frac{1}{2},j+\frac{1}{2}}^{n+\frac{1}{2}-\frac{\sqrt3}{6}}) \right.\\
& \qquad + \left. \mathbf{f}(\mathbf{Q}_{i-\frac{1}{2},j-\frac{1}{2}}^{n+\frac{1}{2}+\frac{\sqrt3}{6}}) + 4 \mathbf{f}(\mathbf{Q} _{i-\frac{1}{2},j}^{n+\frac{1}{2}+\frac{\sqrt3}{6}})+
      \mathbf{f}(\mathbf{Q} _{i-\frac{1}{2},j+\frac{1}{2}}^{n+\frac{1}{2}+\frac{\sqrt3}{6}}) 
    \right)
    \end{split}
\end{equation}
 and analogously for $G_{i,j-\frac{1}{2}}$. 
 
Both rules require new evolved point values at two temporal
quadrature nodes.
However, \eqref{eqn:gaussSimpson} does not require point values at $t_n$,
where the discontinuous reconstruction generally provides
multiple traces. We therefore use \eqref{eqn:gaussSimpson} for the AFCW method.

\subsection{Reconstruction}\label{sec:reconstruction}
To obtain third-order accurate point values,
we use piecewise quadratic reconstructions in primitive
variables.
We denote the reconstruction at time $t_n$ by
$\mathbf{u}_{\mathrm{rec}}^n(x,y)$.
The continuous AF reconstruction uses both point values
and cell averages, whereas the CWENO reconstruction
uses cell averages alone.
In both cases, the reconstruction supplies the initial
data for the evolution step.

Point values in primitive variables can be calculated from point values in conservative variables, and vice versa, without loss of accuracy. 
Since the transformation from conservative to primitive variables is nonlinear, a direct transformation of the conservative cell averages does not, in general, yield third-order accurate primitive cell averages.
For the AF reconstruction, third-order accurate cell average values in primitive variables can be computed without increasing the stencil. This involves first computing a third-order accurate point value of the conserved quantities at the midpoint of a grid cell by solving Simpson's quadrature formula for the
cell-center point value. Next, all point values are transformed into primitive variables. Finally, cell average values in primitive variables are computed using Simpson's rule  (see also \cite{article:CHP2026}).

For the CWENO reconstruction, the required cell averages of the primitive variables are obtained from the conservative cell averages as follows:
We first suppress the time index $n$
and write $m_x=\rho u$ and $m_y=\rho v$.
We compute the primitive variables directly
from the conservative cell averages:
\[
u_{i,j}^{(0)}
=
\frac{\bar{m}_{x,i,j}}{\bar{\rho}_{i,j}},
\qquad
v_{i,j}^{(0)}
=
\frac{\bar{m}_{y,i,j}}{\bar{\rho}_{i,j}},
\]
and
\[
p_{i,j}^{(0)}
=
(\gamma-1)
\left(
\bar{E}_{i,j}
-
\frac{\bar{m}_{x,i,j}^{\,2}+\bar{m}_{y,i,j}^{\,2}}
     {2\bar{\rho}_{i,j}}
\right).
\]
For any cell average \(\bar{a}_{i,j}\), introduce the centered increments
\[
\delta_x\bar{a}_{i,j}
=
\frac{\bar{a}_{i+1,j}-\bar{a}_{i-1,j}}{2},
\qquad
\delta_y\bar{a}_{i,j}
=
\frac{\bar{a}_{i,j+1}-\bar{a}_{i,j-1}}{2}.
\]
For smooth solutions, a Taylor expansion of the conservative-to-primitive transformation yields the following third-order approximations:
\begin{align}
\bar{u}_{i,j}
&=
u_{i,j}^{(0)}
+
\frac{1}{12\bar{\rho}_{i,j}^{\,2}}
\sum_{d\in\{x,y\}}
\delta_d\bar{\rho}_{i,j}
\left(
u_{i,j}^{(0)}\delta_d\bar{\rho}_{i,j}
-
\delta_d\bar{m}_{x,i,j}
\right),
\\
\bar{v}_{i,j}
&=
v_{i,j}^{(0)}
+
\frac{1}{12\bar{\rho}_{i,j}^{\,2}}
\sum_{d\in\{x,y\}}
\delta_d\bar{\rho}_{i,j}
\left(
v_{i,j}^{(0)}\delta_d\bar{\rho}_{i,j}
-
\delta_d\bar{m}_{y,i,j}
\right),
\\
\bar{p}_{i,j}
&=
p_{i,j}^{(0)}
-
\frac{\gamma-1}{24\bar{\rho}_{i,j}}
\sum_{d\in\{x,y\}}
\left[
\left(
\delta_d\bar{m}_{x,i,j}
-
u_{i,j}^{(0)}\delta_d\bar{\rho}_{i,j}
\right)^2
+
\left(
\delta_d\bar{m}_{y,i,j}
-
v_{i,j}^{(0)}\delta_d\bar{\rho}_{i,j}
\right)^2
\right].
\end{align}
Together with $\bar\rho_{i,j}$, these approximations
provide the primitive cell averages used in the
CWENO reconstruction.

Once the required primitive variable data have been obtained, the globally continuous, piecewise-quadratic AF reconstruction and the piecewise-quadratic CWENO reconstruction are constructed as described in \cite[Section~6]{article:HKS2019} and \cite[Section~4.2]{preprint:PTCHL2025}, respectively.

\subsection{Local linearization}\label{sec:LocalLinearization}
The choice of the local linearization, i.e.\ the definition of $\mathbf{u}'$ for the locally linearized problem (\ref{eq:euler_linearized}), plays an important role when constructing third-order accurate methods. On the one hand, we must choose the linearization such that we obtain a third-order accurate approximation for smooth solutions. On the other hand, we must ensure a stable approximation in the vicinity of discontinuities.
	In \cite{article:CHP2026}, Chudzik et al. proposed a local linearization strategy which, together with a correction term, yields a third-order accurate approximation of the primitive point values. Furthermore, it was shown that this approach is not sufficiently robust in the vicinity of shock waves and a modified strategy was proposed for such computations. Helzel and Porfetye adopted this approach in \cite{preprint:HP2026} and proposed a slight modification of the local linearization strategy, which is also used throughout the present paper. The following subsections describe the resulting local linearization strategies, first for the globally continuous AF reconstruction and subsequently for the CWENO reconstruction.

\subsubsection*{Local linearization for smooth solutions}\label{sec:loclinsmooth}
	A third-order accurate approximation of the primitive point values requires an appropriate choice of the local linearization states together with an analysis of the resulting linearization error. Both are described in \cite[Theorem 2.1]{article:CHP2026}. We assume that $\mathbf{u}$ is the solution of the nonlinear Euler system (\ref{eq:euler_primitive}) and that the system is linearized around a state  $$\mathbf{u}' = \mathbf{u}(\bar{x},\bar{y},t_n+ \tau /2) + \mathcal{O}(\tau ^2).$$
Furthermore, let $\mathbf{v}(\bar{x},\bar{y},t_n+ \tau)$ denote the solution of the linearized equations (\ref{eq:euler_linearized}).
        Then the leading term of the linearization error is given by 
	\begin{align*}
		\mathbf{u}(\bar{x},\bar{y},t_n+\tau) &-
		\mathbf{v}(\bar{x},\bar{y},t_n+\tau) = \\
		\frac{1}{2}  &\tau^2\left( A_1(\mathbf{u}) \frac{\partial A_1(\mathbf{u})}{\partial
			\mathbf{u}} \cdot (\partial_x \mathbf{u},\partial_x \mathbf{u}) +
		A_1(\mathbf{u}) \frac{\partial A_2(\mathbf{u})}{\partial
			\mathbf{u}} \cdot (\partial_x \mathbf{u},\partial_y \mathbf{u}) \right.\\
		\qquad \quad +& \left.  A_2(\mathbf{u}) \frac{\partial A_1(\mathbf{u})}{\partial
			\mathbf{u}} \cdot (\partial_y \mathbf{u},\partial_x \mathbf{u}) +
		A_2(\mathbf{u}) \frac{\partial A_2(\mathbf{u})}{\partial
			\mathbf{u}} \cdot (\partial_y \mathbf{u},\partial_y \mathbf{u}) \right) \bigg\vert_{(\bar{x},\bar{y},t_n)}\\
		+& \mathcal{O}(\tau ^3).
	\end{align*}
	This implies that a third-order accurate approximation of the point value $\mathbf{u}(\bar{x},\bar{y},t_n+ \tau )$ can be obtained by linearizing around the state  $\mathbf{u'}$ and adding a correction term, which discretizes the leading term of the linearization error. Details of this correction can be found in \cite[Corollary 2.2]{article:CHP2026}.  The particular choice of the local linearization states depends on the temporal quadrature used in the corresponding method. 
   In the following, let $\hat{\mathbf U}^n$ denote the locally averaged primitive state associated with a fixed spatial quadrature point at time $t_n$. Furthermore, we use the notation
$$\mathcal L(\tau;\mathbf U'),$$
for the point-value evolution over the time interval $\tau$ using the local linearization state $\mathbf U'$. 
For the temporal quadrature based on Simpson's rule, the required point values can be computed successively. Starting from $\hat{\mathbf U}^n$, the point value at the first intermediate time level is computed as
$$
\mathbf U^{n+\frac14}
=
\mathcal L\left(\frac{\Delta t}{4};\hat{\mathbf U}^n\right).
$$
This intermediate point value is not used directly in the flux computation, but serves as the local linearization state for the evolution to $t_n+\Delta t/2$. The resulting point value $\mathbf U^{n+\frac12}$ is then used as the local linearization state for the evolution to $t_n+\Delta t$. Thus, the point values are computed successively, with each newly computed point value serving as the local linearization state for the subsequent evolution.

\ignore{Table~\ref{Tab:linearization_states} summarizes the corresponding choices for the corner point values. The local linearization states for the interface midpoint values can be chosen analogously.}
    \ignore{For the temporal quadrature based on Simpson's rule, Table~\ref{Tab:linearization_states} displays one possible choice of local linearization states for the corner point values. The local linearization states for the interface midpoint values can be chosen analogously.
	\begin{table}[htbp]
		\centering
		\caption{Local linearization states used in smooth regions if Simpson's rule is used in space and time.}
		\label{Tab:linearization_states}
		\begin{tabular}{cc}
			\hline
			Linearization state & Target point value \\
			\hline
			
			$\displaystyle
			\frac14\left(
			\bar{\mathbf U}_{i-1,j-1}+
			\bar{\mathbf U}_{i,j-1}+
			\bar{\mathbf U}_{i-1,j}+
			\bar{\mathbf U}_{i,j}
			\right)$
			&
			$\mathbf U_{i-\frac12,j-\frac12}^{\,n+\frac14}$
			\\[3mm]

			$\mathbf U_{i-\frac12,j-\frac12}^{\,n+\frac14}$
			&
			$\mathbf U_{i-\frac12,j-\frac12}^{\,n+\frac12}$
			\\[3mm]

			$\mathbf U_{i-\frac12,j-\frac12}^{\,n+\frac12}$
			&
			$\mathbf U_{i-\frac12,j-\frac12}^{\,n+1}$
			\\
			\hline
		\end{tabular}
              \end{table}}

If, on the other hand, the discontinuous CWENO reconstruction is used, we employ, as described in Section~\ref{sec:fvm-fluxes}, a different quadrature rule for the temporal integration of the fluxes. While Simpson's rule is retained in space, the temporal integral is approximated by the two-point Gauss--Legendre rule with the nodes \(c_1^{\mathrm G}\) and \(c_2^{\mathrm G}\) defined in \eqref{eq:Gauss-nodes}.
In this case, the local linearization states are therefore constructed differently to provide suitable linearization states for the point-value evolution to the two Gauss nodes.\ignore{Consider a fixed spatial quadrature point.} Suitable linearization states for the two Gauss nodes are given as
\begin{equation}
\mathbf U'_k = (1-c_k^{\mathrm G})\hat{\mathbf U}^n+c_k^{\mathrm G} \mathcal{L}(\frac{\Delta t}{2};\hat{\mathbf U}^n),
\qquad k=1,2.
\label{eq:afcw-shared-linearization}
\end{equation}
Since $\mathcal{L}(\frac{\Delta t}{2};\hat{\mathbf U}^n)$ approximates the point value at $t_n+\Delta t/2$, the
interpolated state $\mathbf U'_k$ approximates the solution at the effective
time $t_n+\frac{c_k^{\mathrm G}\Delta t}{2}.$ This is the temporal midpoint between $t_n$ and the corresponding Gauss time
$t_n+c_k^{\mathrm G}\Delta t$, and therefore satisfies the local
linearization condition required for third-order accuracy. 
Unlike the method based on the globally continuous AF reconstruction, this procedure for computing the local linearization states is applied throughout. In contrast, the local linearization strategy used with the continuous reconstruction is modified below to improve the robustness of the resulting AF methods.

	\subsubsection*{Shock indicator}\label{sec:shock_indicator}
	While the local linearization strategy described in the previous paragraph is suitable for smooth solutions, it does not always lead to a sufficiently robust method in the vicinity of shock waves. We therefore employ a shock indicator to identify such regions. It is used to modify the local linearization states in the vicinity of shock waves and is also a key component of the limiting strategy for the AF method. Further details can be found in \cite[Section~3.2]{article:CHP2026}.
	
	The shock indicator is defined as $$\theta = \exp(-\kappa \phi ^{(1)} \phi ^{(2)}).$$ It detects shock waves and reflects their strength. $ \phi ^{(1)}$ identifies regions with strong pressure gradients using normalized finite-difference approximations of the second derivatives of the pressure. The quantity $ \phi ^{(2)}$ is given by finite-difference approximations of the divergence and the vorticity of the velocity field. $\kappa > 0 $ is a positive parameter that can be chosen freely. In this paper, we use $\kappa = 2$, i.e., the same value as in \cite{article:CHP2026,preprint:HP2026}. Consequently, $\theta$ is close to 1 in smooth regions and decreases towards 0 in the vicinity of shock waves.
	
	\subsubsection*{Local linearization near discontinuous solution structures}\label{sec:final-linearization}
	To improve the robustness of the method in the vicinity of shock waves, we modify the local linearization state, which was introduced above based entirely on accuracy requirements for smooth solutions. Numerical experiments indicate that using the average of the surrounding cell average values as the linearization state yields a robust method even for discontinuous test cases. Since this choice reduces the order of accuracy to second order in smooth regions, we combine both approaches. In the following, we restrict the discussion to the corner point values. Let $\mathbf{U}_{i-\frac{1}{2},j-\frac{1}{2}}^*$ denote the local linearization state used in the method based on the globally continuous AF reconstruction, as constructed in \Cref{sec:loclinsmooth}. In the vicinity of shock waves, this state is blended with the average of the four surrounding cell average values using a shock indicator. We denote by $\theta_{i-\frac{1}{2},j-\frac{1}{2}}$ the shock indicator associated with the grid point $(x_{i-\frac{1}{2}},y_{j-\frac{1}{2}})$. The local linearization state is then defined by
	\begin{equation}\label{eq:final_linearization}
		\begin{split}
			\mathbf{U}'_{i-\frac{1}{2},j-\frac{1}{2}}
			&= s\left(\theta_{i-\frac{1}{2},j-\frac{1}{2}}\right)
			\mathbf{U}_{i-\frac{1}{2},j-\frac{1}{2}}^*\\
			&\quad + \left(1-s\left(\theta_{i-\frac{1}{2},j-\frac{1}{2}}\right)\right)
			\frac{1}{4}\left(\bar{\mathbf{U}}_{i-1,j-1}^n+\bar{\mathbf{U}}_{i,j-1}^n+\bar{\mathbf{U}}_{i-1,j}^n+\bar{\mathbf{U}}_{i,j}^n \right), 
		\end{split}
	\end{equation}
	where $s:[0,1]\rightarrow[0,1]$ is a blending function.  Following \cite{article:CHP2026}, we choose the blending function $$s(\theta) = \max (0, \min(\frac{(\theta - \theta _0)}{(\theta _1 - \theta _0)},1)),$$ with $0 \le \theta _0 < \theta _1 \le 1$, i.e.  $\theta _0 = 0.9$ and $ \theta _1 =0.999$.
        The same blending procedure is applied analogously to the point values $\mathbf{U}_{i-\frac12,j}$ and $\mathbf{U}_{i,j-\frac12}$. However, as in   \cite{preprint:HP2026}, we compute the shock indicator at edge midpoints
that is used to limit the state of the local linearization
        from the indicator of the two neighboring point values via    $\theta_{i-\frac{1}{2},j} := \min \left\{\theta_{i-\frac{1}{2},j-\frac{1}{2}},\theta_{i-\frac{1}{2},j+\frac{1}{2}} \right\}$ and    $\theta_{i,j-\frac{1}{2}} := \min \left\{\theta_{i-\frac{1}{2},j-\frac{1}{2}},\theta_{i+\frac{1}{2},j-\frac{1}{2}} \right\}$. In particular, this modification of the shock indicator  provides a suitable linearization near transonic shock waves, a case that required  a separate treatment in \cite{article:CHP2026}.

In summary, the local linearization is based on the states proposed by the accuracy analysis for sufficiently smooth solutions, whereas near shock waves the linearization is shifted towards the average of the cell average values. Throughout this paper, we use the local linearization defined in \eqref{eq:final_linearization} for all
computations using the globally continuous AF reconstruction.

\subsection{AF Methods for the Linearized Euler Equations}
A simplified version of the method for the nonlinear Euler equations can be used to approximate the linear Euler equations (\ref{eq:euler_linearized}) with constant matrices $A_1$ and $A_2$. In this setting, the flux functions are $A_1 \mathbf{v}$ and $A_2 \mathbf{v}$, and a third-order accurate approximation of the point values yields a third-order accurate AF method. In particular, neither a correction term, as described in Section \ref{sec:correction}, nor local linearizations, as described in Section \ref{sec:LocalLinearization}, are required. Without limiting, the resulting AF methods are linear, allowing stability to be investigated through eigenvalue analysis as described in \cite{article:CHK2021,article:CHL2024,preprint:PTCHL2025}.
The next section presents stability results for AF methods that use the newly derived evolution operators.

\section{New Evolution Operators for the Linearized Euler Equations}\label{sec:4}
In this section we present several new third-order accurate evolution operators for the linearized Euler equations. The first operator, denoted by $\mathrm{EG}^{\text{quad}}$, provides an exact approximation of quadratic plane waves. This operator will later be used in Euler simulations together with the CWENO reconstruction, as our previous studies for acoustics (see \cite{preprint:PTCHL2025}) showed good stability and accuracy for the corresponding method.
In combination with the globally continuous AF reconstruction $\mathrm{EG}^{\text{quad}}$ leads to the same restrictive time step condition as EG2.  
Therefore, in Section \ref{sec:newEG-nu-delta}, we will derive additional evolution operators that improve the stability limit of AF methods using the globally continuous reconstruction.
\subsection{New Evolution Operator Reproducing Exactly Quadratic Plane Waves}
	\subsubsection{Exact 1D Solution for Quadratic Data}
	To establish this design target, we analyze the local wave propagation in the advected reference frame moving with the local background velocity $({u'}, {v'})$. To simplify the mathematical derivation, we assume ${u'} = {v'} = 0$ in this local frame. 
	We consider a planar wave aligned with the $x$-axis and define quadratic initial profiles for all primitive variables. For simplicity, the left state is set to zero:
	\begin{equation}
		\rho(x, 0) = \begin{cases} \rho^R x^2, & x > 0, \\ 0, & x \le 0, \end{cases} \quad 
		u(x, 0) = \begin{cases} u^R x^2, & x > 0, \\ 0, & x \le 0, \end{cases} \quad 
		p(x, 0) = \begin{cases} p^R x^2, & x > 0, \\ 0, & x \le 0, \end{cases}
	\end{equation}
	and $v(x,0) = 0$.
	
	The exact solution evaluates to:
	\begin{align}
		p(x, t) &= \begin{cases} 
			p^R(x^2 + {c'}^2 t^2) - 2{\rho'}{c'}^2 u^R x t, & x > {c'}t, \\ 
			\frac{1}{2}(p^R - {\rho'}{c'} u^R)(x+{c'}t)^2, & -{c'}t < x \le {c'}t, \\ 
			0, & x \le -{c'}t,
		\end{cases} \label{eq:exact_1d_p} \\[1.5ex]
		u(x, t) &= \begin{cases} 
			u^R(x^2 + {c'}^2 t^2) - \frac{2p^R}{{\rho'}} x t, & x > {c'}t, \\ 
			\frac{1}{2}\left(u^R - \frac{p^R}{{\rho'}{c'}}\right)(x+{c'}t)^2, & -{c'}t < x \le {c'}t, \\ 
			0, & x \le -{c'}t,
		\end{cases} \label{eq:exact_1d_u} \\[1.5ex]
		v(x, t) &= 0, \\[1.5ex]
		\rho(x, t) &= \begin{cases} 
			p^R t^2 - 2{\rho'}u^R x t + \rho^R x^2, & x > {c'}t, \\ 
			\frac{1}{2{c'}^2}(p^R - {\rho'}{c'}u^R)(x+{c'}t)^2 + \left(\rho^R - \frac{p^R}{{c'}^2}\right)x^2, & 0 < x \le {c'}t, \\ 
			\frac{1}{2{c'}^2}(p^R - {\rho'}{c'}u^R)(x+{c'}t)^2, & -{c'}t < x \le 0, \\ 
			0, & x \le -{c'}t.
		\end{cases} \label{eq:exact_1d_rho}
	\end{align}
	
	Evaluating these exact solutions at the advected origin $x=0$ yields the required reference states for parameterizing our approximate evolution operator: $p(0, \tau) = \frac{1}{2} (p^R - {\rho'}{c'} u^R) ({c'}\tau)^2$ and $u(0, \tau) = \frac{1}{2} \left(u^R - \frac{p^R}{{\rho'}{c'}}\right) ({c'}\tau)^2$ with $\tau \in (0,\Delta t]$.
	
	\subsubsection{Adjusting the EG2 Operator to Capture the Exact Solution}
	
	Building on the baseline EG2 framework introduced in Section \ref{sec:governing_equations}, we aim to systematically modify its integral formulation so it reproduces the 1D exact analytical solution for the quadratic data derived above. 
	If the standard EG2 operator is applied directly to this specific quadratic initial data, for instance, regarding the pressure \eqref{eq:eg2_p}, the continuous angular integral calculated along the base of the advective acoustic cone yields:
	\begin{equation}\label{eq:naive_integrals}
		\frac{1}{\pi} \int_0^{2\pi} u(Q(\theta)) \cos\theta \,{\rm d}\theta = \frac{4}{3\pi} ({c'}\tau)^2 u^R. 
	\end{equation}
    The corresponding relation in the $y$-direction follows
by considering analogous initial data depending on $y$.
	However, comparing this result to the exact 1D analytical solution in \eqref{eq:exact_1d_p}, it is evident that the standard EG2 operator fails to match the required coefficient of $1/2$. This discrepancy arises because the EG2 operator uses a simple time-rectangle rule for the mantle integrals, which is insufficiently accurate for resolving quadratic spatial variations.
	
	To address this and capture the exact solution, we draw inspiration from the derivation for the acoustic system \cite{preprint:PTCHL2025}. We apply a two-dimensional Taylor expansion to approximate the continuous base integrals using discrete finite differences evaluated at selected points on the circle. This yields:
	\begin{align}\label{eq:taylor_approx}
		\int_0^{2\pi} u(Q(\theta)) \cos\theta \,{\rm d}\theta &= \frac{\pi}{2} \left[ u(Q(0)) - u(Q(\pi)) \right] + \mathcal{O}(\tau ^3), \\
		\int_0^{2\pi} v(Q(\theta)) \sin\theta \,{\rm d}\theta &= \frac{\pi}{2} \left[ v\left(Q\left(\frac{\pi}{2}\right)\right) - v\left(Q\left(\frac{3\pi}{2}\right)\right) \right] + \mathcal{O}(\tau ^3).
	\end{align}
	
	To upgrade the EG2 operator, we replace the continuous angular integrals with a weighted combination of the original continuous integral and the newly derived discrete approximation, parameterized by a weight $\omega_1$. For the $x$-direction term, this adjustment takes the form:
	\begin{equation}
		\omega_1 \int_0^{2\pi} u(Q(\theta)) \cos\theta \,{\rm d}\theta + \left(\frac{1}{\pi} - \omega_1\right) \frac{\pi}{2} \left[ u(Q(0)) - u(Q(\pi)) \right].
	\end{equation}
	By substituting the evaluations for the quadratic data from \eqref{eq:naive_integrals} and \eqref{eq:taylor_approx} into this modified operator, the combined term evaluates to:
	\begin{equation}
		\omega_1 \left( \frac{4}{3} ({c'}\tau)^2 u^R \right) + \left(\frac{1}{\pi} - \omega_1\right) \frac{\pi}{2} \left( u^R ({c'}\tau)^2 \right) = \left( \frac{4}{3}\omega_1 + \frac{1}{2} - \frac{\pi}{2}\omega_1 \right) u^R ({c'}\tau)^2.
	\end{equation}
	
	In order to successfully capture the exact analytical behavior, this combined term must perfectly match the exact evolution term $\frac{1}{2} u^R ({c'}\tau)^2$ derived in \eqref{eq:exact_1d_u}. Equating the coefficients immediately dictates that $\omega_1 = 0$. An identical matching procedure for the transverse $y$-direction yields $\omega_2 = 0$.
	
	This rigorous parameter matching provides a profound structural conclusion: to achieve exactness for quadratic data and thereby construct the improved $\mathrm{EG}^{\text{quad}}$ operator, the continuous trigonometric integrals in the EG2 operator must be entirely replaced by their discrete finite difference counterparts given in \eqref{eq:taylor_approx}.	
	Thus, we obtain the approximate evolution operator $\mathrm{EG}^{\text{quad}}$   for the two-dimensional Euler equations
	\begin{align}
		p(P) &= -p(P') + \frac{1}{\pi} \int_0^{2\pi} p(Q(\theta)) {\rm d}\theta \notag \\
	&\quad - \frac{{\rho'}{c'}}{2} \left[ u(Q(0)) - u(Q(\pi)) \right] - \frac{{\rho'}{c'}}{2} \left[ v\left(Q\left(\frac{\pi}{2}\right)\right) - v\left(Q\left(\frac{3\pi}{2}\right)\right) \right] + \mathcal{O}(\tau^3), \label{egquad_p} \\[1.5ex]
			u(P) &= \frac{1}{\pi} \int_0^{2\pi} \left[ u(Q(\theta))\left(2\cos^2\theta - \frac{1}{2}\right) + 2v(Q(\theta))\sin\theta\cos\theta \right] {\rm d}\theta \notag \\
		&\quad - \frac{1}{2{\rho'}{c'}} \left[ p(Q(0)) - p(Q(\pi)) \right] + \mathcal{O}(\tau^3), \label{egquad_u} \\[1.5ex]
		v(P) &= \frac{1}{\pi} \int_0^{2\pi} \left[ v(Q(\theta))\left(2\sin^2\theta - \frac{1}{2}\right) + 2u(Q(\theta))\sin\theta\cos\theta \right] {\rm d}\theta \notag \\
		&\quad - \frac{1}{2{\rho'}{c'}} \left[ p\left(Q\left(\frac{\pi}{2}\right)\right) - p\left(Q\left(\frac{3\pi}{2}\right)\right) \right] + \mathcal{O}(\tau^3), \label{egquad_v} 
	\end{align}
	together with \eqref{eq:rho-exact-prelation}. The AFCW method considered below uses the $\mathrm{EG}^{\mathrm{quad}}$
evolution operator.

	\subsection{New Evolution Operator with Increased Stability}\label{sec:newEG-nu-delta}
    \ignore{
	As mentioned above, the stability limit of AFEG2 is significantly lower than the
necessary stability limit dictated by the compact stencil. 
        The aim is therefore to derive new approximate evolution operators that retain the third-order accuracy of EG2 while increasing the admissible $CFL$ number. To this end, we follow the approach introduced in \cite{preprint:PTCHL2025} for the acoustic equations and apply it to the linearized Euler equations.
	For this purpose, we use the following lemma.}
    As mentioned above, the AF method combining the exact evolution operator with the globally continuous reconstruction is stable for CFL numbers up to $0.5$, whereas the corresponding method using EG2, referred to as AFEG2 below, is stable only up to CFL$=0.279$. More generally, throughout the remainder of this paper, we use the notation AF followed by the name of an evolution operator to denote its combination with the globally continuous reconstruction. This section aims to derive new approximate evolution operators that retain the third-order accuracy of EG2 while increasing the stability limit of the resulting AF methods. To this end, we follow the approach introduced in \cite{preprint:PTCHL2025} for the acoustic equations and apply it to the linearized Euler equations.
For this purpose, we use the following lemma.

	\begin{lemma}\label{th:approximation_of_pv}
		Let $ f \in C^3(\Omega)$, where $ \Omega \subset \mathbb{R} ^2$ and contains the closed disk $\overline{B_R(x_0,y_0)}$. For $R>0$ define the circle parametrization $\mathbf{Q}_R(\theta)=(x_0+R\cos(\theta), y_0+R\sin(\theta))$, $\theta \in [0, 2\pi]$ and 
		\begin{equation*}
			\hat{f}_R(x_0,y_0):= \frac{1}{3}\left(\frac{4}{2\pi}\int_{0}^{2 \pi} f(\mathbf{Q}_{\frac{R}{2}}(\theta)) {\rm d} \theta - \frac{1}{2\pi}\int_{0}^{2 \pi} f(\mathbf{Q}_{R}(\theta)) {\rm d} \theta \right).
		\end{equation*}
		Then 
		\begin{equation*}
			f(x_0,y_0) = \hat{f}_R(x_0,y_0) + \mathcal{O} (R^3)
		\end{equation*}
		for $ R > 0 $ sufficiently small.
	\end{lemma}
	The proof is given in \cite[Lemma 4.1]{preprint:PTCHL2025}.
	To derive a new evolution formula for $p$, we modify the corresponding EG2 evolution formula by replacing the point value $p(P')$ with the third-order accurate approximation \(p_{\delta c' \tau}(P')\) obtained from Lemma~\ref{th:approximation_of_pv} using $R = \delta c' \tau$ with $\delta \ge 0$ and $\tau \in (0,\Delta t]$. The condition $\delta \cdot CFL \le 0.5$ must be satisfied.
	To preserve the relation between density and pressure given by \eqref{eq:rho-exact-prelation}, we modify the density evolution equation accordingly when changing the pressure evolution equation. The new equations of pressure and density read 
	\begin{equation} \label{eq:Eg2delta}
		\begin{split}
			p(P) &= - \hat{p}_{\delta c' \tau}(P') + \frac{1}{\pi} \int_{0}^{2 \pi} \left( p (\mathbf{Q}(\theta)) - \rho ' c' u(\mathbf{Q}(\theta))\cos \theta - \rho ' c' v(\mathbf{Q}(\theta))\sin \theta \right) \D \theta\\&  + \mathcal{O} (\tau ^3),\\
			\rho(P) &= \rho(P') + \frac{1}{c'^2}\left(p(P)-p(P')\right).
		\end{split}
	\end{equation}
	The combination of the new formulas for updating $p(P)$ and $\rho (P)$ and the original EG2 formulas for $u(P)$ and $v(P)$, given in \eqref{eq:eg2_u} and \eqref{eq:eg2_v},  defines a new family of evolution operators, which we denote by EG2$_\delta$. Note that for $\delta = 0$,
	EG2$_0$ reduces to EG2.
	
	Furthermore, the exact equation \eqref{eq:exact_u} for $u(P)$ contains the term $\frac{1}{2}u(P')$. This term cancels out in the derivation of EG2. It motivates us to add the term $0=u(P')-u(P')$ to the EG2 formula \eqref{eq:eg2_u} and replace one occurrence of $u(P')$ by an approximation obtained from Lemma~\ref{th:approximation_of_pv}. Numerical experiments indicate improved stability when the term $u(P')$ is replaced by $u_{\nu c'\tau}(P')$. We apply the same modification to the update of $v(P)$. The resulting formulas are
	\begin{equation} \label{eq:UundV}
		\begin{split}
			u(P) &= \frac{1}{\pi} \int_{0}^{2 \pi} \left( - \frac{p(\mathbf{Q}(\theta))}{\rho ' c'} \cos \theta + u(\mathbf{Q}(\theta)) (2 \cos ^2\theta - 0.5) +2 v(\mathbf{Q}(\theta))\sin \theta \cos \theta\right) \D \theta\\ &-u(P') + \hat{u}_{\nu c' \tau}(P')+ \mathcal{O} (\tau ^3), \\
			v(P) &= \frac{1}{\pi} \int_{0}^{2 \pi} \left( - \frac{p(\mathbf{Q}(\theta))}{\rho ' c'} \sin \theta + v(\mathbf{Q}(\theta)) (2 \sin ^2\theta - 0.5) +2 u(\mathbf{Q}(\theta))\sin \theta \cos \theta\right) \D \theta\\ &-v(P') + \hat{v}_{\nu c' \tau}(P')+ \mathcal{O} (\tau ^3).
		\end{split}
	\end{equation}
	Again, $\nu \ge 0$ and the condition $ \nu \cdot CFL \le 0.5$ must be satisfied. We refer to the evolution operator that uses \eqref{eq:Eg2delta} for updating $p(P)$ and $\rho (P)$ and \eqref{eq:UundV} for updating $u(P)$ and $v(P)$ as EG2$_{\delta,\nu}$. The corresponding AF methods are referred to as AFEG2$_{\delta}$
	and AFEG2$_{\delta, \nu}$, respectively. The parameters $\delta,\nu \ge 0$ may be chosen freely, provided that $\delta \cdot CFL \le 0.5$ and $\nu \cdot CFL  \le 0.5$. Since the stability limit of the AF methods is below $0.5$, these conditions are automatically fulfilled for $\delta,\nu \le 1$. Therefore, we initially restrict ourselves to the case $\delta,\nu \le 1$.
	
	\subsection{Computational Costs}\label{sec:sompCosts}
	A larger admissible time step does not necessarily make a numerical method more efficient if it significantly increases computational cost per time step. Therefore, besides the stability properties discussed in Section~\ref{sec:lin-stability}, it is also important to compare the computational costs of the different evolution operators. As discussed below, one particular choice of the free parameters $\delta, \nu$ requires only a small amount of additional work compared to EG2. As in the acoustic case, the computational costs of EG$^{\text{quad}}$ are comparable to those of EG2. Therefore, we restrict our considerations to EG2$_\delta$ and EG2$_{\delta,\nu}$ in the following.
	
	The update of point values using EG2 involves integrals of the form
	\begin{equation}
		\int_{0}^{2\pi}w^{\mathrm{rec}}(\mathbf{Q}(\theta))
		\sin^{\alpha}\theta\cos^{\beta}\theta\,{\rm d}\theta,
	\end{equation}
	where $w\in\{p,u,v\}$ and $\alpha,\beta\in\{0,1,2\}$. No integration is required for the density component.
	
\begin{figure}
		\centering
		\includegraphics[scale=0.65]{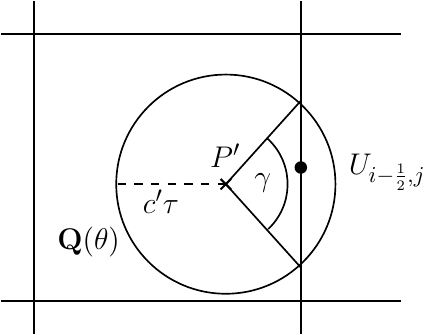}
		\includegraphics[scale=0.65]{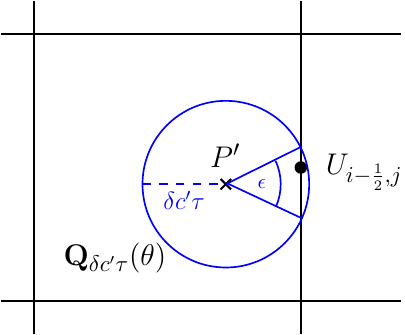}
		\includegraphics[scale=0.65]{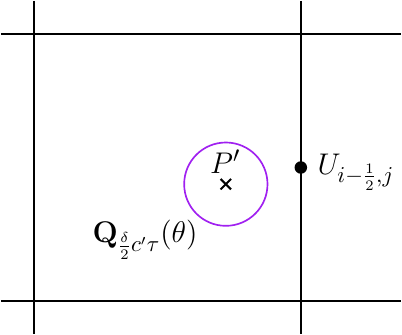}
		\caption{\label{fig:EG-Up} Update of a left point value using the method of bicharacteristics.}
	\end{figure}

	As an example, we consider the update of the left face point value
	$\mathbf{U}^n_{i-\frac12,j}$. Assuming that the bicharacteristic circle intersects only the two neighboring cells $(i-1,j)$ and $(i,j)$, the integrals take the form
	\begin{equation}\label{eqn:EG2-comp}
		\begin{aligned}
			&\sum_{k,l\in\{0,1,2\}}\tilde{w}_{i,j}^{k,l}
			\int_{-\frac{\gamma}{2}}^{\frac{\gamma}{2}}
			\sin^{\alpha+k}\theta\cos^{\beta+l}\theta\,{\rm d}\theta \\
			+{}&
			\sum_{k,l\in\{0,1,2\}}\tilde{w}_{i-1,j}^{k,l}
			\int_{\frac{\gamma}{2}}^{2\pi-\frac{\gamma}{2}}
			\sin^{\alpha+k}\theta\cos^{\beta+l}\theta\,{\rm d}\theta,
		\end{aligned}
	\end{equation}
	where $\tilde{w}_{i,j}^{k,l}$, $k,l\in\{0,1,2\}$, are linear combinations of the degrees of freedom of cell $(i,j)$ depending on $x$, $y$, $u'$, $v'$, $c'$ and $\tau$. We compute these coefficients exactly using the SymPy module in Python. The angle $\gamma$ also depends on $x$, $y$, $u'$, $v'$, $c'$ and $\tau$. The integrals for the remaining point values are obtained analogously.
	Compared to EG2, the evolution operators EG2$_\delta$ and EG2$_{\delta,\nu}$ additionally involve integrals of the form
	\begin{equation}
		\begin{split}
			&\int_{0}^{2\pi}
			w^{\mathrm{rec}}(\mathbf{Q}_{\varrho c'\tau}(\theta))
			\sin^{\alpha}\theta\cos^{\beta}\theta\,{\rm d}\theta,\\
			&\int_{0}^{2\pi}
			w^{\mathrm{rec}}(\mathbf{Q}_{\frac{\varrho}{2}c'\tau}(\theta))
			\sin^{\alpha}\theta\cos^{\beta}\theta\,{\rm d}\theta,
		\end{split}
	\end{equation}
	with $\varrho=\delta,\nu$. These integrals can again be evaluated in the form of \eqref{eqn:EG2-comp}. In general, however, neither the previously computed coefficients $\tilde{w}_{i,j}^{k,l}$ nor the corresponding trigonometric integrals can be reused. Consequently, we must evaluate six additional integrals for each point-value update. As an example, Figure~\ref{fig:EG-Up} illustrates three possible configurations of the bicharacteristic circles used to evolve the point value of one of the primitive variables.
	
	A particularly attractive special case is obtained by choosing
	$\delta=\nu=1$. Then three of the additional circular averages coincide with averages already required in the original EG2 evolution operator.
\ignore{
        For the corresponding evolution operators for the acoustic equations, the situation is even simpler. The configuration of the  circles used during a time step is fixed and therefore independent of the local solution. As a result, the required trigonometric integrals can be precomputed, and the corresponding coefficients can be reused for all integrals involved in the update of a given point value. Consequently, this reduces the computational cost. We refer to \cite[Section~4.5]{preprint:PTCHL2025} for further details.}

	\subsection{Investigation of linear stability}\label{sec:lin-stability}
	In this section, we analyze the linear stability of AF methods when applied to the linearized Euler equations. These results indicate the stability of the corresponding methods for the Euler equations.
	In \cite{article:CHL2024}, a fully discrete AF method based on EG2 for the linearized Euler equations was proposed. However, the corresponding eigenvalue analysis was only performed for the AF method applied to the acoustic equations, confirming stability for CFL$\le 0.279$.
	In this section we will conduct a stability analysis for the linearized Euler equations using EG2, as well as several new evolution operators. 
	To analyze the linear stability, we consider discretizations on a square domain with periodic boundary conditions in $x-$ and $y-$ direction. The domain is discretized using a grid with $m\times m$ grid cells. The linear method can then be written in the form
	\[
	U^{n+1}=B(\Delta t)U^n,
	\]
	where the vector $U^n$ contains all degrees of freedom on the two-dimensional Cartesian grid, and the matrix $B$ describes the evolution of the degrees of freedom over one time step of size $\Delta t$. A necessary condition for stability is that all eigenvalues of $B(\Delta t)$ lie within the unit circle.
	\ignore{\sout{
	\cite{article:CHL2024} showed that the AF method for the acoustic equations based on EG2 is stable under the condition $CFL\le0.279$. While EG$^{\mathrm{quad}}$ does not improve this stability limit, \cite{preprint:PTCHL2025} showed that EG2$_{0.7}$ increases it to $CFL\le0.418$ and EG2$_{0.8,0.2}$ further to $CFL\le0.440$.}}
        We now investigate the eigenvalues of the various AF methods for the linearized Euler equations with fixed background velocities $u'=v'=0$.
	
	\begin{figure}
		\centering
		\includegraphics[width=0.24\textwidth]{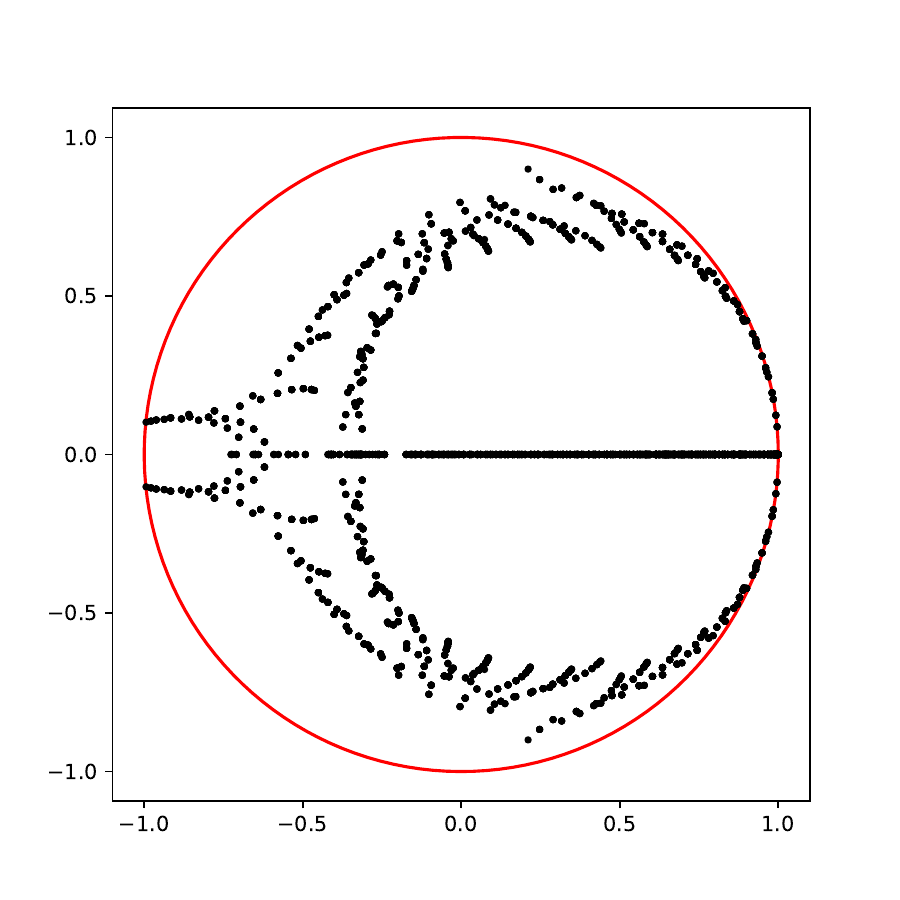}
		\includegraphics[width=0.24\textwidth]{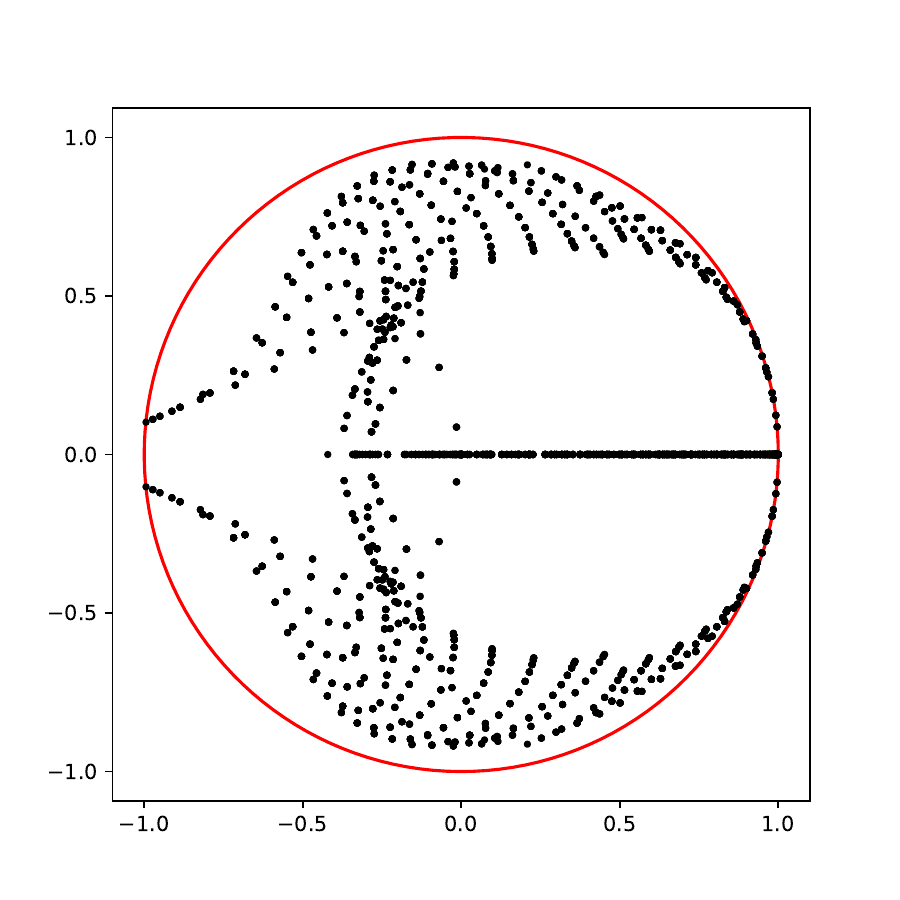}
		\includegraphics[width=0.24\textwidth]{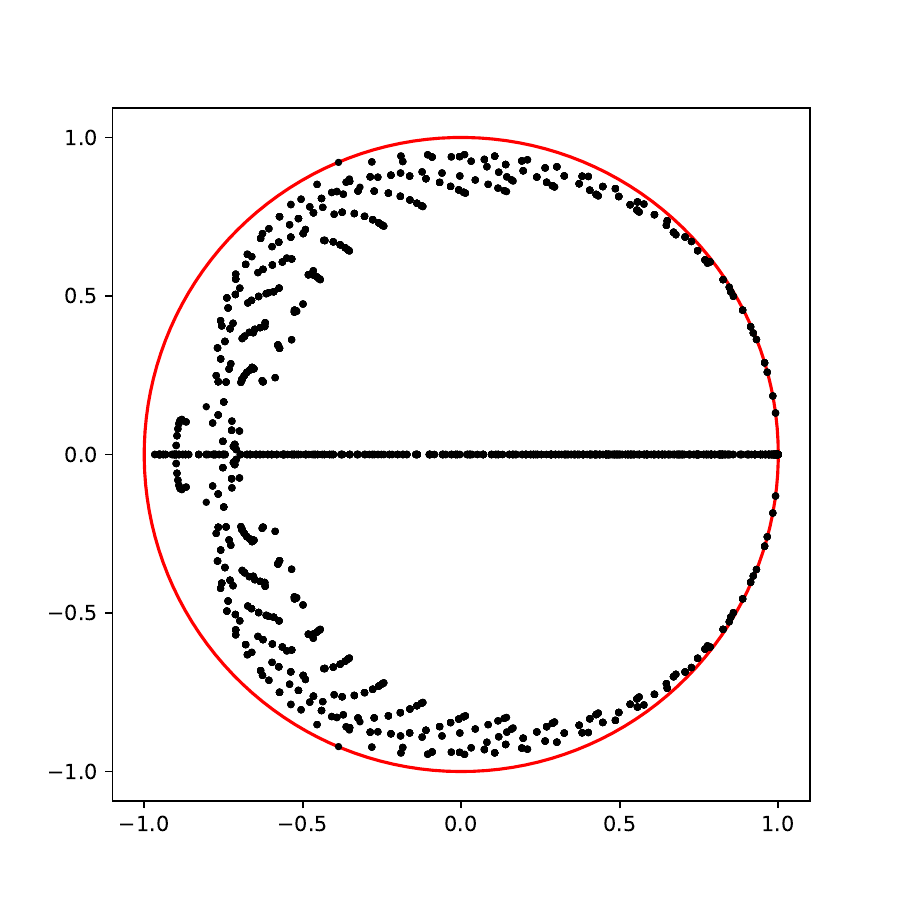}
		\includegraphics[width=0.24\textwidth]{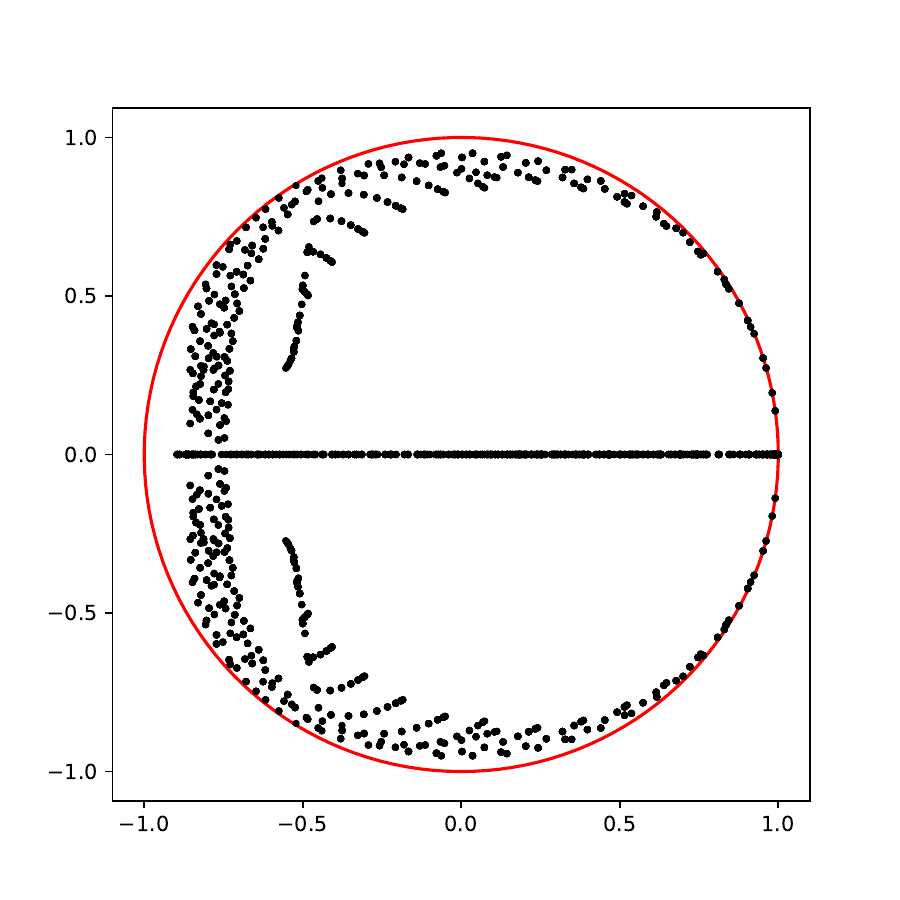}
		\caption{\label{fig:eigenvaluesacoustics} Eigenvalues of the matrix $B$ for the AF method for linearized Euler using, from left to right, EG2 and EG$^{\mathrm{quad}}$ with $CFL=0.279$, EG2$_{0.7}$ with $CFL=0.418$, and EG2$_{1.0,0.4}$ with $CFL=0.440$. The eigenvalues were computed on a $20^2$ grid with periodic boundary conditions.}
	\end{figure}
	
	In Figure~\ref{fig:eigenvaluesacoustics} we present the eigenvalues of the AF methods based on EG2, EG$^{\mathrm{quad}}$, EG2$_{0.7}$ and EG2$_{1.0,0.4}$. AFEG2 is stable for $CFL\le0.279$. The same stability limit is obtained for AFEG$^{\mathrm{quad}}$. Using AFEG2$_{0.7}$ increases the stability limit to $CFL\le0.418$, while AFEG2$_{1.0,0.4}$ further increases it to $CFL\le0.440$.
	As discussed above, the parameter choice $\delta=\nu=1$ is particularly attractive from a computational efficiency viewpoint, since only three additional integrals must be evaluated. However, this choice does not provide the largest admissible time step. While AFEG2$_{1.0}$ is stable up to $CFL\le0.418$, the stability limit of AFEG2$_{1.0,1.0}$ is reduced to $CFL\le0.378$, whereas AFEG2$_{1.0,0.4}$ remains stable up to $CFL\le0.440$. This observation suggests that the correction terms $u(P')-\hat{u}_{\nu c' \tau}$ and  $v(P')-\hat{v}_{\nu c' \tau}$
	significantly influence the stability of the fully discrete method. To control their contribution while retaining the computational advantages of the choice $\delta=\nu=1$, we introduce a weighting parameter $\omega\in[0,1]$ and replace \eqref{eq:UundV} by
	\begin{equation} \label{eq:UundVnew}
		\begin{split}
			u(P) &= \frac{1}{\pi} \int_{0}^{2 \pi} \left( - \frac{p(\mathbf{Q}(\theta))}{\rho ' c'} \cos \theta + u(\mathbf{Q}(\theta)) (2 \cos ^2\theta - 0.5) +2 v(\mathbf{Q}(\theta))\sin \theta \cos \theta\right) \D \theta\\ &-\omega(u(P') - \hat{u}_{\nu c' \tau}(P'))+ \mathcal{O} (\tau ^3), \\
			v(P) &= \frac{1}{\pi} \int_{0}^{2 \pi} \left( - \frac{p(\mathbf{Q}(\theta))}{\rho ' c'} \sin \theta + v(\mathbf{Q}(\theta)) (2 \sin ^2\theta - 0.5) +2 u(\mathbf{Q}(\theta))\sin \theta \cos \theta\right) \D \theta\\ &-\omega(v(P') - \hat{v}_{\nu c' \tau}(P'))+ \mathcal{O} (\tau ^3)
		\end{split}
	\end{equation}
	with $\omega\in[0,1]$. The resulting evolution operators, defined by \eqref{eq:Eg2delta} and \eqref{eq:UundVnew}, are denoted by EG2$_{\delta,\nu}^{\omega}$. Note that $\mathrm{EG2}_{\delta, \nu}^0 = \mathrm{EG2}_{\delta}$ and $\mathrm{EG2}_{\delta, \nu}^1 = \mathrm{EG2}_{\delta,\nu}$. 
	Table~\ref{Tab:CFL_values_EG2-delta-nu-omega} lists the maximum admissible CFL numbers of AFEG2$_{1.0,1.0}^{\omega}$ for different values of $\omega$. The largest stability limit occurs at $\omega=0.4$, yielding $CFL=0.444$. Hence, AFEG2$_{1.0,1.0}^{0.4}$ combines the computational advantages of the choice $\delta=\nu=1$ with a stability limit that is slightly higher than that of AFEG2$_{1.0,0.4}$.
    \ignore{
	\begin{remark}
 Modifying only the evolution formulas for $\rho$ and $p$ yields exactly the same stability limits as for the acoustic equations. Differences arise only when we also modify the evolution formulas for $u$ and $v$. For comparison, we refer to Table~7 in \cite{preprint:PTCHL2025}, which reports the corresponding limits for the AF methods applied to the acoustic equations.
 Even in the case $u'=v'=0$, the stability limits for the acoustic equations and the linearized Euler equations are not identical. In this case, the optimal parameter choice changes slightly. For example, AFEG2$_{0.8,0.2}$ applied to the acoustic equations is stable up to $CFL = 0.44$, whereas the corresponding method for the linearized Euler equations is stable only up to $CFL = 0.432$. On the other hand, AFEG2$_{1.0,0.4}$ applied to the acoustic equations is stable only up to $CFL = 0.43$, while the corresponding method for the linearized Euler equations remains stable up to $CFL = 0.44$. This suggests that the additional modification of the velocity evolution affects the stability properties and leads to a slightly different optimal choice of the parameters $\delta$ and $\nu$. A possible explanation is that the modified evolution of the velocity components changes their coupling with the pressure and density equations, thereby altering the overall stability behavior of the fully discrete scheme. 
	\end{remark}}
	
	\begin{table}[!ht]
		\caption{Maximum admissible CFL numbers for a stable method with EG2$_{1.0,1.0}^{\omega}$, $\omega = 0, 0.1, \dots, 1$.}
		\vspace*{0.1cm}
		\sisetup{round-mode=places,round-precision=3}
		\hspace*{+0.2cm}
		\begin{minipage}{0.95\textwidth}
			\centering
			\begin{tabular}{l *{11}{S[table-column-width=0.8cm,table-text-alignment=center]}}
				\toprule
				{$\omega$}& {{0.0}} & {0.1} & {0.2} & {0.3} & {0.4} & {0.5} & {0.6} & {0.7} & {0.8} & {0.9} & {1.0} \\
				\midrule
				CFL & \num{0.418} & \num{0.426} & \num{0.432} & \num{0.439} & \num{0.444} & \num{0.442} & \num{0.430} & \num{0.426} & \num{0.402} & \num{0.387} & \num{0.378} \\
				\bottomrule
			\end{tabular}
			\label{Tab:CFL_values_EG2-delta-nu-omega}
		\end{minipage}
	\end{table}
	
	The presented analysis provides only an indication of the stability properties of the corresponding AF methods for the nonlinear Euler equations. Nevertheless, the numerical results suggest that the stability limits of the Euler methods are indeed close to those predicted by the linear analysis. 
\section{Approximation of Transonic Rarefaction Waves}\label{sec:efix}
Our AF methods for the nonlinear Euler equations use local
linearizations and truly multi-dimensional approximate evolution
operators for  the linear Euler equations. This approach raises
the question of whether all nonlinear solution structures can be
correctly represented. It is well-known that classical finite volume
methods, which are based on local linearizations and the solution
of one-dimensional Riemann problems, require an entropy fix.
Roe (see \cite{article:Roe2017})
pointed out that, by avoiding discontinuous reconstructions,
AF eliminates the need for Riemann problems and hence avoids
the entropy fixes required by certain Riemann
solvers. 

For first- and second-order accurate FVEG methods, Luk\'a\v{c}ov\'a et al.\
\cite{article:LM2010}
observed the typical unphysical jump discontinuities for
one-dimensional transonic rarefaction waves and proposed modifications
of the evolution operators which changed the geometry of the circles
to avoid the problem.  \ignore{add more references, \cite{article:Bollermann2011}}
They also observed that the numerical problems
are much less pronounced in second-order accurate FVEG methods.
Thus, it was not surprising that the third-order accurate AF method
performed well in approximating one-dimensional
Riemann problems involving transonic rarefaction waves, as demonstrated in  \cite{article:CHP2026}. 

For the AF method, it was previously observed that the
approximation of nonlinear problems requires a special treatment when the
characteristic speed changes sign, to avoid a decoupling of the
degrees of freedom
\cite{article:HKS2019,article:Barsukow2021,article:DBK2025,article:CHP2026}. However,
a decoupling was not observed 
for transonic rarefaction waves, but rather for shocks across which the
characteristic speeds change sign. We have observed that the choice of
the local linearization is essential to avoid related
problems in fully discrete AF methods, cf.~Section \ref{sec:LocalLinearization}, where
a detailed description of our local linearization strategy is presented.   
Recently, Abgrall and Liu \cite{preprint:AbgrallLiu2026} proposed an entropy correction for semi-discrete AF-type methods based on a blending of the high-order numerical flux with a first-order entropy-satisfying flux. Their approach enforces the entropy inequality through updates to the cell averages, without requiring an entropy correction to the point-value update.

We now discuss a test problem which illustrates the importance of an
entropy-fix in a two-dimensional Euler simulation when using the AF
method with globally continuous reconstruction and the EG2 evolution
operator. Numerical problems are
observed in flow regimes with $|u'|<c'$, $|v'|<c'$ and
$\sqrt{u'^2+v'^2}>c'$, thus a situation that can be related to transonic
rarefaction waves. However, we do not observe an unphysical
approximation of a rarefaction wave but instead observe that AFEG2 is
unable to accurately approximate the Shu vortex problem  for
certain background speeds. We 
show that a version of the new EG2$_{\delta, \nu}^\omega$ evolution
operator avoids this problem by introducing additional circles which
are larger than the circle suggested by EG2.  
It has been reported in \cite{article:Spiegel2015} that the Shu vortex problem often leads to problems
when used to predict the order of convergence of high-order numerical methods
with small dissipation. 

In this context, we refer to \cite{LukacovaTangYuan2026}, where we established entropy stability for a modified
semi-discrete FVEG method with additional numerical diffusion.

The following example presents the Shu-vortex problem.

\begin{example}[Isentropic vortex]
	\label{ex:shu_vortex}
	
	We consider the advection of a smooth isentropic vortex. The initial condition is given by
	\begin{align*}
		\rho(x,y,0) &= \left(1+\delta T\right)^{\frac{1}{\gamma-1}},\\
		u(x,y,0) &= 1-\frac{\sigma}{2\pi}y\exp\!\left(\frac{1-r}{2}\right),\\
		v(x,y,0) &= 1+\frac{\sigma}{2\pi}x\exp\!\left(\frac{1-r}{2}\right),\\
		p(x,y,0) &= \left(1+\delta T\right)^{\frac{\gamma}{\gamma-1}},
	\end{align*}
	where
	\[
	r=x^2+y^2,\qquad
	\delta T=-\frac{(\gamma-1)\sigma^2}{8\gamma\pi^2}\exp(1-r),
	\]
	and $\sigma=5$. The computational domain is defined as $[-10,10]^2$ with periodic boundary conditions in both coordinate directions. 
\end{example}

\begin{figure}
	
	\makebox[\textwidth][c]{\includegraphics[width=0.25\textwidth]{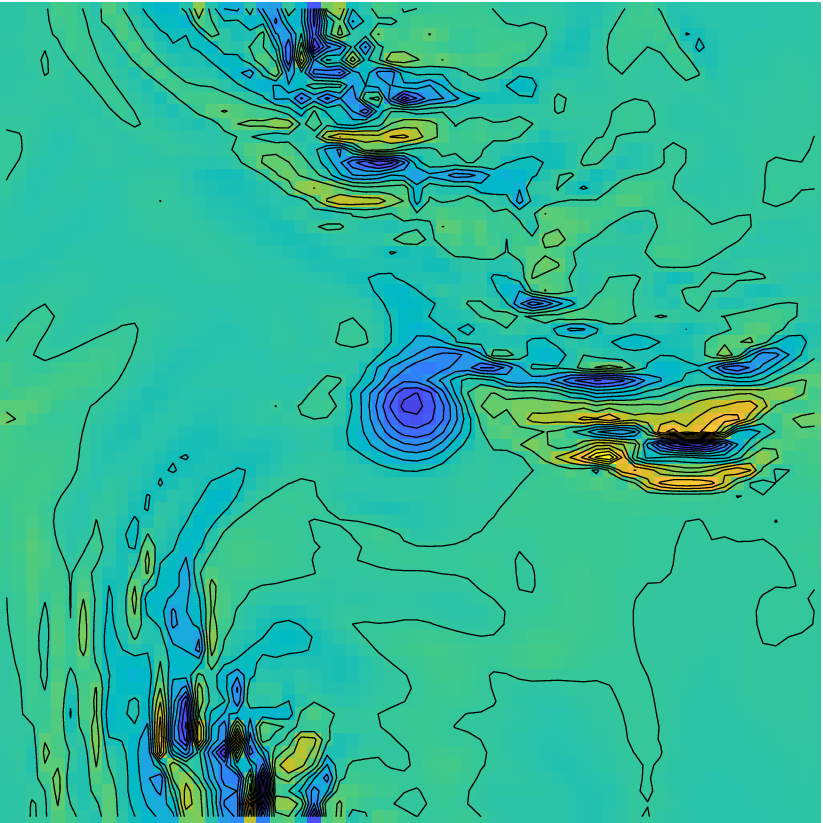}
		\includegraphics[width=0.25\textwidth]{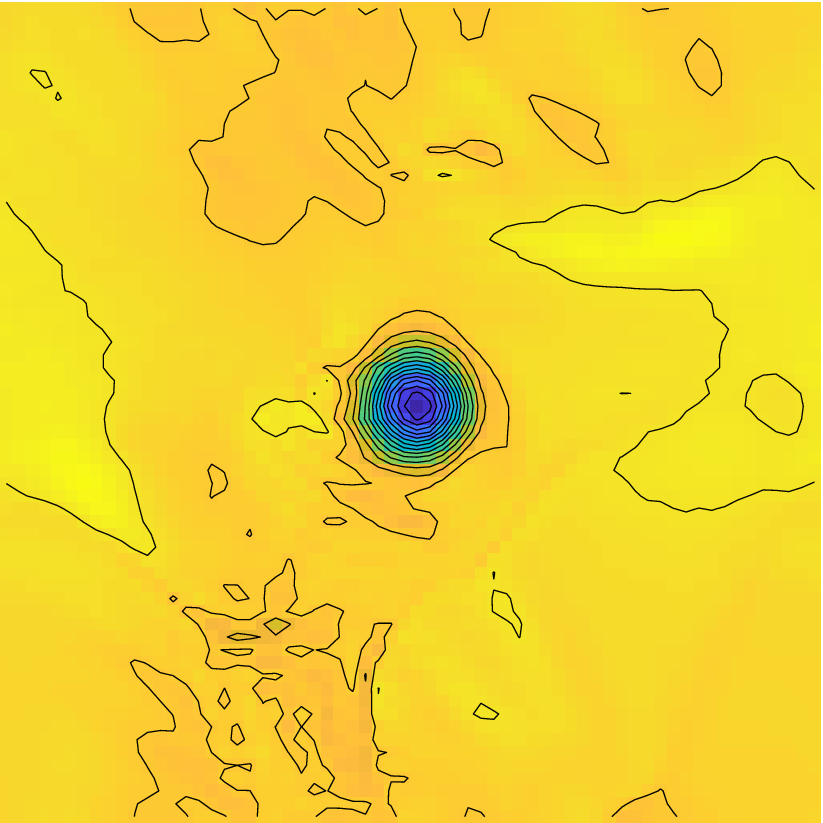}
		\includegraphics[width=0.25\textwidth]{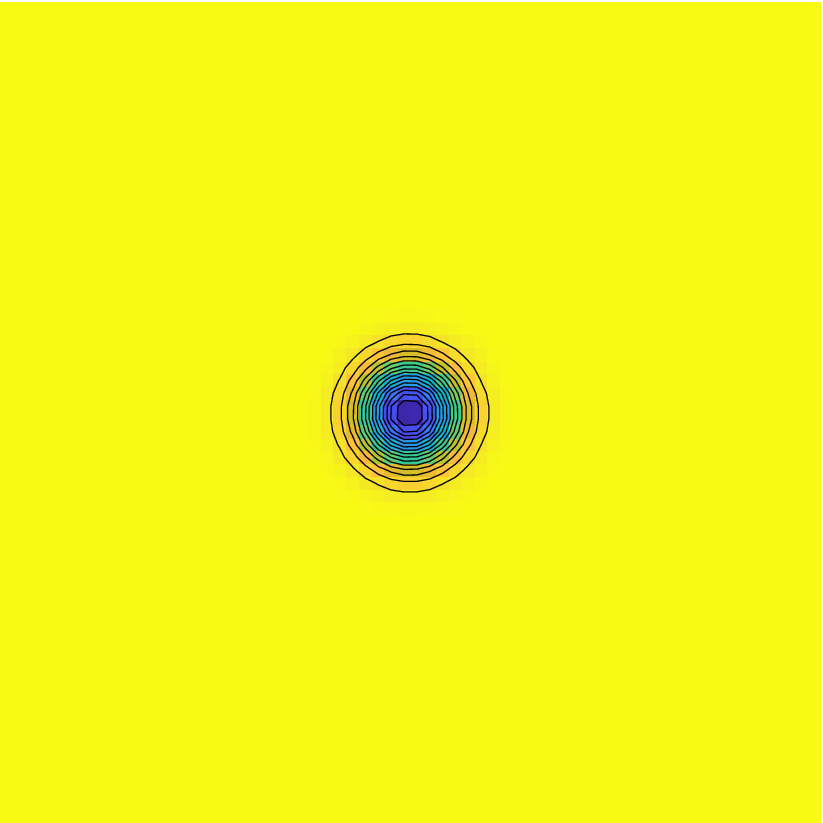}
		\includegraphics[width=0.25\textwidth]{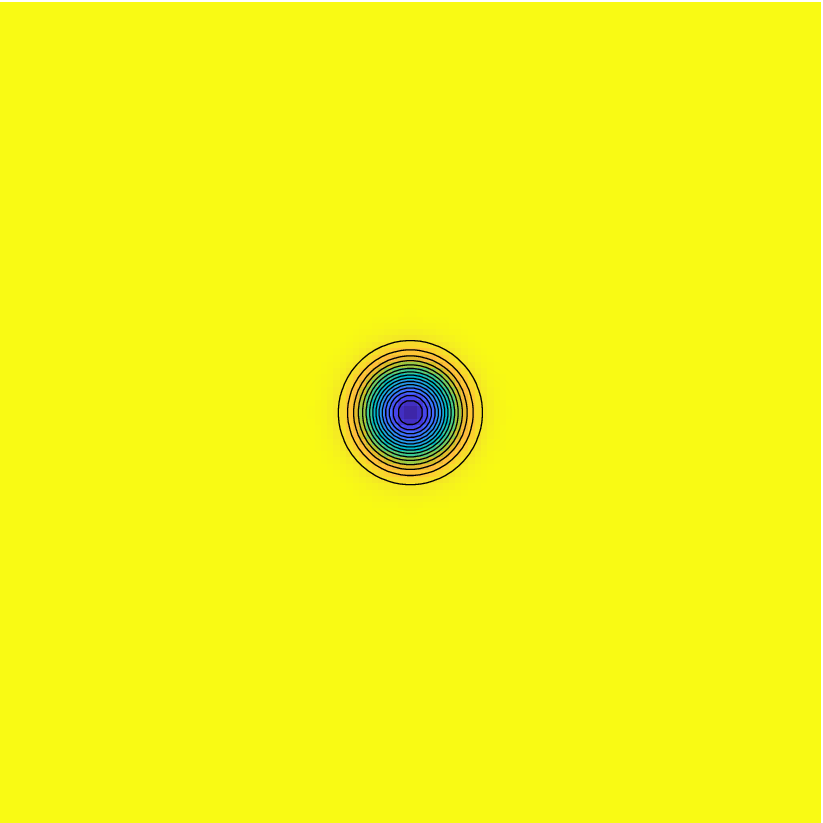}}

	\caption{\label{fig:ShuVortex} Approximation of the Shu-vortex problem at $t = 60$ on grids with $64^2$ (left, middle left, middle right) and $128^2$ (right) cells using AFEG2 (left, middle left) and AFEG2$_{2.0}$ (middle right, right). We performed the middle-left computation with point value limiting from \cite{article:CHP2026}, while we used no limiting for the other computations. }
\end{figure}

In our numerical experiments, we observe the same unphysical behavior for all evolution operators considered in this work within the AF method with the continuous AF reconstruction, namely EG2, EG$^{\mathrm{quad}}$, EG2$_\delta$, EG2$_{\delta,\nu}$, and EG2$_{\delta,\nu}^{\omega}$ with $\delta,\nu \le 1$. In all cases, the computations eventually become unstable.
Figure~\ref{fig:ShuVortex} (left) shows the solution at $t=60$ obtained using AFEG2. The vortex structure is no longer preserved, and spurious artefacts have developed in the solution. Increasing the computational domain does not eliminate this behavior, indicating that it is not caused by boundary effects. Applying the point value limiter does not prevent the instability, as we show in the center-left plot. In contrast to the decoupling near shocks, the problem can not be avoided by choosing different local linearization states. 

\begin{figure}
    \centering
    \includegraphics[width=0.49\linewidth]{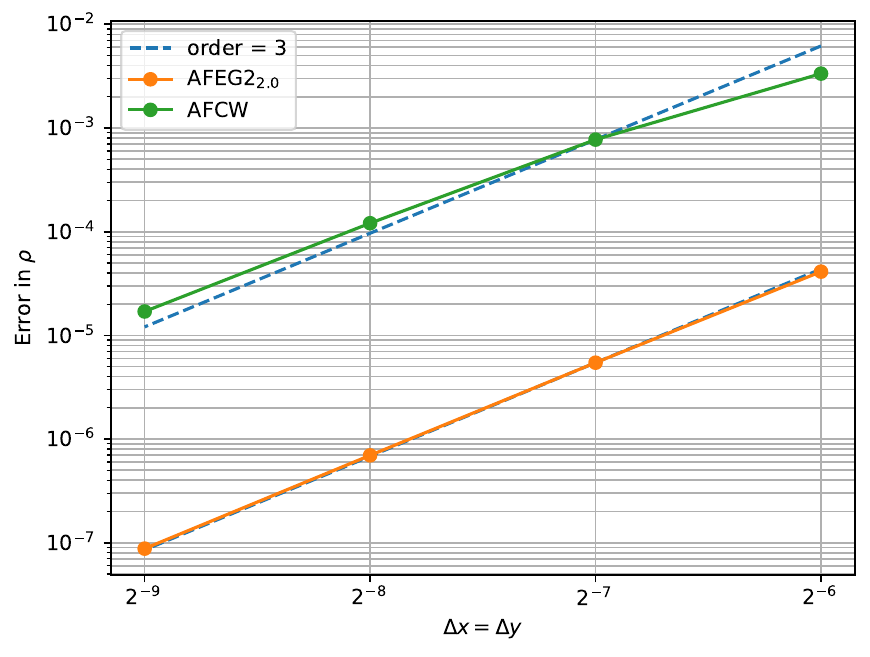}
        \includegraphics[width=0.49\linewidth]{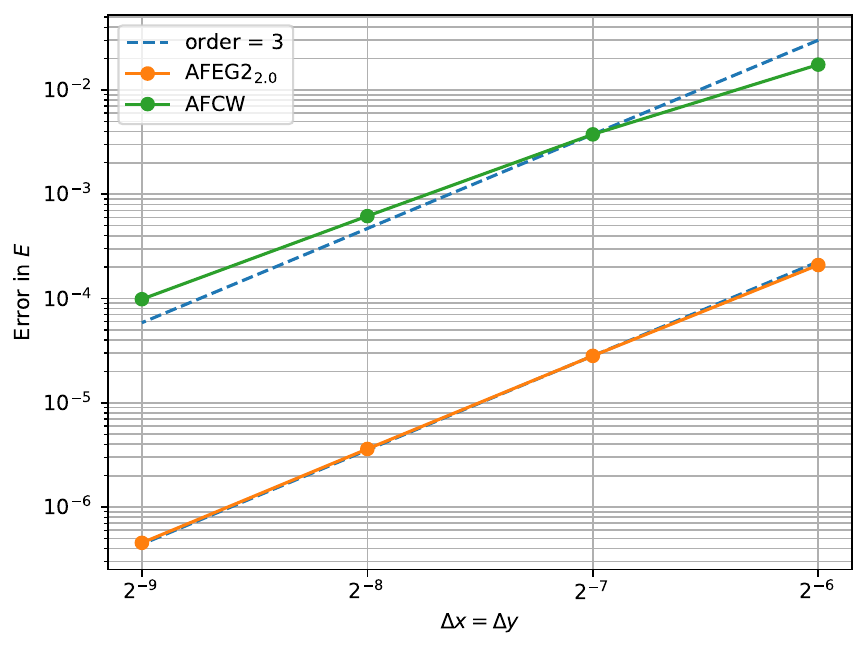}
    \caption{$L_1$-error of the density (left) and total energy (right) for the Shu-vortex problem in Example~\ref{ex:shu_vortex} as a function of the grid resolution at $t=20$, using AFEG2$_{2.0}$ and AFCW. The dashed line indicates third-order convergence.
}
    \label{fig:shu-convergence}
\end{figure}

The connection to transonic flow regimes described above suggests that the integration circles used by the evolution operator do not capture sufficient information from the relevant domain of dependence. This motivates enlarging the integration radii by choosing parameters $\delta, \nu$ greater than one.
For this test case, choosing $\delta \ge 1.8$ was sufficient to prevent the unphysical phenomenon. To reduce the computational cost, we choose $\delta = 2$, which, for the same reasons as the choice $\delta = 1$ allows for a particularly efficient implementation. In Figure~\ref{fig:ShuVortex} (center right and right), we show the numerical solution of the Shu vortex problem obtained with AFEG2$_{2.0}$ on two different grids. The vortex structure is well preserved, and no spurious artefacts are visible. Due to the construction of the evolution operator, the entropy-fix does not reduce the order of accuracy of the method. As shown in \Cref{fig:shu-convergence} \ignore{Table \ref{Tab:shuvortexEG2delta}}, AFEG2$_{2.0}$ retains the expected third-order convergence.

\ignore{
\begin{table}[!ht]
	\caption{Error  measured in the $L_1$-norm and EOC for Example \ref{ex:shu_vortex} 
		using AFEG2$_{2.0}$ with $\mbox{CFL} \le 0.25$ at $t=20$. We added a correction term to the point-value update to eliminate the linearization error.}
	
	\sisetup{
		scientific-notation = true, 
		round-mode=places
	}
	\begin{center}
		\begin{minipage}{0.8\textwidth}
			\begin{tabular}{
					*1{S[table-column-width=0.7cm,table-text-alignment=center]}
					*4{S[table-column-width=1.5cm,table-text-alignment=center]}
					*4{S[table-column-width=0.7cm,table-text-alignment=center, round-precision=4]}
				}
				\toprule
				{Res.}      & \multicolumn{4}{c}{Error
					in }      & \multicolumn{4}{c}{EOC in}  \\
				\midrule
				& {$\bar{\rho}$}   & {$\bar{\rho u}$}  & {$\bar{\rho v}$}        &
				{$\bar{E}$}
				& {$\bar{\rho}$}   & {$\bar{\rho u}$}  & {$\bar{\rho v}$}        &	{$\bar{E}$}\\			 
				\midrule
				{64} & \num{4.11653458e-05} & \num{9.11100803e-05}  & \num{8.50038341e-05} & \num{2.09450041e-04}    & {---} &  {---}   & {---} & {---}\\
				{128} & \num{5.45122192e-06} & \num{1.21054625e-05}  & \num{1.15999800e-05} & \num{2.81902620e-05}    & {2.92} &  {2.91}   & {2.87} & {2.89}\\
				{256} & \num{6.96272269e-07} & \num{1.54470223e-06}  & \num{1.49626906e-06} & \num{3.60635566e-06}    & {2.97} &  {2.97}   & {2.95} & {2.97}\\
				{512} & \num{8.77518045e-08} & \num{1.94228598e-07}  & \num{1.89315045e-07} & \num{4.53786796e-07}    & {2.99} &  {2.99}   & {2.98} & {2.99}\\
				\bottomrule
				
			\end{tabular}
			\label{Tab:shuvortexEG2delta}
		\end{minipage}
	\end{center}
\end{table}}
Furthermore, the numerical results show that modifying only the pressure integration radius is sufficient, i.e., no enlargement of the velocity integration circles is required. For AF methods with the most compact stencil, $\delta$ needs to satisfy $\delta \cdot CFL \le 0.5$. However, we can use $\delta=2$ only locally as an entropy-fix for our newly derived evolution operators that allow larger time steps away from transonic rarefaction waves.

\FloatBarrier

\begin{figure}[htbp]
 \centering
 \makebox[\textwidth][c]{%
  \includegraphics[width=0.30\textwidth]
   {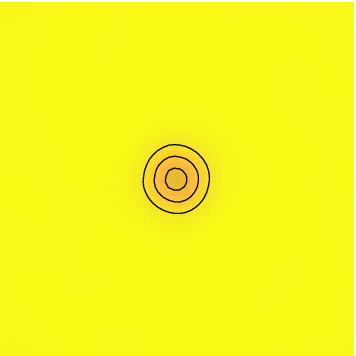}%
  \hspace{0.025\textwidth}%
  \includegraphics[width=0.30\textwidth]
   {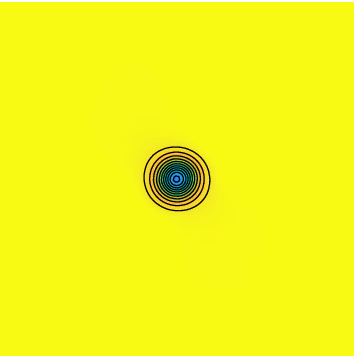}%
  \hspace{0.025\textwidth}%
  \includegraphics[width=0.30\textwidth]
   {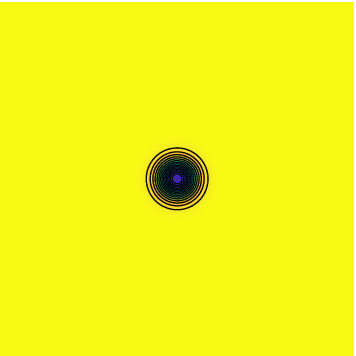}%
 }
 \caption{AFCW approximations of the density field for the isentropic Shu-vortex
 at $t=60$ on
 $64^2$ (left), $128^2$ (middle), and $256^2$ (right) grids, computed with $\mathrm{CFL}\leq0.7$ and without any limiting.
All panels use the same color scale and contour levels.}
 \label{fig:afcw-gauss2-shu-vortex}
\end{figure}

To assess whether the deterioration observed for AFEG2 also occurs
for this fully discrete configuration, We repeat Example~\ref{ex:shu_vortex} using AFCW with
$\mathrm{CFL}\leq0.7$, without additional point-value or flux limiting and no entropy fix is applied in these computations.
Figure~\ref{fig:afcw-gauss2-shu-vortex} shows the numerical density
fields at $t=60$ on different grids. On all three grids, the vortex remains localized and retains
nested, nearly circular contours. No nonlocal spurious structures comparable
to those in the AFEG2 results in Figure~\ref{fig:ShuVortex} are visible.

The vortex becomes progressively sharper and more compact as the grid is
refined from $64^2$ to $256^2$. Because the same color scale and contour levels
are used in all panels, this behavior is consistent with reduced numerical
diffusion, which broadens the vortex on coarse grids, rather than a physical
contraction of the vortex. The $256^2$ solution therefore provides a
better-resolved approximation of the localized vortex profile. These results
demonstrate that the complete AFCW configuration remains stable up to $t=60$. \ignore{Since several components of the discretizations differ from
AFEG2, the comparison does not identify a single component responsible for the
improved long-time behavior.} 
Numerical results presented in Figure~\ref{fig:afcw-gauss2-shu-vortex}  indicate  that the 
numerical viscosity of AFCW provides sufficient stabilization for the Shu-vortex problem.

\ignore{
\begin{table}[!ht]
	\caption{Error measured in the $L_1$-norm and EOC for Example~\ref{ex:shu_vortex} 
		using AFCW with $\mbox{CFL} \le 0.70$ at $t=20$. We added a correction term to the point value update to eliminate the linearization error.}
	
	\sisetup{
		scientific-notation = true, 
		round-mode=places
	}
	
	\begin{center}
		\begin{minipage}{0.8\textwidth}
			\begin{tabular}{
					*1{S[table-column-width=0.7cm,table-text-alignment=center]}
					*4{S[table-column-width=1.5cm,table-text-alignment=center]}
					*4{S[table-column-width=0.7cm,table-text-alignment=center, round-precision=4]}
				}
				\toprule
				{Res.} & \multicolumn{4}{c}{Error in}
				& \multicolumn{4}{c}{EOC in} \\
				\midrule
				& {$\bar{\rho}$} & {$\bar{\rho u}$} & {$\bar{\rho v}$} & {$\bar{E}$}
				& {$\bar{\rho}$} & {$\bar{\rho u}$} & {$\bar{\rho v}$} & {$\bar{E}$} \\
				\midrule
				
				{64}
			& \num{3.34343815e-03}
			& \num{7.16723094e-03}
			& \num{7.21486121e-03}
			& \num{1.75003765e-02}
			& {---} & {---} & {---} & {---} \\
			
			{128}
			& \num{7.72581252e-04}
			& \num{1.42329747e-03}
			& \num{1.36728516e-03}
			& \num{3.75037189e-03}
			& {2.11} & {2.33} & {2.40} & {2.22} \\
			
			{256}
			& \num{1.20693974e-04}
			& \num{2.64241534e-04}
			& \num{2.64205575e-04}
			& \num{6.17096113e-04}
			& {2.68} & {2.43} & {2.37} & {2.60} \\
			
			{512}
			& \num{1.69927075e-05}
			& \num{4.49885272e-05}
			& \num{4.51033938e-05}
			& \num{9.84998819e-05}
			& {2.83} & {2.55} & {2.55} & {2.65} \\
				
				\bottomrule
			\end{tabular}
			\label{tab:afcw-shu-vortex-convergence}
		\end{minipage}
	\end{center}
\end{table}}

The accuracy of AFCW is examined separately at $t=20$.
\ignore{
Table~\ref{tab:afcw-shu-vortex-convergence} shows that, between the $256^2$ and
$512^2$ grids, the density EOC reaches $2.83$, indicating that the density
approximation is approaching its third-order asymptotic regime. Over the same
refinement step, the EOCs for the momentum and total energy components range
from $2.55$ to $2.65$, showing that these components have not yet fully reached
the asymptotic regime on the grids considered.}
Figure~\ref{fig:shu-convergence} illustrates the convergence behavior of AFCW. Between the $256^2$ and $512^2$ grids, the corresponding  experimental order of convergence (EOC) for the density is $2.83$, indicating that the density approximation is approaching its third-order asymptotic regime. For the total energy, the EOC is slightly lower. The momentum components exhibit a similar convergence behavior to the total energy. Over the same refinement step, the EOCs for the momentum and total energy components range from $2.55$ to $2.65$, indicating that these components have not yet fully reached the asymptotic regime on the grids considered.

\FloatBarrier

\section{Further Numerical Results}\label{sec:numerical-results}
In this section, we present numerical results for fully discrete AF and AFCW methods for the Euler equations. We consider both the original EG2 evolution operator and the new evolution operators introduced in \Cref{sec:4}. We also compare numerical results obtained with continuous AF reconstruction with those obtained with discontinuous CWENO reconstruction. First, we perform two further convergence studies to compare the accuracy of the resulting methods. We then consider test problems involving shocks and other interesting flow configurations. All computations are performed using time steps close to the stability limit suggested by the linear stability analysis in \Cref{sec:lin-stability}. The numerical examples in Sections~\ref{sec:RP}, \ref{sec:SodShockTube}, and \ref{sec:jet} further demonstrate that the limiting strategy introduced in \cite{article:CHP2026} can be successfully combined with the new evolution operators. This strategy includes bound-preserving limiting of point values, shock-indicator-based limiting to reduce unphysical oscillations, and flux limiting to ensure bound-preserving cell-average values.


\subsection{Investigation of Accuracy}
In \cite{article:CHP2026}, we showed that AFEG2 achieves third-order accuracy for the Euler equations. The convergence studies presented here confirm this result and further demonstrate that the method using the new evolution operators also achieves the expected third-order accuracy. 
We first consider the advection of a smooth density pulse.

\begin{example}[Density advection]\label{ex:advection}
	\label{ex:density_advection}
	We consider a smooth density advection problem. The initial condition is given by
	\begin{align*}
		\rho(x,y,0) &= 2.5 \exp\left(-40\left((x-x_a)^2 + (y-y_a)^2\right)\right) + 0.1, \\
		u(x,y,0)&=v(x,y,0)=w_0,\qquad p(x,y,0)=1,\quad w_0\in\{0.5,1.0\},
	\end{align*}
	where the initial pulse is centered at $x_a = y_a = -0.31$. The computational domain is $[-1,1]^2$, with periodic
boundary conditions in both coordinate directions. The exact solution is obtained by translating the initial density profile with the constant velocity $(w_0,w_0)$, i.e., $$\rho(x,y,t) = \rho (x-w_0t,y-w_0t,0),$$
	while the velocity and pressure remain constant.
\end{example}

Depending on the choice of reconstruction, this test case requires the entropy fix proposed in \Cref{sec:efix}. For the continuous AF reconstruction, we therefore use AFEG2$_{2.0}$ with $CFL \le 0.25$. We compare the resulting approximations with those obtained using AFCW with $CFL \le 0.7$.

\begin{figure}
    \centering
    \includegraphics[width=0.49\linewidth]{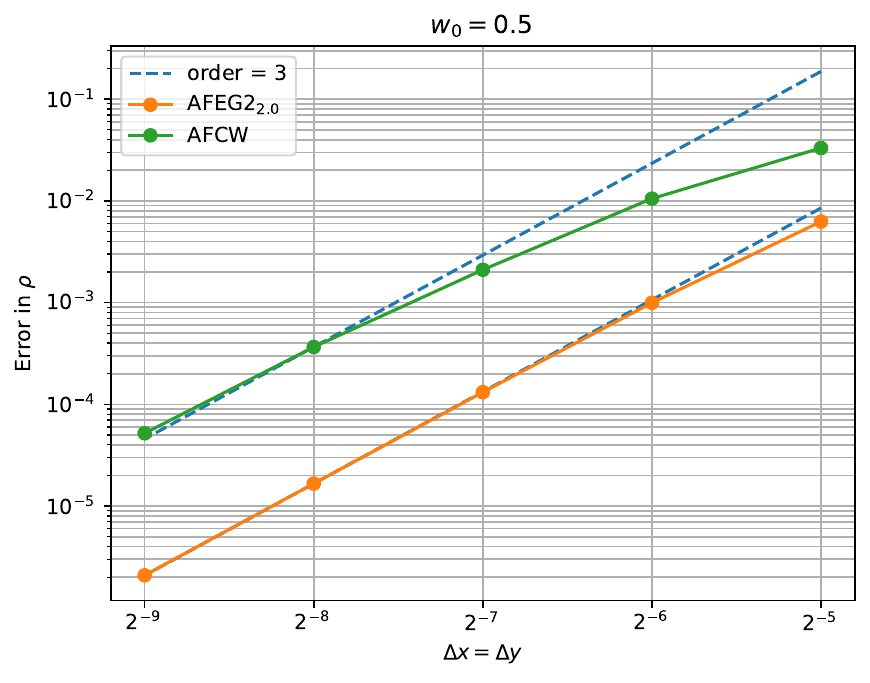}
        \includegraphics[width=0.49\linewidth]{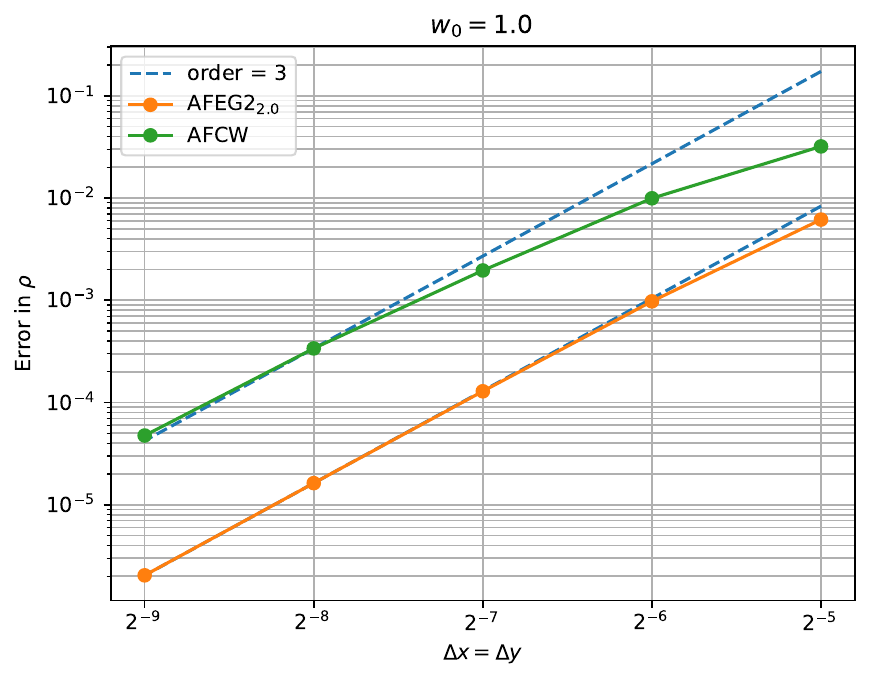}
    \caption{$L_1$-error for the density advection problem in Example~\ref{ex:density_advection} as a function of the grid resolution for $w_0=0.5$ at $t=4$ (left) and $w_0=1.0$ at $t=2$ (right), using AFEG2$_{2.0}$ and AFCW. The dashed line indicates third-order convergence.
}
    \label{fig:density-advection}
\end{figure}

In \Cref{fig:density-advection}, we plot the $L_1$-error of the density against the grid resolution for $w_0=0.5$ at time $t=4$ (left) and for $w_0=1$ at time $t=2$ (right), corresponding to one complete rotation in both cases. 
Thus, the error is computed by comparing the numerical solution with the exact solution.
A reference line corresponding to third-order convergence is included in both plots.

For both methods, we observe similar convergence behavior for the two choices of $w_0$. AFEG2$_{2.0}$ shows the expected third-order convergence. For AFCW, the observed convergence order is slightly below three on the considered grids, but the results indicate that the method approaches third-order convergence in the asymptotic regime. The errors obtained with AFCW are approximately one order of magnitude larger than those obtained with AFEG2$_{2.0}$.

If the convergence study is performed using a method based on the continuous AF reconstruction without applying the entropy fix, such as AFEG2, no apparent problems are initially observed in the convergence behavior of any of the solution components. Only for $w_0=1$, a slight reduction in the observed order becomes visible on the $512^2$ grid. However, examining the pressure point values reveals the underlying issue already on much coarser grids. In the region where transonic-rarefaction-like configurations arise, the pressure point values do not remain constant, although they should. When the entropy fix is applied, the pressure point values remain constant up to machine precision. In contrast, no entropy fix is required when using the AFCW reconstruction.

As a last convergence study, we consider the moving vortex problem proposed by Kadioglu et al. \cite{article:KRM2008}.
	\begin{example}[Moving Vortex]\label{ex:MovVortEuler}
		We consider a smooth traveling vortex for the Euler equations. A rotating vortex is initially located at $(0.5,0.5)$ and propagates with constant speed $(u_c,v_c)$. The initial values have the form 
		\begin{align*}
			\rho(x,y,0) &=
			\begin{cases}
				\rho_c + \dfrac{1}{2}\left(1 - r^2\right)^6, & r < 1, \\
				\rho_c, & \text{otherwise},
			\end{cases} \\[1ex]
			u(x,y,0) &=
			\begin{cases}
				u_c - 1024 \sin(\theta)\,(1 - r)^6 r^6, & r < 1, \\
				u_c, & \text{otherwise},
			\end{cases} \\[1ex]
			v(x,y,0) &=
			\begin{cases}
				v_c + 1024 \cos(\theta)\,(1 - r)^6 r^6, & r < 1, \\
				v_c, & \text{otherwise},
			\end{cases} \\[1ex]
			p(x,y,0) &=
			\begin{cases}
				p_c + \bigl(p(r) - p(1)\bigr), & r < 1, \\
				p_c, & \text{otherwise}.
			\end{cases}
		\end{align*}
		Here $r$ is the scaled distance from the initial center of the vortex, i.e.
		\begin{align*}
			r = \sqrt{(x-0.5)^2 + (y-0.5)^2}/R,
		\end{align*}
         $\theta$ denotes the polar angle about the
initial vortex center $(0.5,0.5)$.
		The function $p(r )$ is a polynomial  of degree 36. In our computation we use $R = 0.4, \rho_c = 0.5, u_c = v_c = 1$ and $p_c = 0.1. $ We simulate the vortex on the domain $\left[0, 1\right]^2$ using periodic boundary conditions. 
	\end{example}

\begin{figure}
    \centering
    \includegraphics[width=0.49\linewidth]{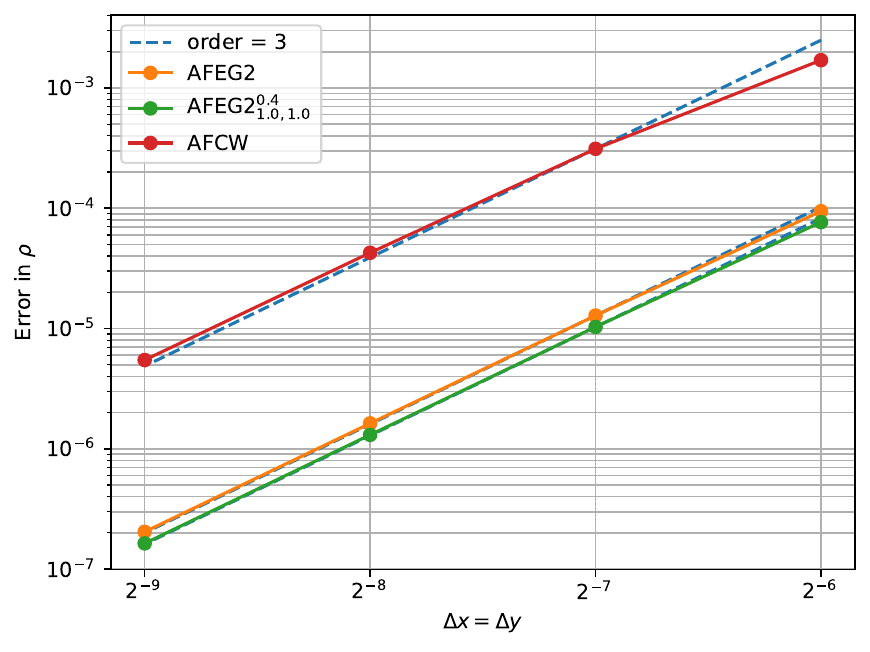}
        \includegraphics[width=0.49\linewidth]{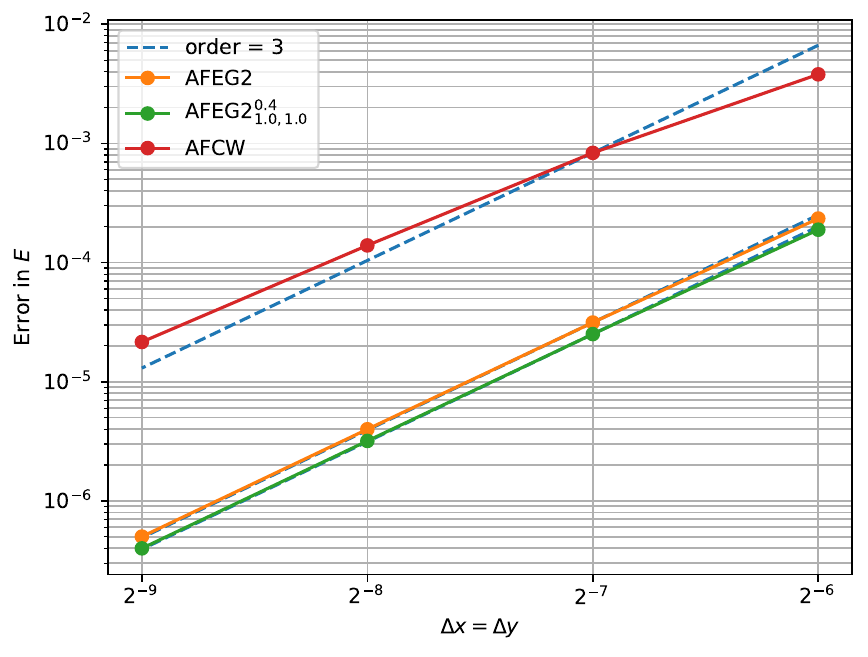}
    \caption{$L_1$-error of the density (left) and total energy (right) for Example~\ref{ex:MovVortEuler} as a function of the grid resolution at $t=1$ , using AFEG2, AFEG2$_{1.0,1.0}^{0.4}$, and AFCW. The dashed lines indicate third-order convergence.
}
    \label{fig:mov-vortex-convergence}
\end{figure}

\Cref{fig:mov-vortex-convergence} shows the $L_1$-errors of the density and total energy plotted against the grid resolution at time $t=1$ for AFEG2, AFEG2$_{1.0,1.0}^{0.4}$, and AFCW. For both AFEG2 and AFEG2$_{1.0,1.0}^{0.4}$, the results confirm the expected third-order convergence. For AFCW, the density approximation also approaches the third-order asymptotic regime, with an EOC of $2.96$ between the $256^2$ and $512^2$ grids. Over the same refinement step, the EOCs for the momentum and total-energy components range from $2.68$ to $2.69$, indicating that these components are still in the pre-asymptotic regime on the grids considered. Similar to the previous test case, the errors obtained with AFCW are more than an order of magnitude larger than those obtained with AFEG2 and AFEG2$_{1.0,1.0}^{0.4}$.

\ignore{
	\begin{table}[!ht]
		\caption{Error  measured in the $L_1$-norm and EOC for Example \ref{ex:MovVortEuler} 
			using {\color{afcwpurple}AFEG2} with $\mbox{CFL} \le 0.279$. We added a correction term to the point-value update to eliminate the linearization error.}
		
		\sisetup{
			scientific-notation = true, 
			round-mode=places
		}
		\begin{center}
			\begin{minipage}{0.8\textwidth}
				\begin{tabular}{
						*1{S[table-column-width=0.7cm,table-text-alignment=center]}
						*4{S[table-column-width=1.5cm,table-text-alignment=center]}
						*4{S[table-column-width=0.7cm,table-text-alignment=center, round-precision=4]}
					}
					\toprule
					{Res.}      & \multicolumn{4}{c}{Error
						in }      & \multicolumn{4}{c}{EOC in}  \\
					\midrule
					& {$\bar{\rho}$}   & {$\bar{\rho u}$}  & {$\bar{\rho v}$}        &
					{$\bar{E}$}
					& {$\bar{\rho}$}   & {$\bar{\rho u}$}  & {$\bar{\rho v}$}        &	{$\bar{E}$}\\			 
					\midrule
					{64} & \num{9.43943157e-05} & \num{1.60600269e-04}  & \num{1.63710053e-04} & \num{2.34253665e-04}    & {---} &  {---}   & {---} & {---}\\
					{128} & \num{1.28056525e-05} & \num{2.15412502e-05}  & \num{2.21045884e-05} & \num{3.13976525e-05}    & {2.88} &  {2.90}   & {2.89} & {2.90}\\
					{256} & \num{1.62647545e-06} & \num{2.73376687e-06}  & \num{2.81378995e-06} & \num{3.99066618e-06}    & {2.98} &  {2.98}   & {2.97} & {2.98}\\
					{512} & \num{2.03890945e-07} & \num{3.42907710e-07}  & \num{3.53619307e-07} & \num{5.01228872e-07}    & {3.00} &  {2.99}   & {2.99} & {2.99}\\
					\bottomrule
					
				\end{tabular}
				\label{Tab:movingvortexEG2delta}
			\end{minipage}
			\end{center}
		\end{table}}
		
\ignore{
	\begin{table}[!ht]
		\caption{Error  measured in the $L_1$-norm and EOC for Example \ref{ex:MovVortEuler} 
			using {\color{blue}AFEG2$_{1.0,1.0}^{0.4}$ }with $\mbox{CFL} \le 0.444$. We added a correction term to the point-value update to eliminate the linearization error.}
		
		\sisetup{
			scientific-notation = true, 
			round-mode=places
		}
		\begin{center}
			\begin{minipage}{0.8\textwidth}
				\begin{tabular}{
						*1{S[table-column-width=0.7cm,table-text-alignment=center]}
						*4{S[table-column-width=1.5cm,table-text-alignment=center]}
						*4{S[table-column-width=0.7cm,table-text-alignment=center, round-precision=4]}
					}
					\toprule
					{Res.}      & \multicolumn{4}{c}{Error
						in }      & \multicolumn{4}{c}{EOC in}  \\
					\midrule
					& {$\bar{\rho}$}   & {$\bar{\rho u}$}  & {$\bar{\rho v}$}        &
					{$\bar{E}$}
					& {$\bar{\rho}$}   & {$\bar{\rho u}$}  & {$\bar{\rho v}$}        &	{$\bar{E}$}\\			 
					\midrule
					{64} & \num{7.66105193e-05} & \num{1.29157458e-04}  & \num{1.33779527e-04} & \num{1.88802961e-04}    & {---} &  {---}   & {---} & {---}\\
					{128} & \num{1.03008789e-05} & \num{1.72360207e-05}  & \num{1.78896999e-05} & \num{2.51080472e-05}    & {2.89} &  {2.91}   & {2.90} & {2.91}\\
					{256} & \num{1.30586624e-06} & \num{ 2.18591663e-06}  & \num{2.27073389e-06} & \num{3.18513355e-06}    & {2.98} &  {2.98}   & {2.98} & {2.98}\\
					{512} & \num{1.63458886e-07} & \num{2.73929997e-07}  & \num{2.84776738e-07} & \num{3.99390772e-07}    & {3.00} &  {3.00}   & {3.00} & {3.00}\\
					\bottomrule
					
				\end{tabular}
				\label{Tab:movingvortexEG2deltanunew}
			\end{minipage}
			\end{center}
		\end{table}}

\ignore{
\begin{table}[!ht]
	\caption{Error measured in the $L_1$-norm and EOC for Example \ref{ex:MovVortEuler} 
		using AFCW with $\mbox{CFL} \le 0.70$. We added a correction term to the point-value update to eliminate the linearization error.}
	
	\sisetup{
		scientific-notation = true, 
		round-mode=places
	}
	
	\begin{center}
		\begin{minipage}{0.8\textwidth}
			\begin{tabular}{
					*1{S[table-column-width=0.7cm,table-text-alignment=center]}
					*4{S[table-column-width=1.5cm,table-text-alignment=center]}
					*4{S[table-column-width=0.7cm,table-text-alignment=center, round-precision=4]}
				}
				\toprule
				{Res.} & \multicolumn{4}{c}{Error in} 
				& \multicolumn{4}{c}{EOC in} \\
				\midrule
				& {$\bar{\rho}$} & {$\bar{\rho u}$} & {$\bar{\rho v}$} & {$\bar{E}$}
				& {$\bar{\rho}$} & {$\bar{\rho u}$} & {$\bar{\rho v}$} & {$\bar{E}$} \\
				\midrule
				
				{64} 
				& \num{1.70121070e-03} 
				& \num{2.61363280e-03} 
				& \num{2.59111853e-03} 
				& \num{3.80177303e-03}
				& {---} & {---} & {---} & {---} \\
				
				{128} 
				& \num{3.11672517e-04} 
				& \num{5.97601586e-04} 
				& \num{6.00286513e-04} 
				& \num{8.33092856e-04}
				& {2.45} & {2.13} & {2.11} & {2.19} \\
				
				{256} 
				& \num{4.25950872e-05} 
				& \num{1.02348157e-04} 
				& \num{1.03191769e-04} 
				& \num{1.39348766e-04}
				& {2.87} & {2.55} & {2.54} & {2.58} \\
				
				{512} 
				& \num{5.46645101e-06} 
				& \num{1.59365142e-05} 
				& \num{1.60028906e-05} 
				& \num{2.15559954e-05}
				& {2.96} & {2.68} & {2.69} & {2.69} \\
				
				\bottomrule
			\end{tabular}
			\label{tab:afcw-gauss2-moving-vortex}
		\end{minipage}
	\end{center}
\end{table}}

	\ignore{
		Next, we will consider a Gresho Vortex Problem 
		
		\begin{example}[Gresho Vortex]\label{ex:Gresho}
			The Gresho problem is a stationary rotating vortex problem. 
			The velocity field coincides with that of Example \ref{ex:StatVort} 
			The pressure is prescribed such that the centrifugal force is balanced by the pressure gradient and is given by
			\begin{equation}\label{eq:gresho_pressure}
				p(r) =
				\begin{cases}
					5 + \dfrac{25}{2} r^{2}, & 0 \le r < 0.2, \\[0.5ex]
					9 - 4 \ln 0.2 + \dfrac{25}{2} r^{2} - 20 r + 4 \ln r, & 0.2 \le r < 0.4, \\[0.5ex]
					3 + 4 \ln 2, & r \ge 0.4 .
				\end{cases}
			\end{equation}
			The density is set to one. The computational domain is $[-0.75,0.75]^2$, with periodic boundary conditions.
	\end{example}}

\FloatBarrier
	
	\subsection{Two-dimensional Riemann problems}\label{sec:RP}
	\ignore{In the last convergence studies presented above, both methods achieve the expected third-order accuracy and exhibit only minor differences in numerical errors. In this section, we therefore compare the methods for more challenging test cases.}
    We next consider three classical two-dimensional Riemann
problems described in \cite{article:RCG1993,lax1998solution}
to compare the methods for solutions involving discontinuities.
	\begin{example}
		We consider classical two-dimensional Riemann
		problems.
		In the computational domain $[0,1]^2$, the initial values in
		primitive variables $\mathbf{u}=(\rho,u,v,p)$ are given by
		\begin{equation*}
			\mathbf{u}(x,y,0)=
			\begin{cases}
				\mathbf{u}_1, & \quad x > x_0,\ y > y_0,\\
				\mathbf{u}_2, & \quad x \le x_0,\ y > y_0,\\
				\mathbf{u}_3, & \quad x \le x_0,\ y \le y_0,\\
				\mathbf{u}_4, & \quad x > x_0,\ y \le y_0,
			\end{cases}
		\end{equation*}
		where $\mathbf{u}_1,\dots,\mathbf{u}_4$ denote
		constant states
		prescribed in the four quadrants.
		The initial condition is therefore composed of four
		constant states
		separated by discontinuities aligned with the
		coordinate directions,
		intersecting at $(x_0,y_0)$. We use outflow boundary
		conditions.
		
		We consider the following three configurations from
		\cite{lax1998solution}:

		\begin{alignat*}{2}
			&\text{Configuration 3:}\qquad x_0=y_0=0.8, \\
			&\mathbf{u}_1 = (1.5,0,0,1.5),      &\quad
			&\mathbf{u}_2 = (0.5323,1.206,0,0.3),\\
			&\mathbf{u}_3 = (0.138,1.206,1.206,0.029), &\quad
			&\mathbf{u}_4 = (0.5323,0,1.206,0.3).\\[0.5em]
			&\text{Configuration 4:}\qquad x_0=y_0=0.5, \\
			&\mathbf{u}_1 = (1.1,0,0,1.1),   &\quad
			&\mathbf{u}_2 = (0.5065,0.8939,0,0.35),\\
			&\mathbf{u}_3 = (1.1,0.8939,0.8939,1.1),        &\quad
			&\mathbf{u}_4 = (0.5065,0,0.8939,0.35).\\[0.5em]
			&\text{Configuration 12:}\qquad x_0=y_0=0.5, \\
			&\mathbf{u}_1 = (0.5313,0,0,0.4),   &\quad
			&\mathbf{u}_2 = (1,0.7276,0,1),\\
			&\mathbf{u}_3 = (0.8,0,0,1),        &\quad
			&\mathbf{u}_4 = (1,0,0.7276,1).
		\end{alignat*}
		
	\end{example}

\begin{figure}
		
		\makebox[\textwidth][c]{\includegraphics[width=0.3\textwidth]{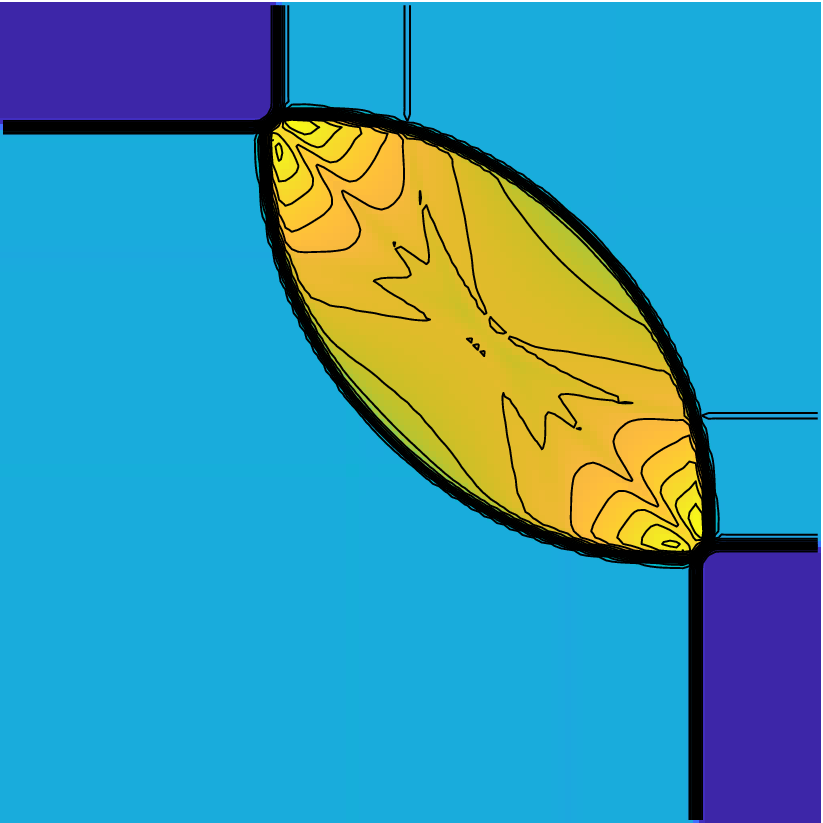}
			\includegraphics[width=0.3\textwidth]{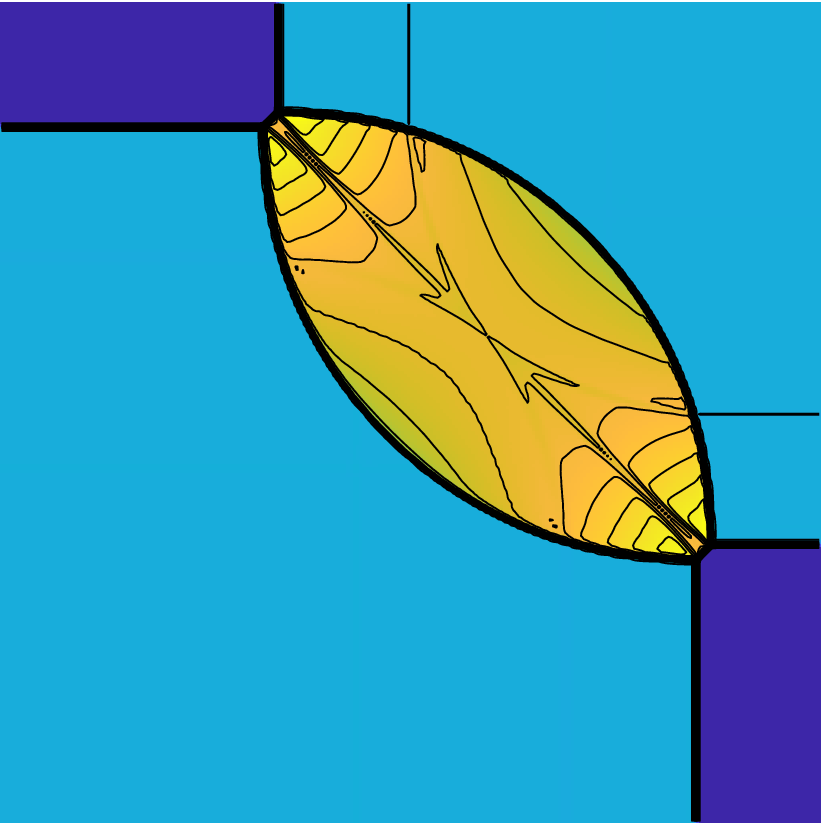}
			\includegraphics[width=0.3\textwidth]{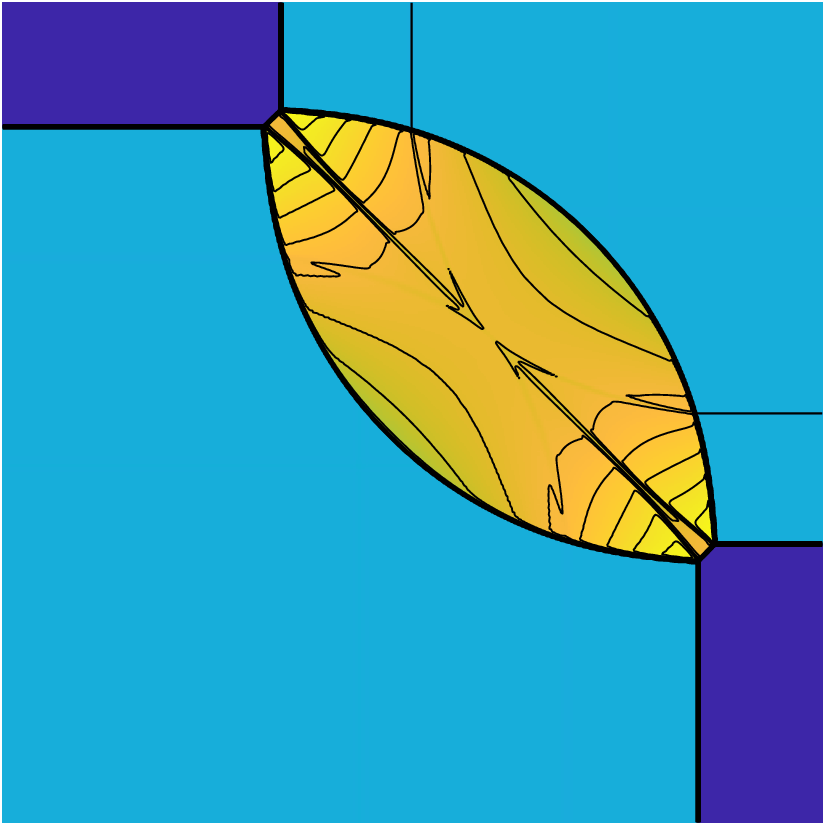}}\\[4pt]
		\makebox[\textwidth][c]{\includegraphics[width=0.3\textwidth]{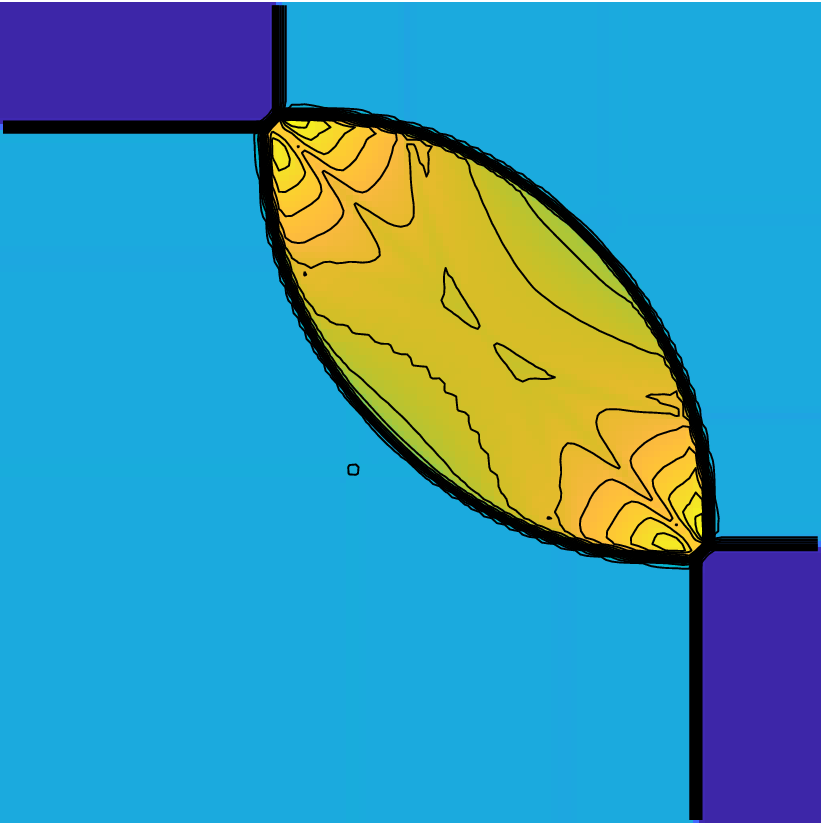}
			\includegraphics[width=0.3\textwidth]{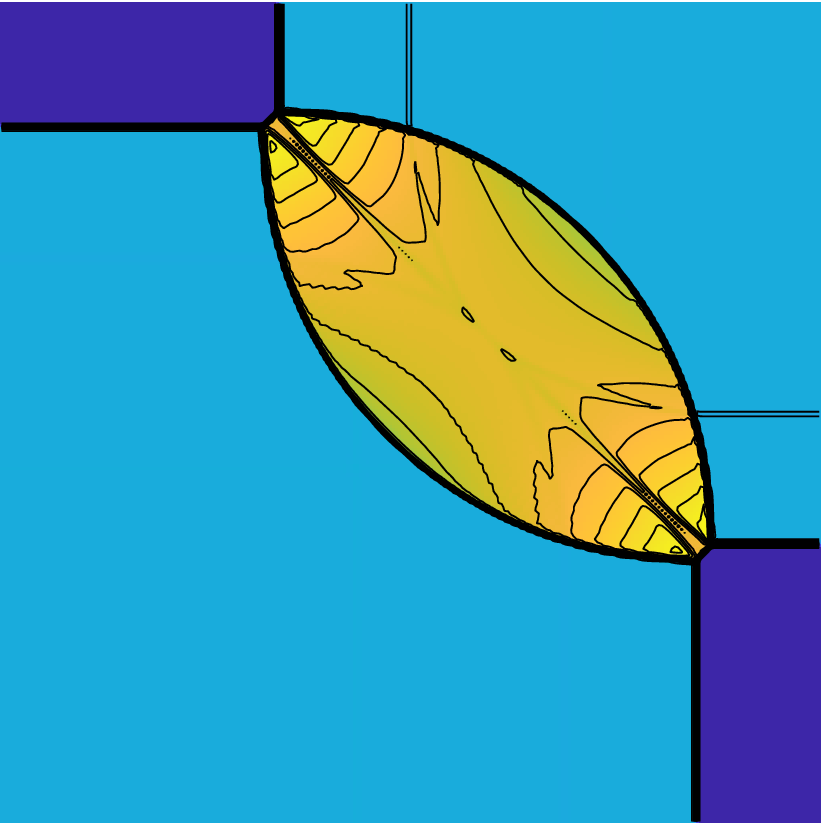}
			\includegraphics[width=0.3\textwidth]{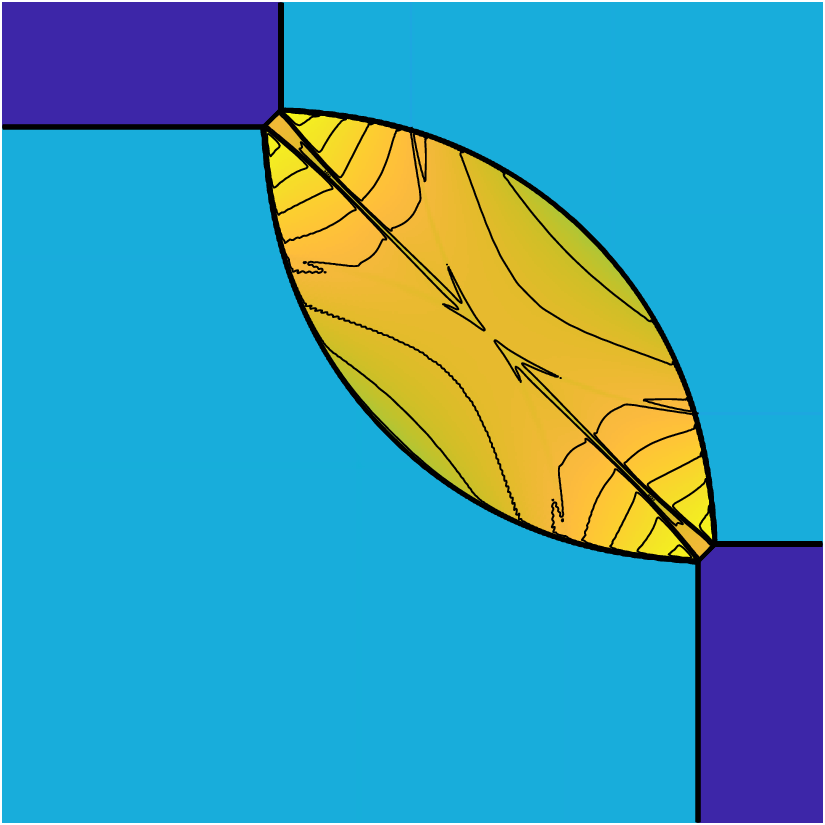}}\\[4pt]
		\makebox[\textwidth][c]{\includegraphics[width=0.3\textwidth]{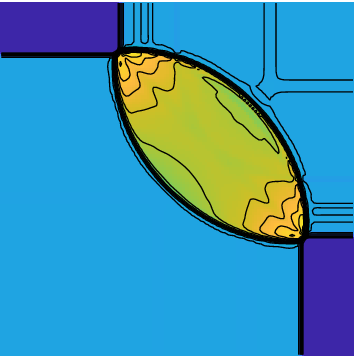}
			\includegraphics[width=0.3\textwidth]{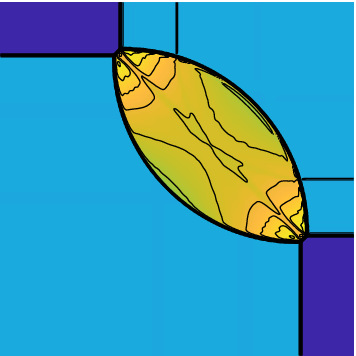}
			\includegraphics[width=0.3\textwidth]{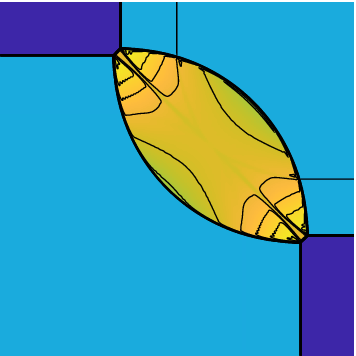}}
		
		\caption{\label{fig:C4} Approximation of a two dimensional Riemann problem (configuration 4) at $t=0.21$ on a grid with $128^2$ (left), $256^2$ (middle), and $512^2$ (right) grid cells with AFEG2 (top row, $CFL \le 0.279$), AFEG2$_{1.0,1.0}^{0.4}$ (middle row, $CFL \le 0.444$) and AFCW (bottom row, $CFL \le 0.7$).}
		\end{figure}

\begin{figure}[htb]
		
		\makebox[\textwidth][c]{\includegraphics[width=0.32\textwidth]{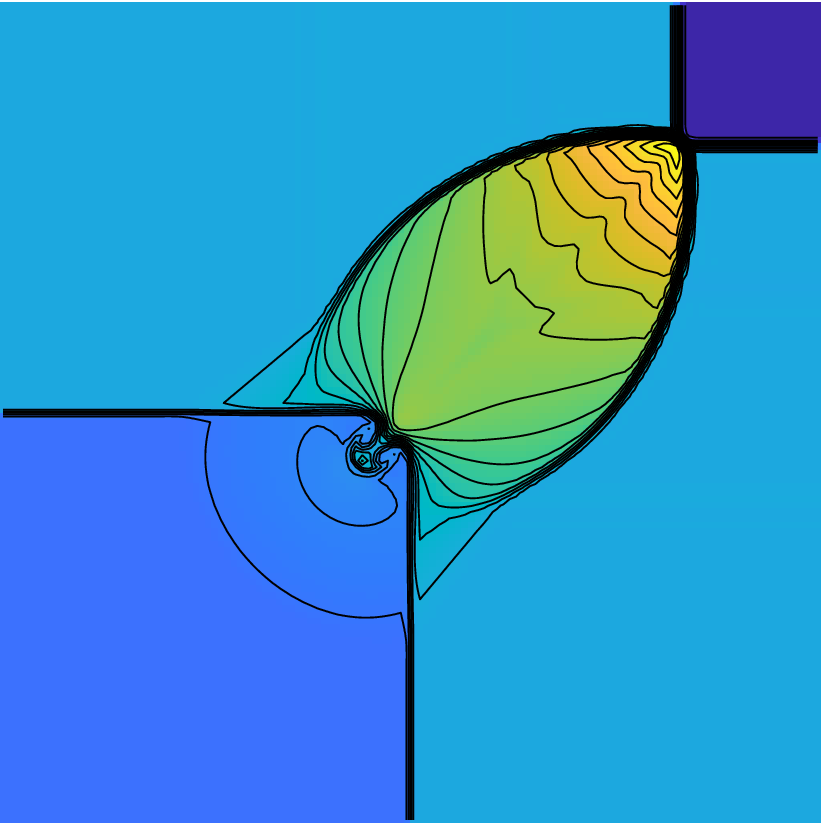}
			\includegraphics[width=0.32\textwidth]{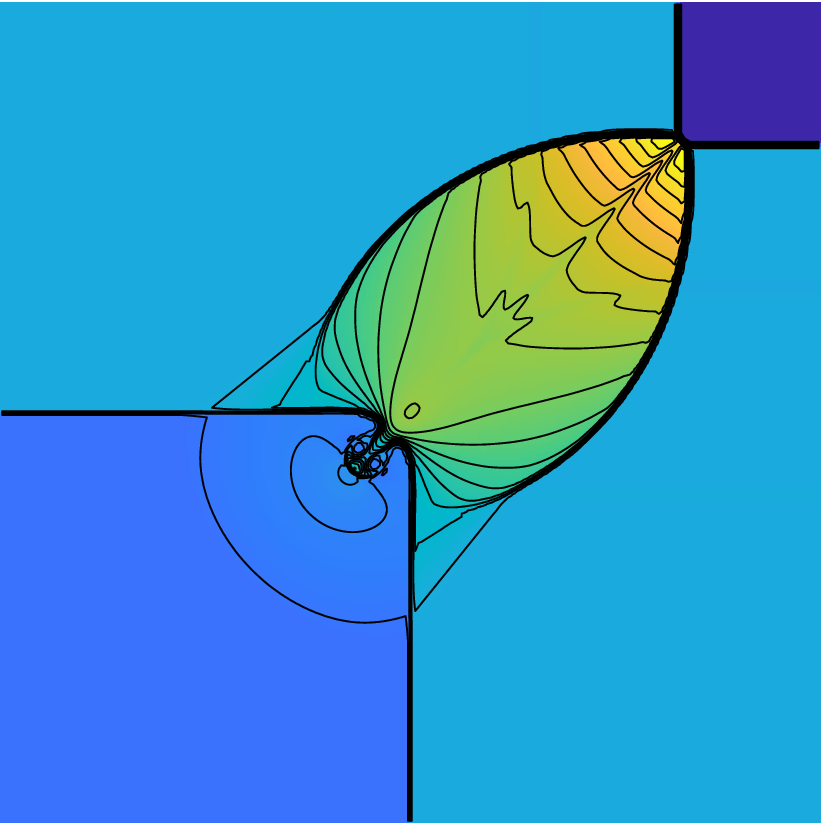}
			\includegraphics[width=0.32\textwidth]{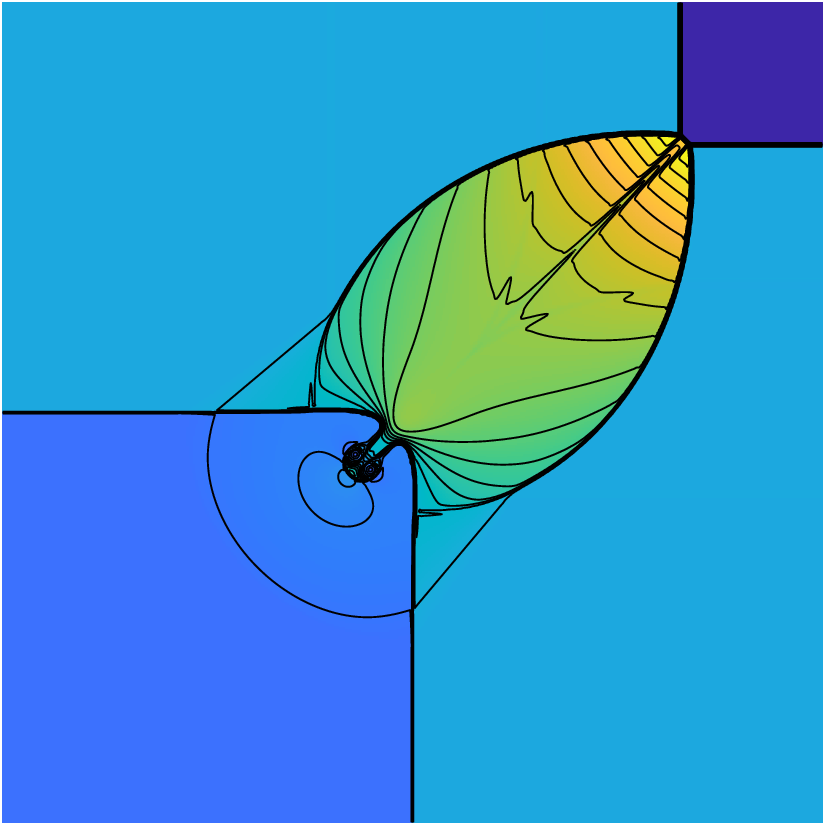}}\\[4pt]
		\makebox[\textwidth][c]{\includegraphics[width=0.32\textwidth]{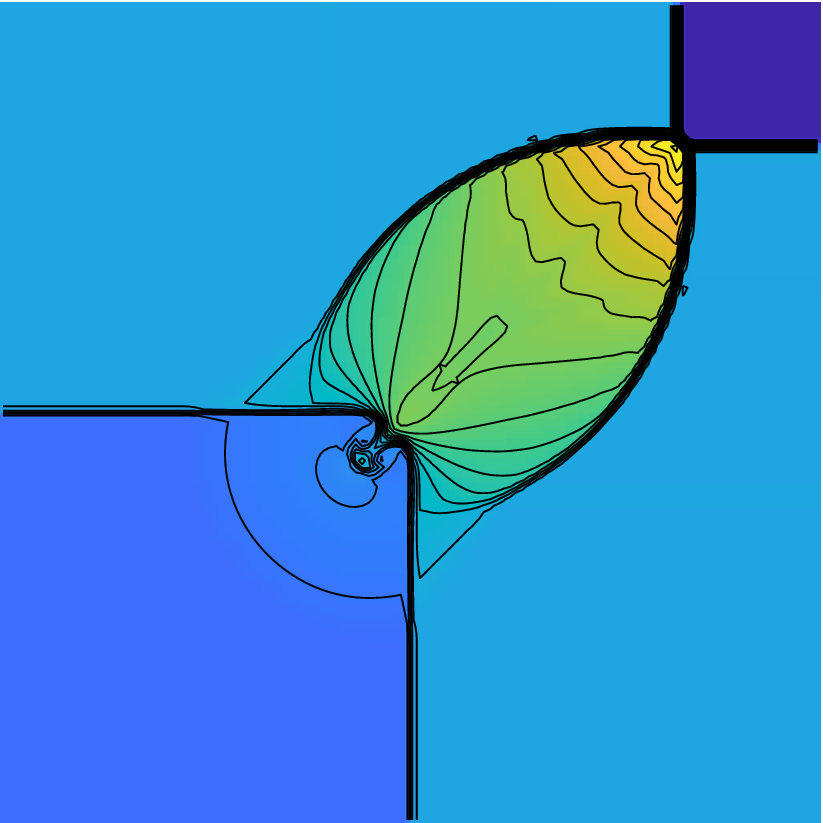}
			\includegraphics[width=0.32\textwidth]{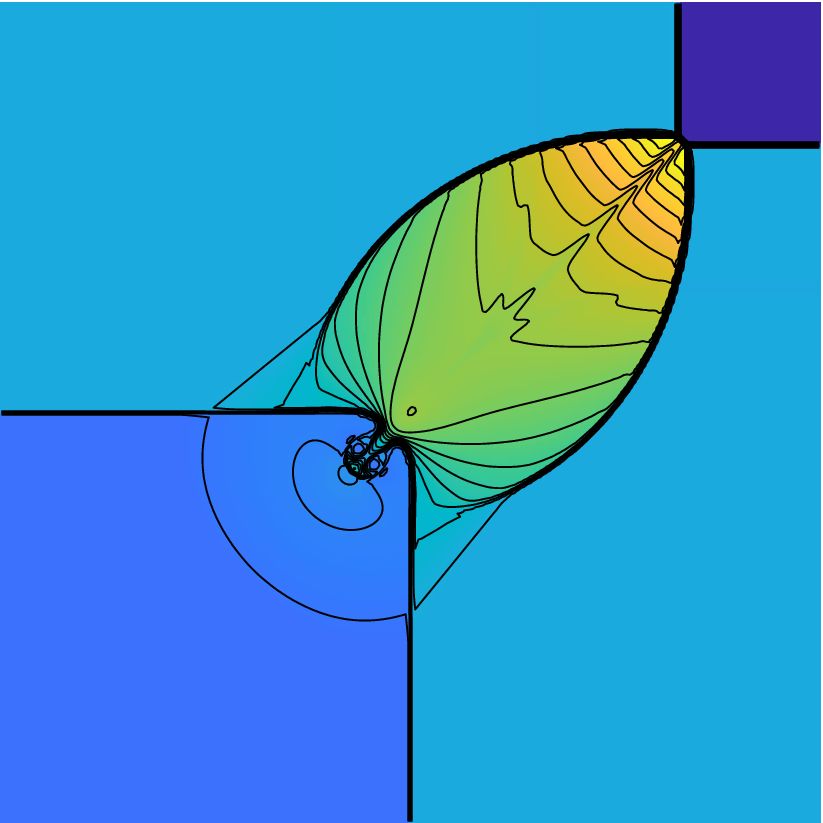}
			\includegraphics[width=0.32\textwidth]{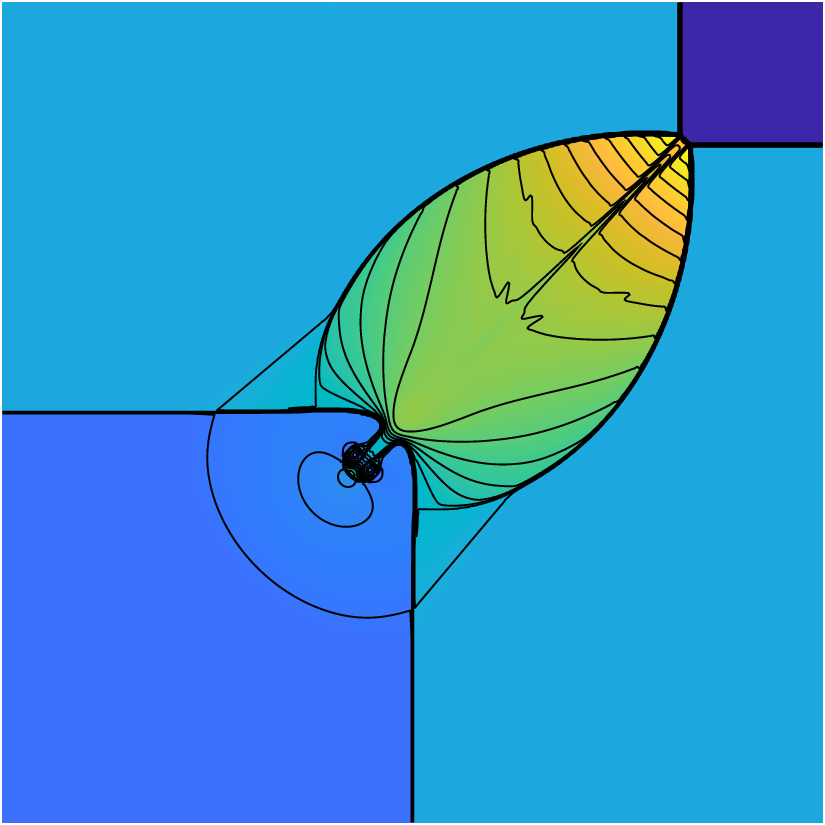}}\\[4pt]
		\makebox[\textwidth][c]{\includegraphics[width=0.32\textwidth]{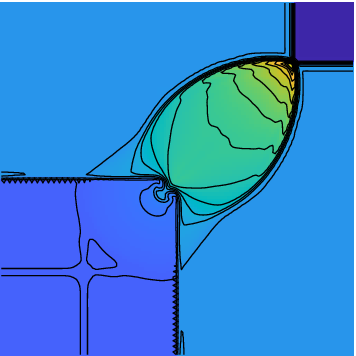}
			\includegraphics[width=0.32\textwidth]{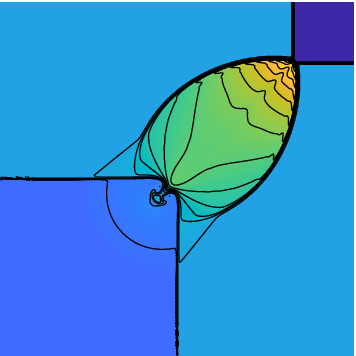}
			\includegraphics[width=0.32\textwidth]{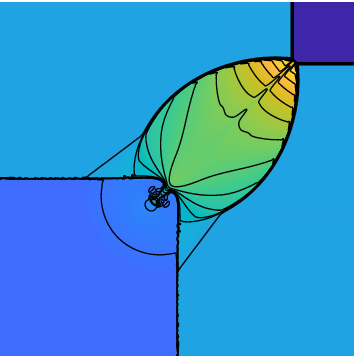}}
		\caption{\label{fig:CF12} Approximation of a two dimensional Riemann problem (configuration 12) at $t=0.21$ on a grid with $128^2$ (left), $256^2$ (middle), and $512^2$ (right) grid cells with AFEG2 (top row, $CFL \le 0.279$), AFEG2$_{1.0,1.0}^{0.4}$ (middle row, $CFL \le 0.444$), and AFCW (bottom row, $CFL \le 0.7$).}
		\end{figure}

    For computations using AFEG2 and AFEG2$_{1.0,1.0}^{0.4}$, the point-value limiter is applied with the shock indicator described above, and both methods are used with CFL numbers close to their respective stability limits. For AFCW, Configurations 4 and 12 are computed with
$\mathrm{CFL}\leq0.7$ without additional point-value
or flux limiting. For Configuration~3, the same shock-indicator-based point-value limiting strategy as in the AFEG2 computations is applied, and the calculation is performed with $\mathrm{CFL}\le 0.5$.
	 Figure~\ref{fig:C4} shows numerical solutions for Configuration 4 obtained with AFEG2, AFEG2$_{1.0,1.0}^{0.4}$ and AFCW on grids with different resolutions. For all grids, the methods provide a good resolution of the solution structures. The AF methods with globally continuous reconstruction provide a slightly sharper approximation of the contact lines forming at the Mach stem. A similar behavior is observed for Configuration~12 shown in Figure \ref{fig:CF12}. 

  \ignore{However, as for Configuration~4, the AFCW solutions exhibit slightly more numerical diffusion, resulting in a somewhat less sharp resolution of the solution structures compared with AFEG2 and AFEG2$_{1.0,1.0}^{0.4}$.}
	
\begin{figure}
		
		\makebox[\textwidth][c]{\includegraphics[width=0.3\textwidth]{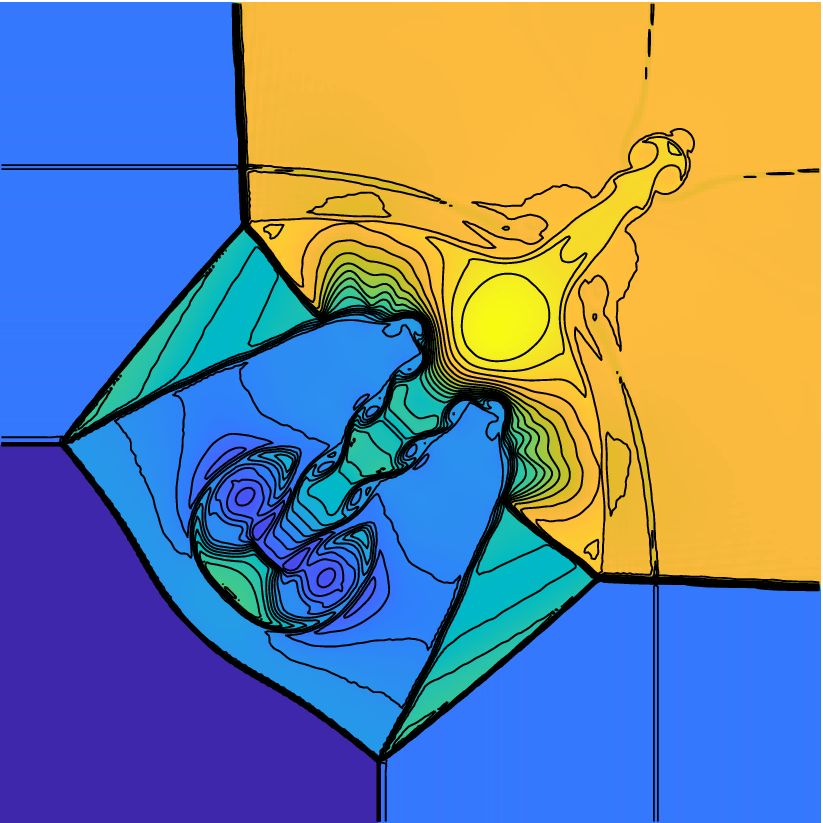}
			\includegraphics[width=0.3\textwidth]{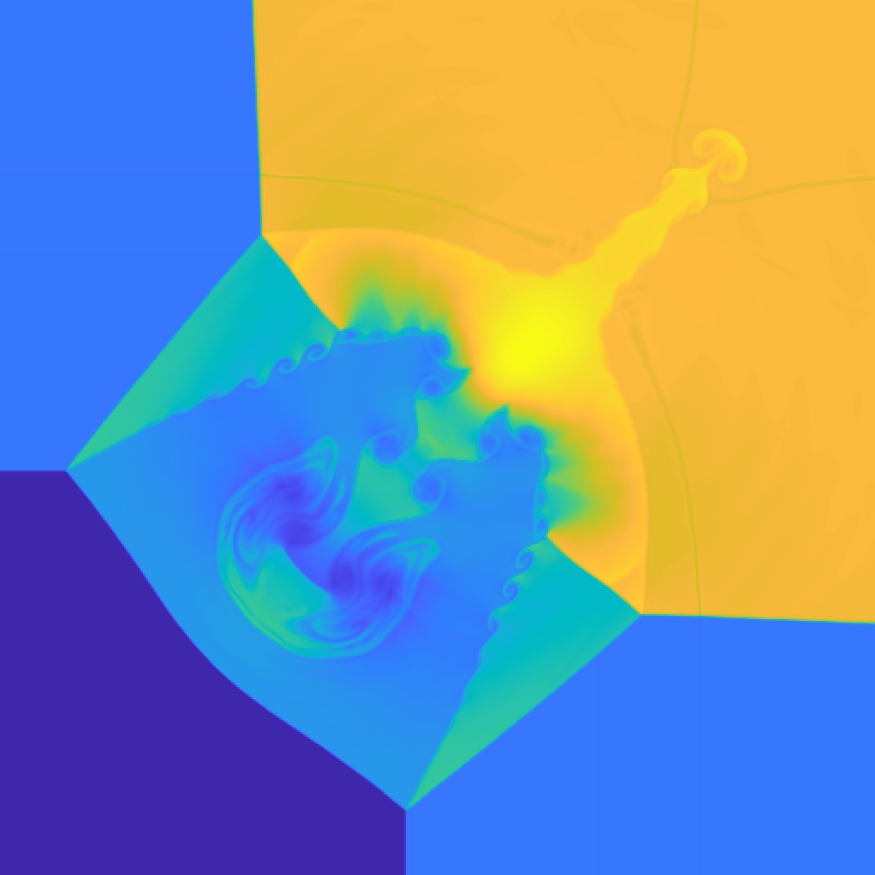} 
			\includegraphics[width=0.3\textwidth]{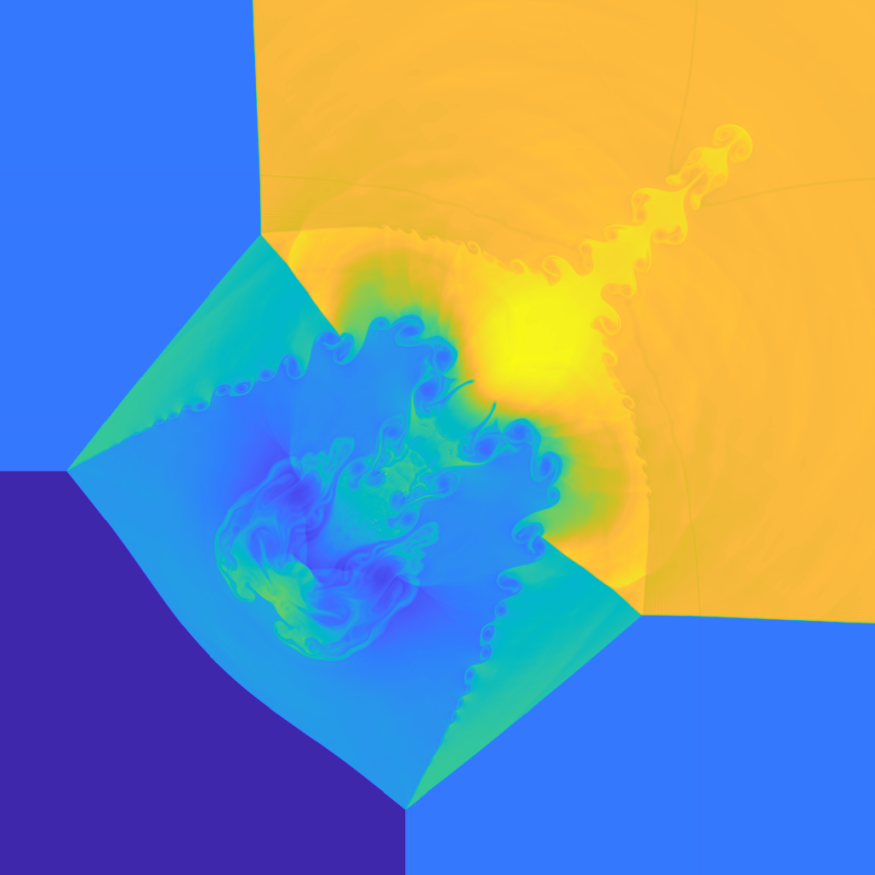}}\\[4pt]
		\makebox[\textwidth][c]{\includegraphics[width=0.3\textwidth]{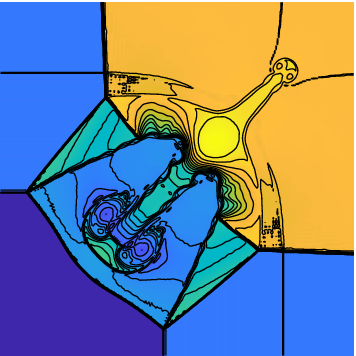}
			\includegraphics[width=0.3\textwidth]{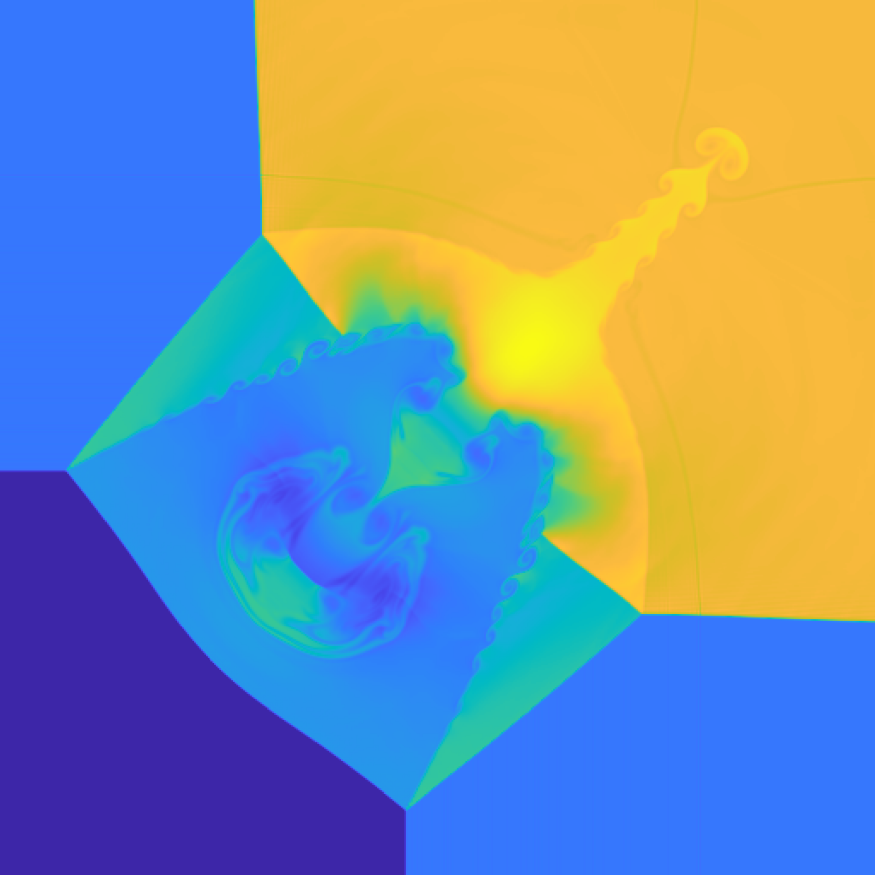}
			\includegraphics[width=0.3\textwidth]{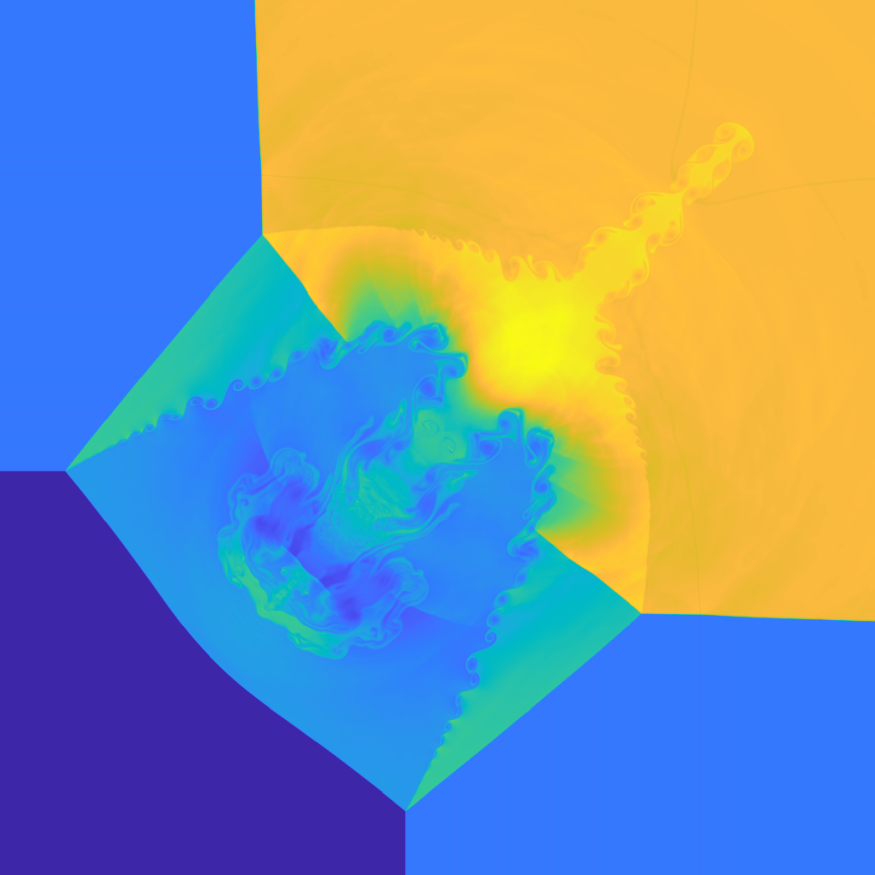}}\\[4pt]
		\makebox[\textwidth][c]{\includegraphics[width=0.3\textwidth]{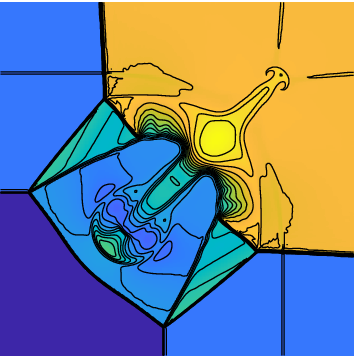}
			\includegraphics[width=0.3\textwidth]{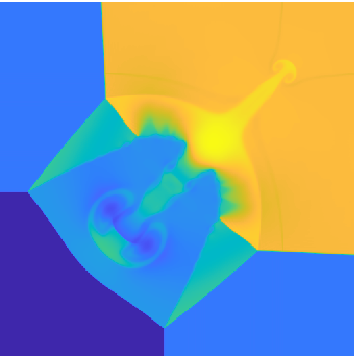}
			\includegraphics[width=0.3\textwidth,trim=70 0 70 0,clip]{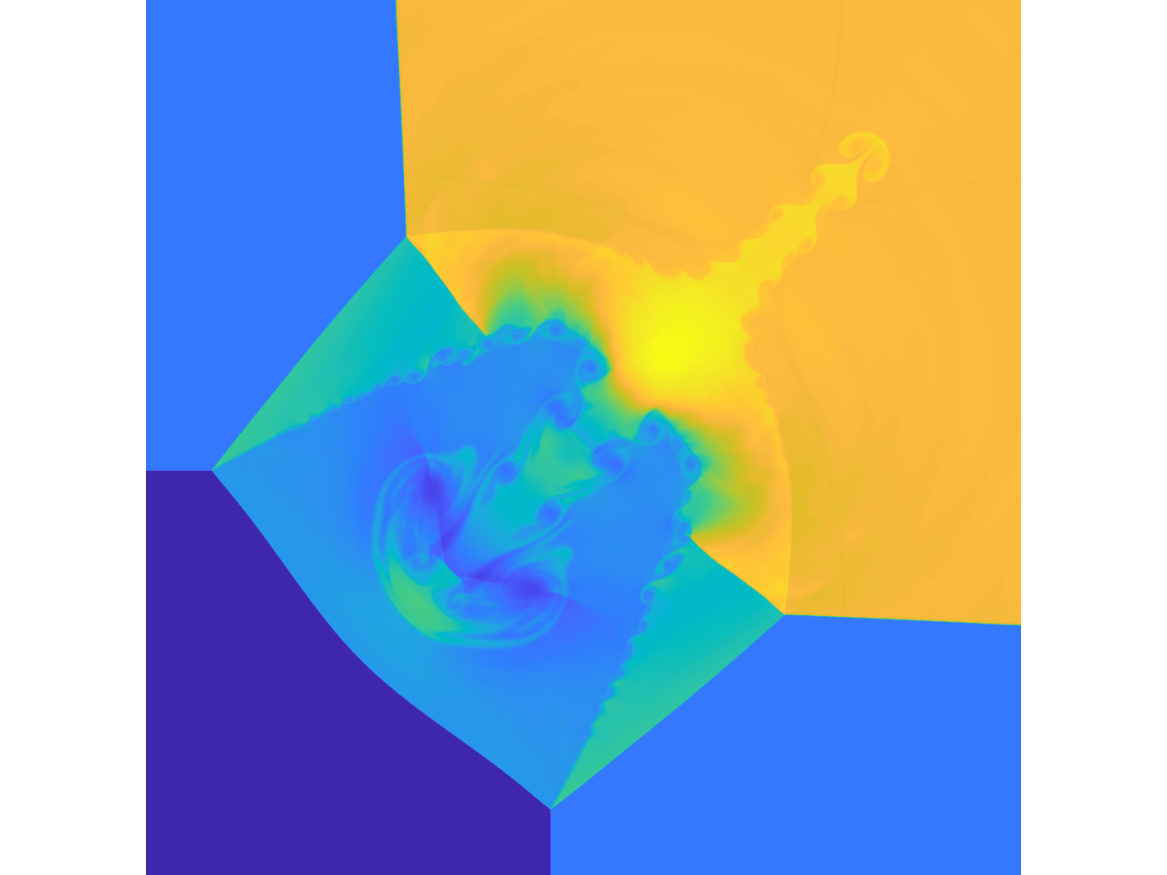}}
		\caption{\label{fig:C3} Approximation of a two dimensional Riemann problem (configuration 3) at $t=0.8$ on a grid with $256^2$ (left), $512^2$ (middle) and $1024^2$ (right) grid cells with AFEG2 (top row, $CFL \le 0.279$), AFEG2$_{1.0,1.0}^{0.4}$ (middle row, $CFL \le 0.444$), and AFCW (bottom row, $\mathrm{CFL}\le 0.5$).}
		\end{figure}

	\ignore{
	In Figure~\ref{fig:C3} we present results of Configuration 3 at time $t=0.8$ calculated with two different methods, AFEG2 and AFEG2$_{1.0,1.0}^{0.4}$ on grids with $256^2$, $512^2$ and $1024^2$ grid cells.  Again, we applied bound-preserving and shock-indicated point-value limiting. Although slight differences in the solution structure are visible, both methods provide comparably accurate approximations.
	\textcolor{afcwpurple}{With the point-value limiter activated, the AFCW
	calculations reproduce the global interacting-wave pattern on all three
	grids. Mesh refinement progressively sharpens the narrow wave structures and
	reveals smaller-scale features in the central interaction region, without
	introducing prominent grid-scale oscillations.}}
In Figure~\ref{fig:C3}, we present the results for Configuration~3 at time $t=0.8$ obtained with AFEG2, AFEG2$_{1.0,1.0}^{0.4}$, and AFCW on grids with $256^2$, $512^2$, and $1024^2$ cells. All three methods reproduce the complex wave pattern well.
    The two AF methods with globally continuous reconstruction show a slightly more pronounced roll-up of the contact lines.
		\subsection{Sod shock tube}\label{sec:SodShockTube}
\begin{figure}[htbp]
		\centering
		\makebox[\textwidth][c]{%
			\includegraphics[width=0.30\textwidth,trim=83.5bp 26bp 69.5bp 41.5bp,clip]{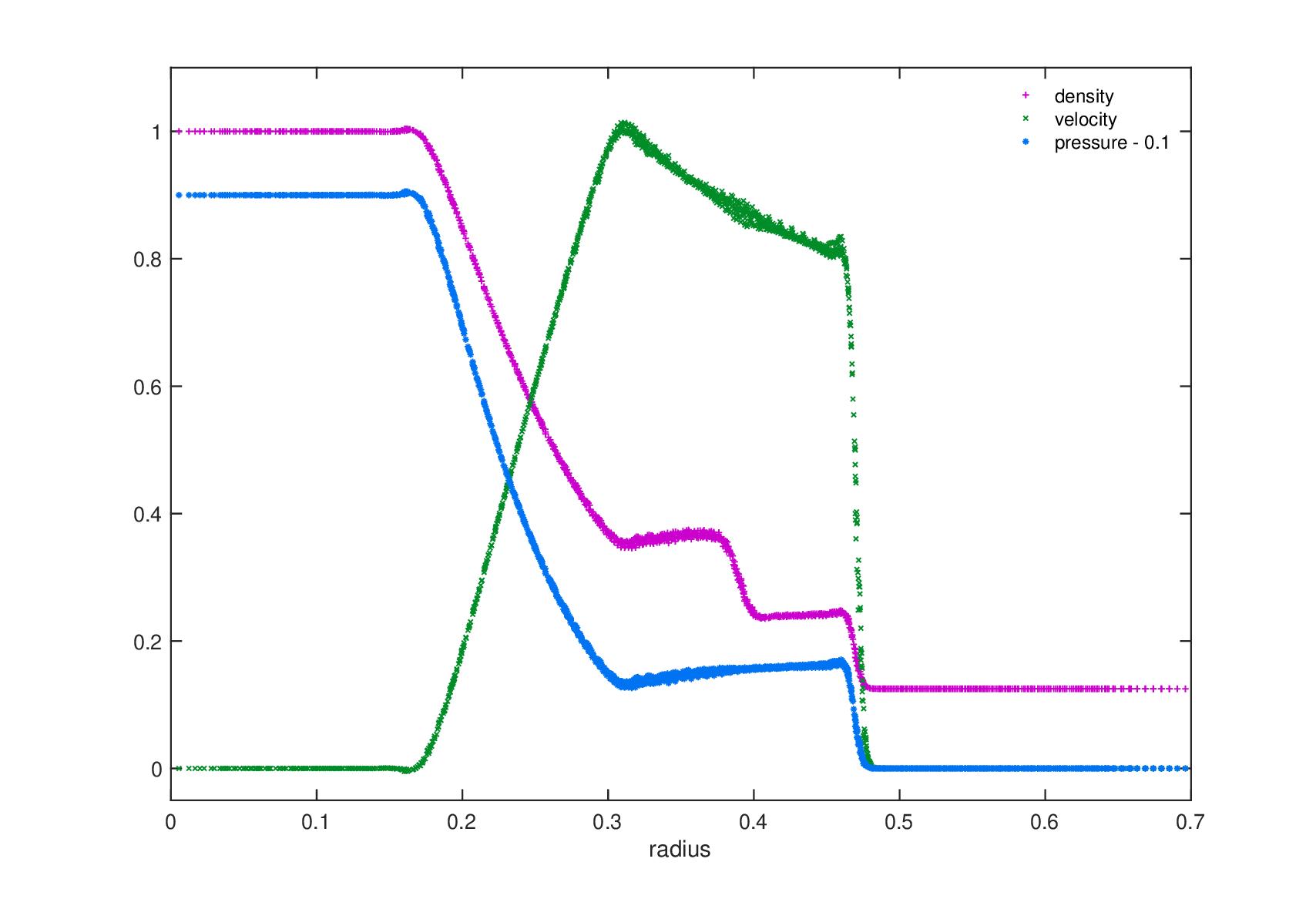}\hfill
			\includegraphics[width=0.30\textwidth,trim=83.5bp 26bp 69.5bp 41.5bp,clip]{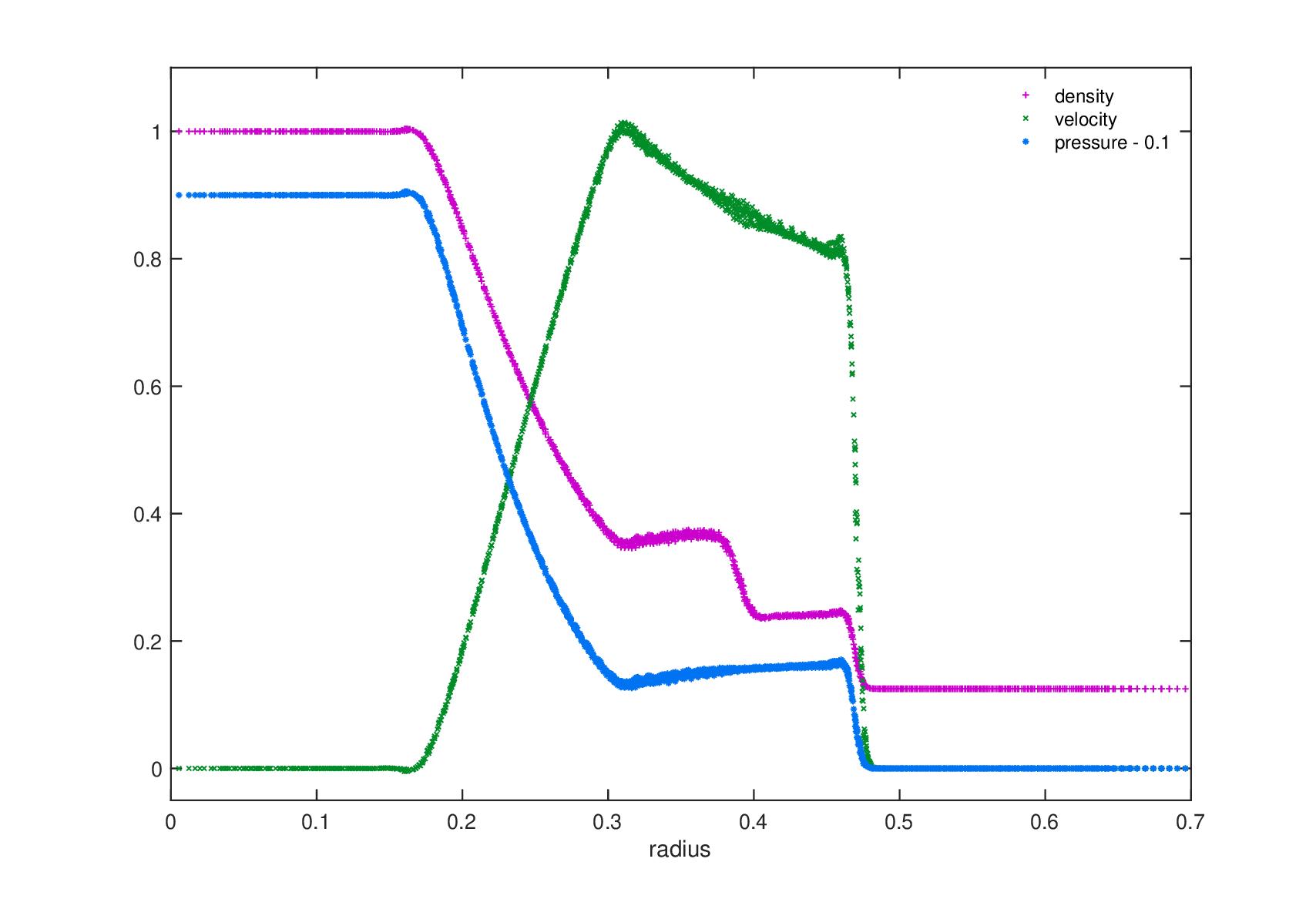}\hfill
			\includegraphics[width=0.30\textwidth]{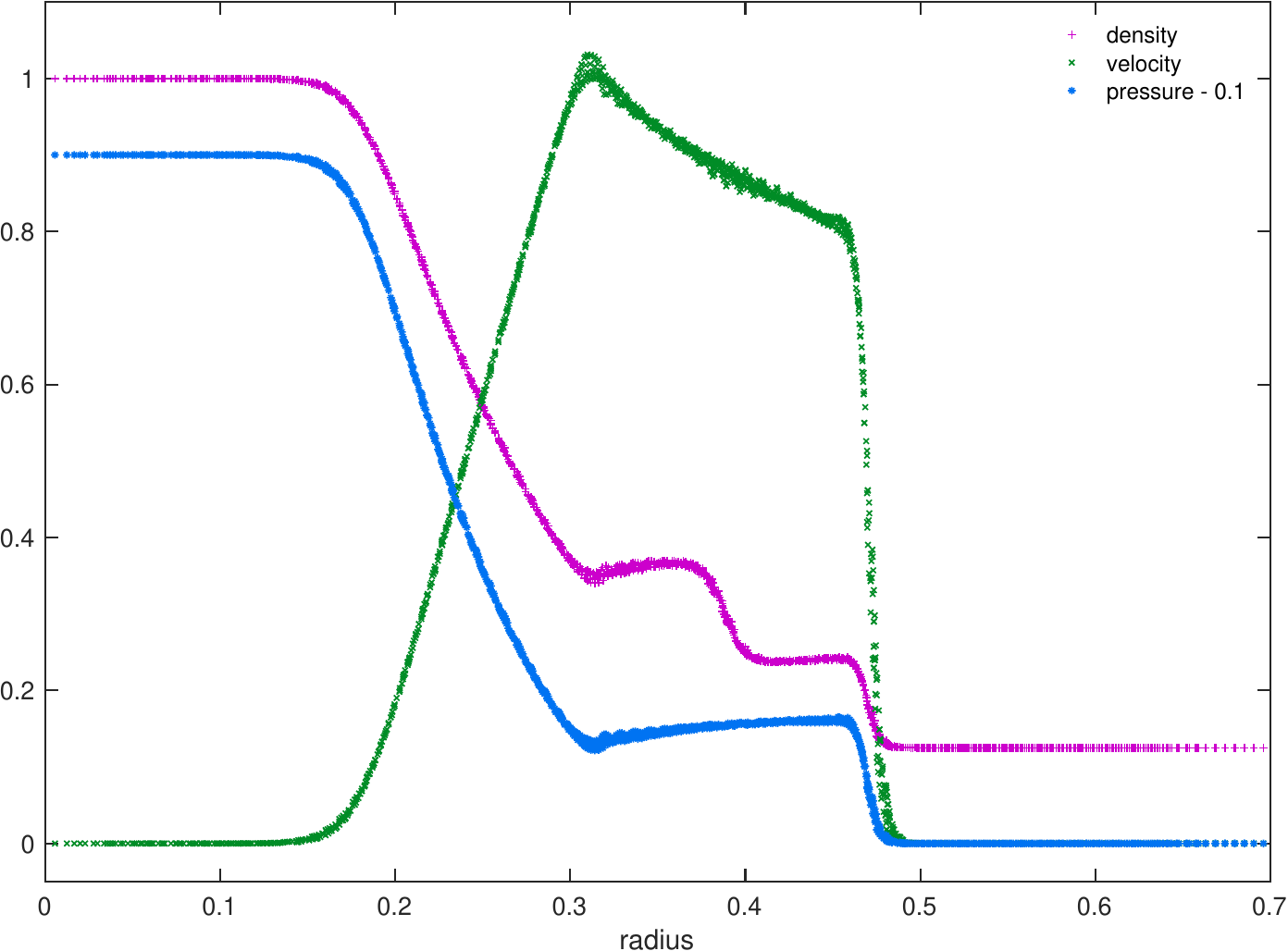}}
		\caption{\label{fig:SOD}
			Radial scatter plot for the two-dimensional Sod shock tube problem at $t=0.1$, computed on a $128^2$ grid using AFEG2 (left, $\mathrm{CFL}\le 0.279$), AFEG2$_{1.0,1.0}^{0.4}$ (middle, $\mathrm{CFL}\le 0.444$), and AFCW (right, $\mathrm{CFL}\le 0.5$). The plotted quantities are density, radial velocity, and shifted pressure $p-0.1$.}
		
		\end{figure}

	As a further test case, we consider the two-dimensional Sod shock tube problem. 
	\begin{example}[Circular Sod shock tube]
		\label{ex:circular_sod}
		We consider the two-dimensional circular Sod shock tube problem on
		$[0,1]^2$. The initial data are
		\[
		(\rho,u,v,p)(x,y,0)=
		\begin{cases}
			(1,0,0,1), & r<0.3,\\
			(0.125,0,0,0.1), & r\ge 0.3,
		\end{cases}
		\]
		where
		\[
		r=\sqrt{\left(x-\frac{1}{2}\right)^2+
			\left(y-\frac{1}{2}\right)^2}.
		\]
		This problem produces an outward-moving circular shock, a contact wave, and an inward-moving rarefaction. 
	\end{example}
	For this test case, we applied the bound-preserving and shock-indicating limiting strategy from \cite{article:CHP2026}. Again, no flux limiting was necessary. Figure~\ref{fig:SOD} shows the solution at $t=0.1$. The three methods produce qualitatively similar radial profiles.
	
		\subsection{Mach 80 Jet}\label{sec:jet}
	We consider the Mach 80 jet problem as an extreme robustness test for the fully discrete AF schemes. For these computations, we use the point-value limiting strategy together with the flux limiter to obtain bound-preserving cell averages.
	
	\begin{example}[Mach 80 jet]
		\label{ex:mach80_jet}

		The computational domain is
		\[
		\Omega=[0,2]\times[-0.5,0.5],
		\]
		and the ratio of specific heats is set to
		\[
		\gamma=\frac{5}{3}.
		\]
		Initially, the domain is filled with a stationary ambient gas,
		\[
		(\rho,u,v,p)(x,y,0)=(0.5,0,0,0.4127).
		\]
		At the left boundary, a high-speed jet is imposed on the interval
		\(|y|<0.05\):
		\[
		(\rho,u,v,p)(0,y,t)=
		\begin{cases}
			(5,30,0,0.4127), & |y|<0.05,\\
			(0.5,0,0,0.4127), & |y|\geq 0.05.
		\end{cases}
		\]
		Outflow boundary conditions are prescribed at the top, bottom, and right
		boundaries. The numerical solution is evolved until \(T=0.07\).
	\end{example}

	\begin{figure}[!htbp]
		\includegraphics[width=0.46\linewidth]{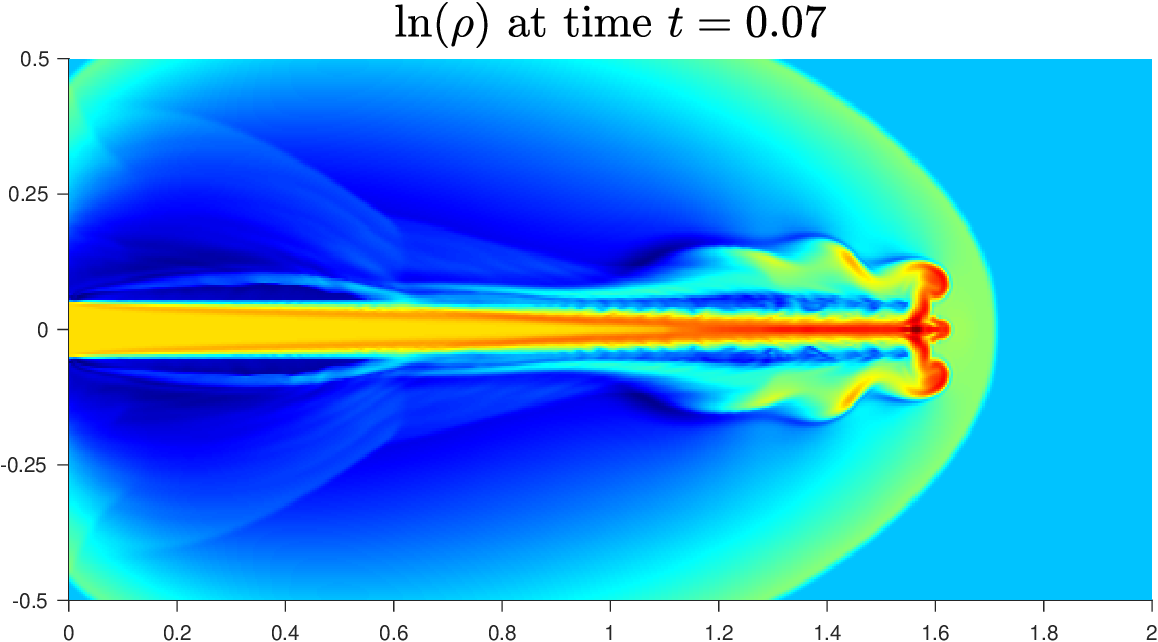}\hspace{-1.8mm} 
		\raisebox{+0.074\height}{
			\includegraphics[width=0.0224\linewidth]{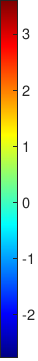}
		}\hspace{-1mm}
		\includegraphics[width=0.46\linewidth]{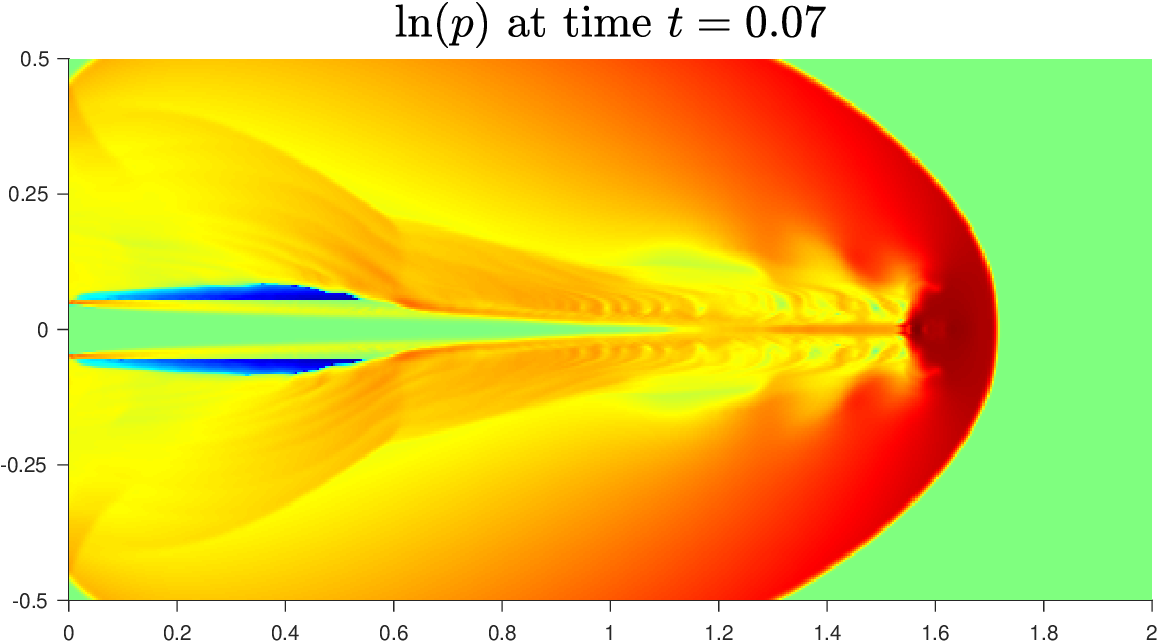}\hspace{-1.8mm} 
		\raisebox{+0.074\height}{
			\includegraphics[width=0.0224\linewidth]{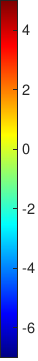}
		}\\
		\includegraphics[width=0.46\linewidth]{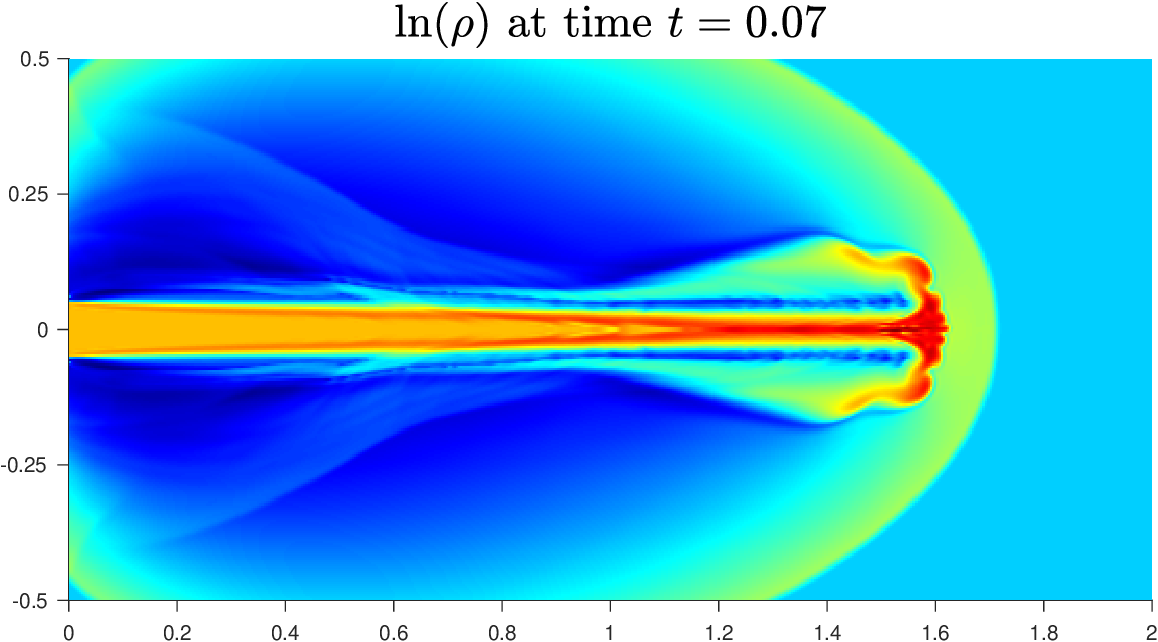}\hspace{-1.8mm} 
		\raisebox{+0.074\height}{
			\includegraphics[width=0.0224\linewidth]{./figs/EG2/Jet/colorbar_logdensity_jet.pdf}
		}\hspace{-1mm}
		\includegraphics[width=0.46\linewidth]{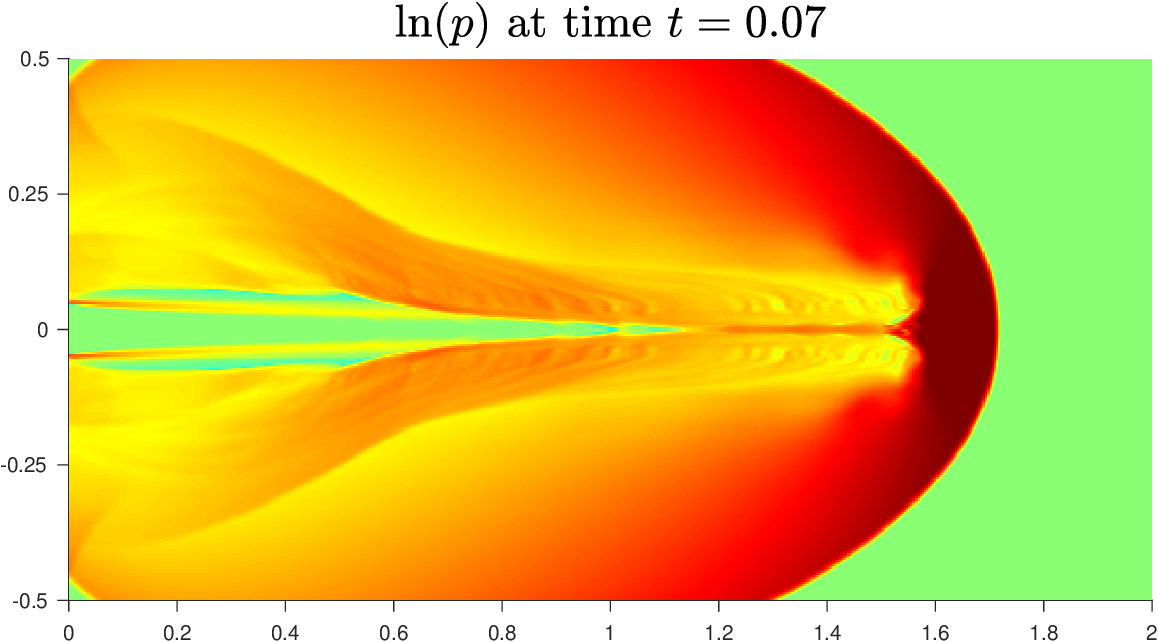}\hspace{-1.8mm} 
		\raisebox{+0.074\height}{
			\includegraphics[width=0.0224\linewidth]{./figs/EG2/Jet/colorbar_logpressure_jet.pdf}
		}\\
        		\includegraphics[width=0.46\linewidth]{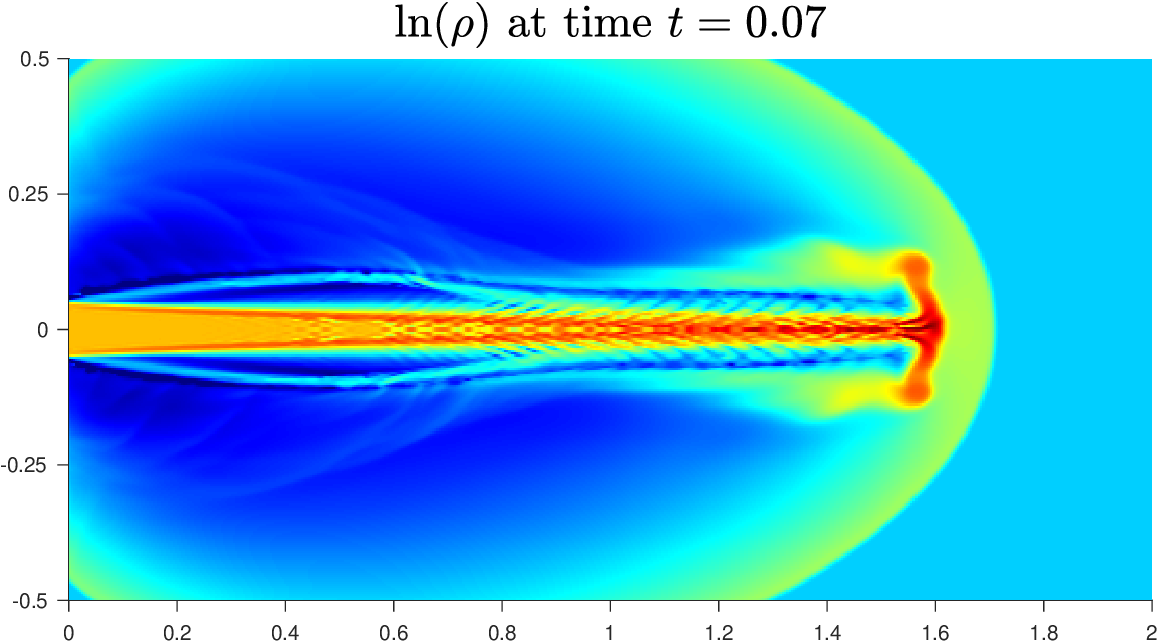}\hspace{-1.8mm} 
		\raisebox{+0.074\height}{
			\includegraphics[width=0.0224\linewidth]{./figs/EG2/Jet/colorbar_logdensity_jet.pdf}
		}\hspace{-1mm}
		\includegraphics[width=0.46\linewidth]{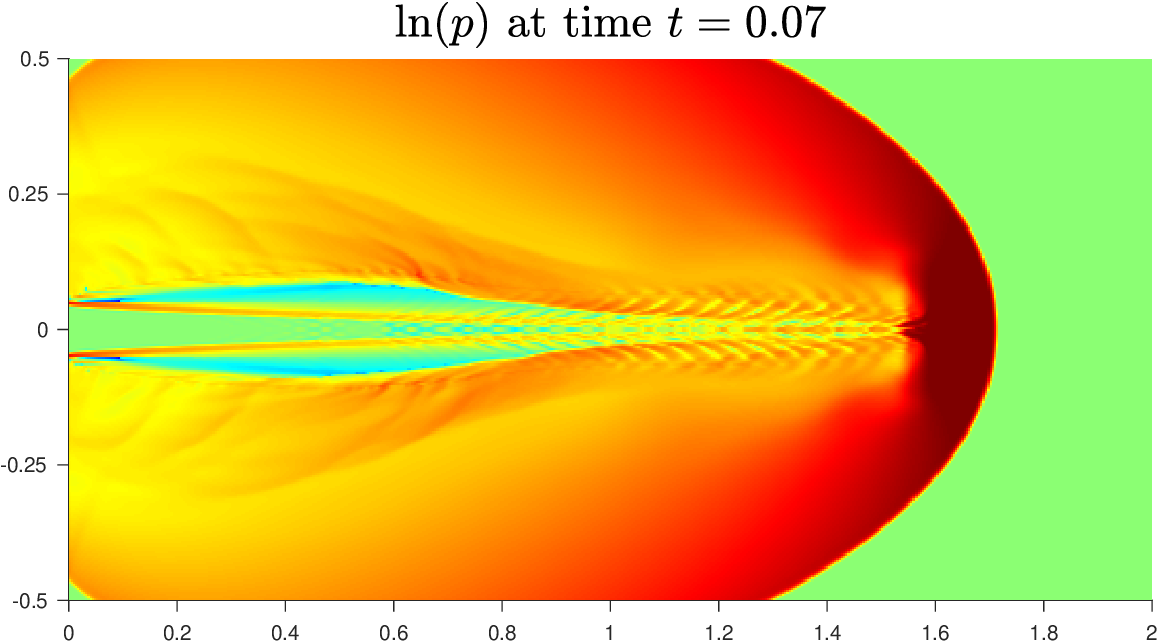}\hspace{-1.8mm} 
		\raisebox{+0.074\height}{
			\includegraphics[width=0.0224\linewidth]{./figs/EG2/Jet/colorbar_logpressure_jet.pdf}
		}
		\caption{\label{fig:Jet} Approximations of the Mach 80 jet at time $t=0.07$ with $512\times256$ cells. We show the logarithm of the density (left) and the logarithm of the pressure (right). The calculations were performed using AFEG2 (top), AFEG2$_{2.0}$ (middle), and AFCW (bottom) with $\mathrm{CFL}\le 0.23$.}
	\end{figure}
	
In Figure~\ref{fig:Jet}, we present the solutions computed with AFEG2, AFEG2$_{2.0}$, and AFCW. Since flux limiting is applied, we use the reduced CFL number $\mathrm{CFL}\leq 0.23$ for these methods.

		\subsection{Kelvin-Helmholtz instability}
	In this section, we compare the performance of AFEG2, AFEG2$_{1.0,1.0}^{0.4}$, and AFCW for a Kelvin-Helmholtz instability test problem. 
	We consider the following example.
	\begin{example}
		\label{ex:kelvin_helmholtz}
		We consider a rigorous multi-mode perturbation setup on the computational domain $[0,1]^2$. The initial condition is given by:
		\begin{equation}
			(\rho, u, v, p)(x, y, 0) = \begin{cases}
				(2, -0.5, 0, 2.5) & \text{if } I_1(x) \le y \le I_2(x), \\
				(1, 0.5, 0, 2.5) & \text{otherwise},
			\end{cases}
		\end{equation}
		where the interface profiles $I_j(x) = J_j + \epsilon Y_j(x)$ for $j=1,2$ are chosen as small perturbations around the lower interface $J_1 = 0.25$ and the upper interface $J_2 = 0.75$.
		
		To trigger the development of complex turbulent structures, a multi-mode perturbation function is applied to the interfaces:
		\begin{equation}
			Y_j(x) = \sum_{k=1}^m a_j^k \cos(b_j^k + 2k\pi x), \quad j = 1, 2.
		\end{equation}
		Here, the mode parameters $a_j^k \in [0,1]$ and phase shifts $b_j^k$ are arbitrary but fixed numbers. The amplitude coefficients $a_j^k$ are normalized such that $\sum_{k=1}^m a_j^k = 1$, which guarantees that the interface displacement strictly satisfies $|I_j(x) - J_j| \le \epsilon$. For this simulation, we use $m = 10$ perturbation modes and a perturbation amplitude $\epsilon = 0.01$. We use periodic boundary conditions in both spatial directions.
	\end{example}

	\begin{figure}
		\makebox[\textwidth][c]{\includegraphics[width=0.3\textwidth]{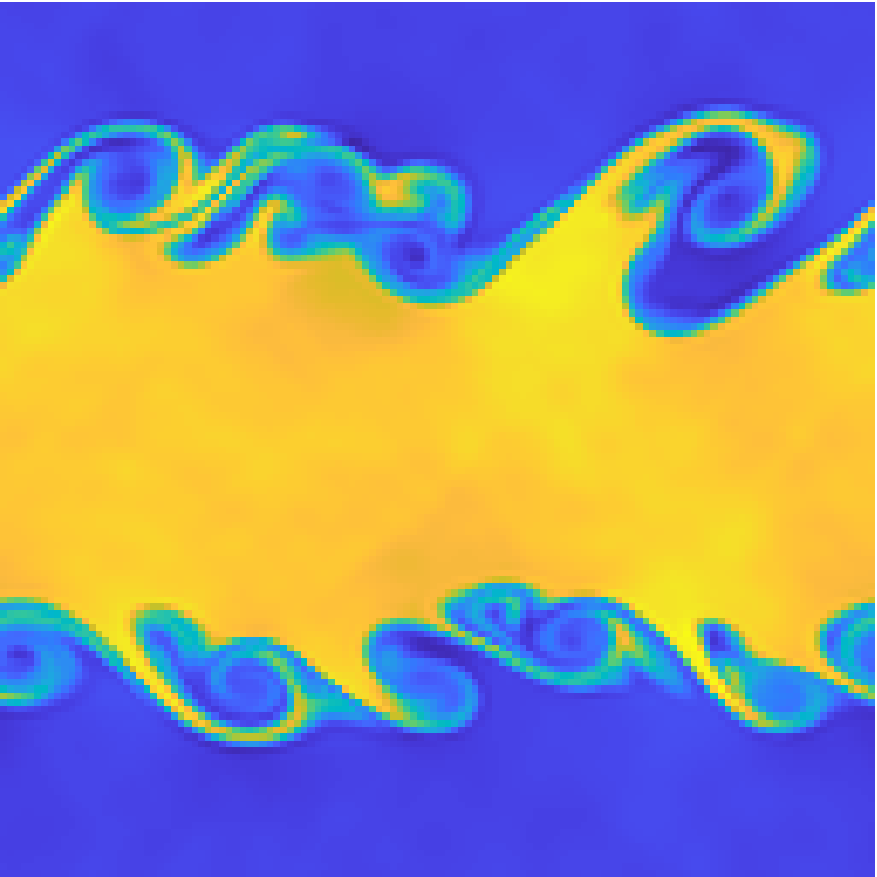}
			\includegraphics[width=0.3\textwidth]{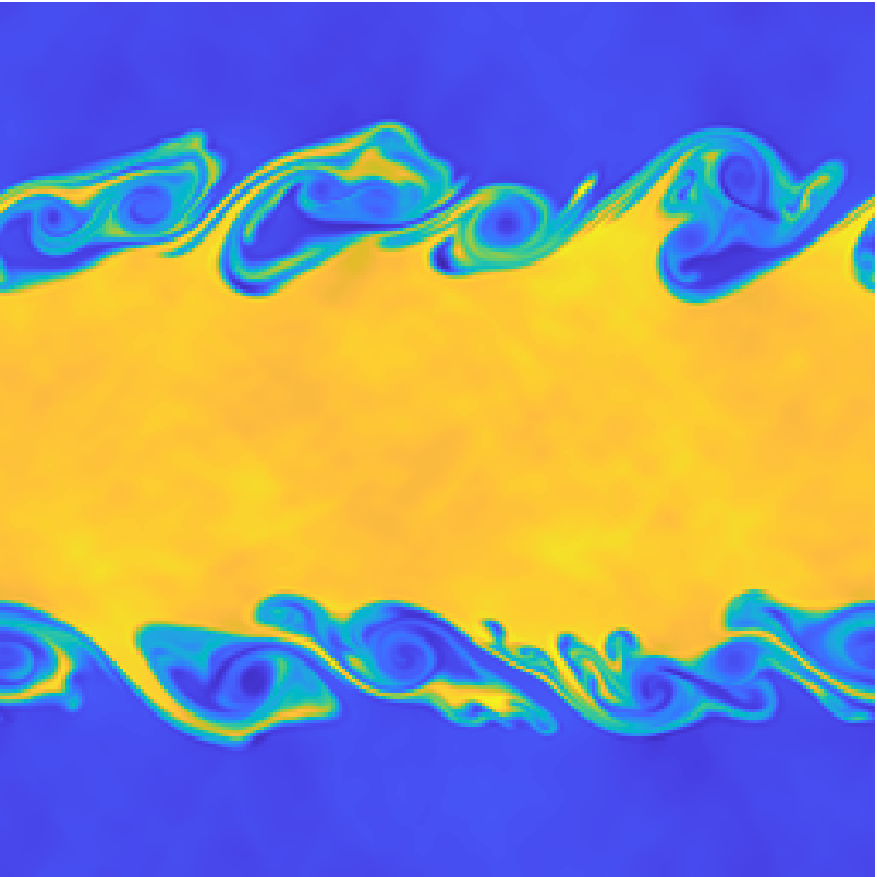} 
			\includegraphics[width=0.3\textwidth]{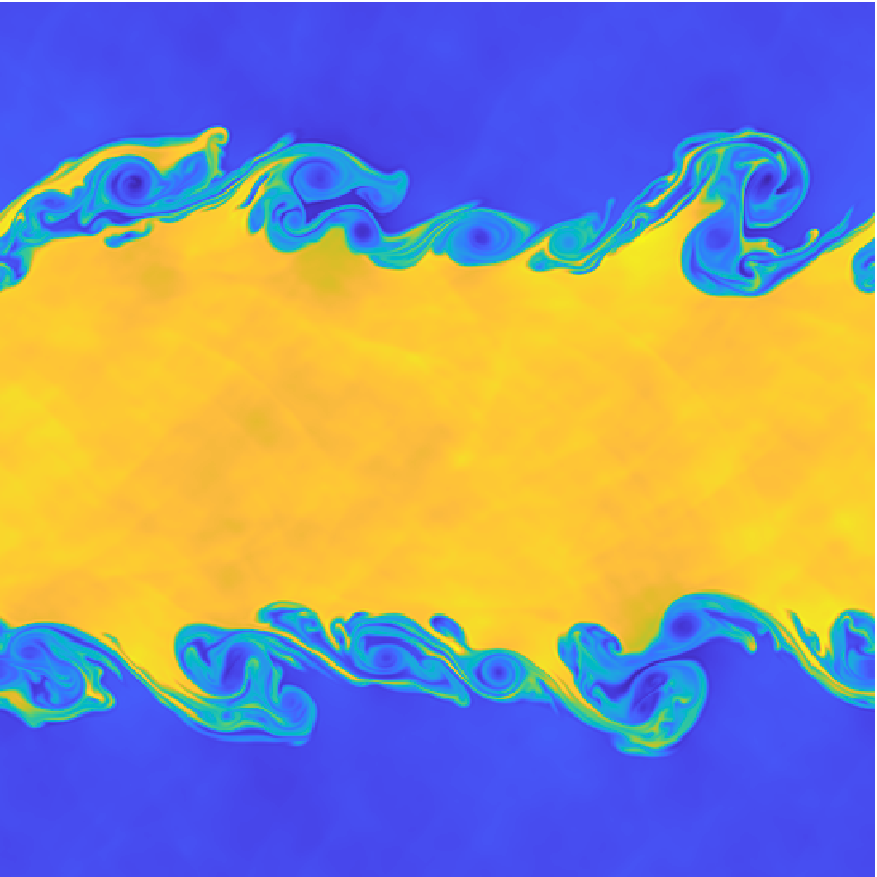}}\\[4pt]
		\makebox[\textwidth][c]{\includegraphics[width=0.3\textwidth]{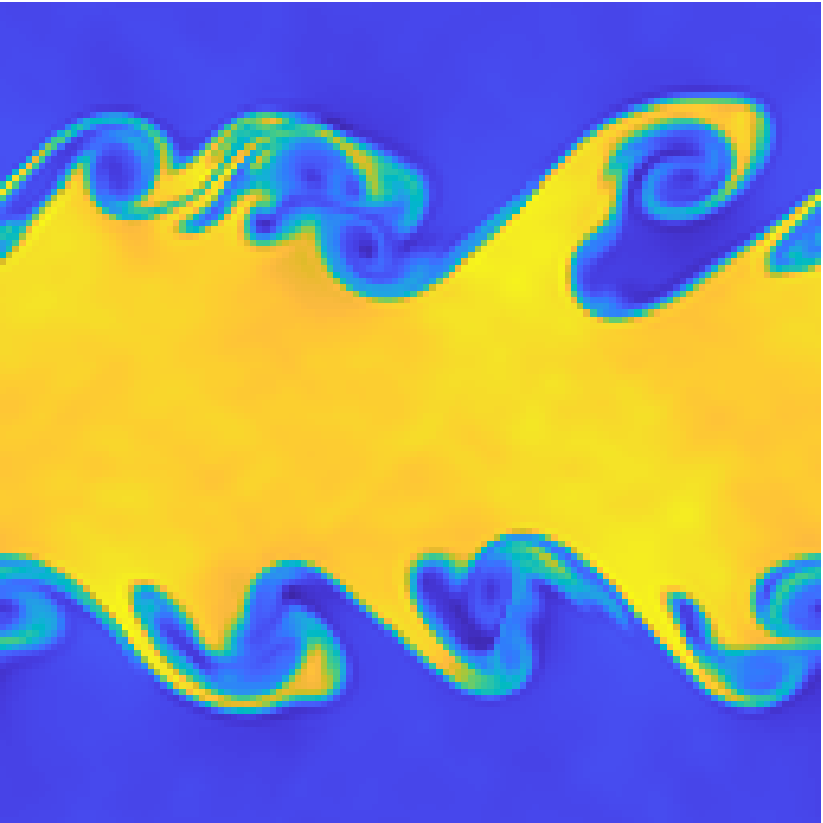}
			\includegraphics[width=0.3\textwidth]{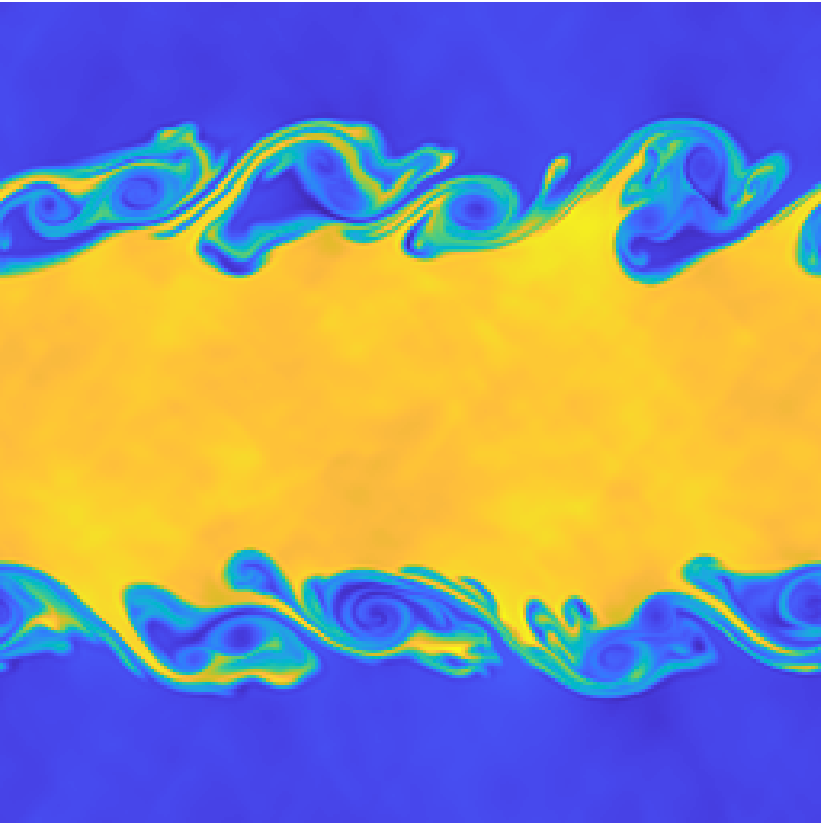}
			\includegraphics[width=0.3\textwidth]{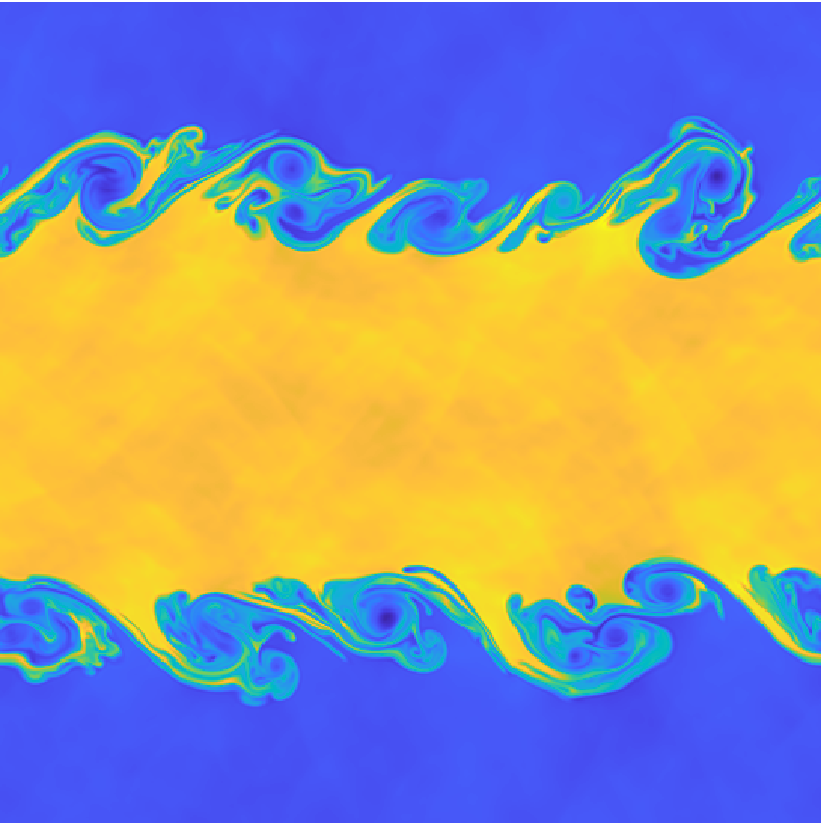}}\\[4pt]
		\makebox[\textwidth][c]{\includegraphics[width=0.3\textwidth]{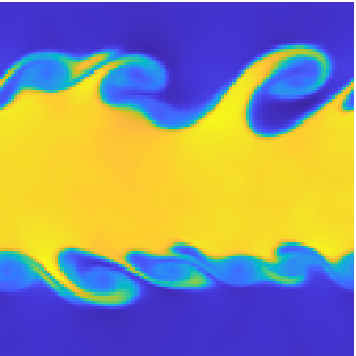}
			\includegraphics[width=0.3\textwidth]{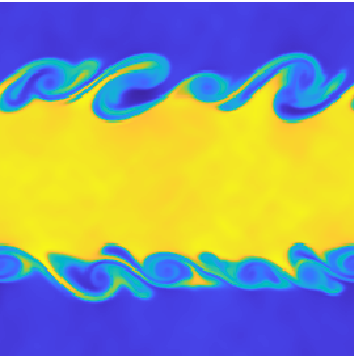}
			\includegraphics[width=0.3\textwidth]{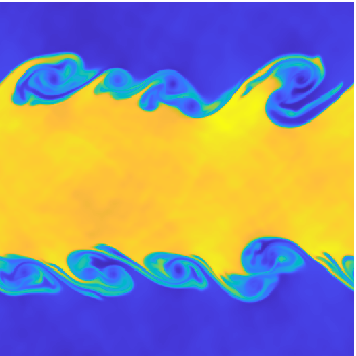}}
		\caption{\label{fig:KHrand1} Approximation of Example \ref{ex:kelvin_helmholtz} at $t=1$ on grids with $128^2$ (left), $256^2$ (middle), and $512^2$ (right) cells, computed using AFEG2 (top row, $\mathrm{CFL}\le 0.279$), AFEG2$_{1.0,1.0}^{0.4}$ (middle row, $\mathrm{CFL}\le 0.444$), and AFCW (bottom row, $\mathrm{CFL}\le 0.7$). No limiting was used.}
		\end{figure}
		
		\begin{figure}
		
		\makebox[\textwidth][c]{\includegraphics[width=0.3\textwidth]{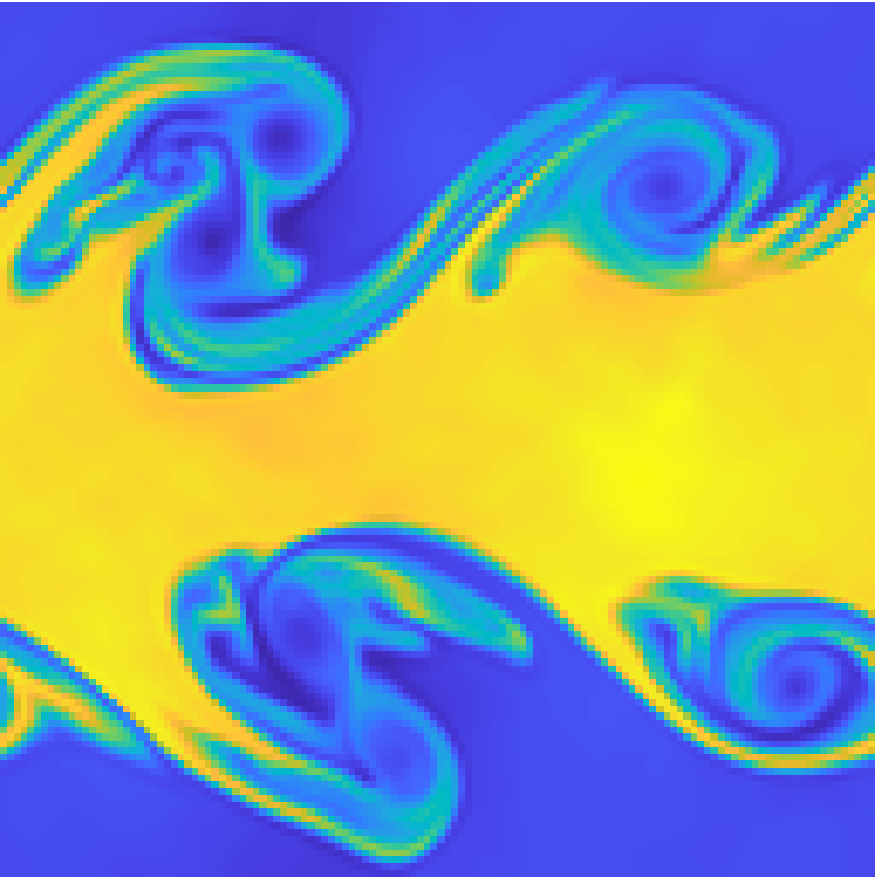}
			\includegraphics[width=0.3\textwidth]{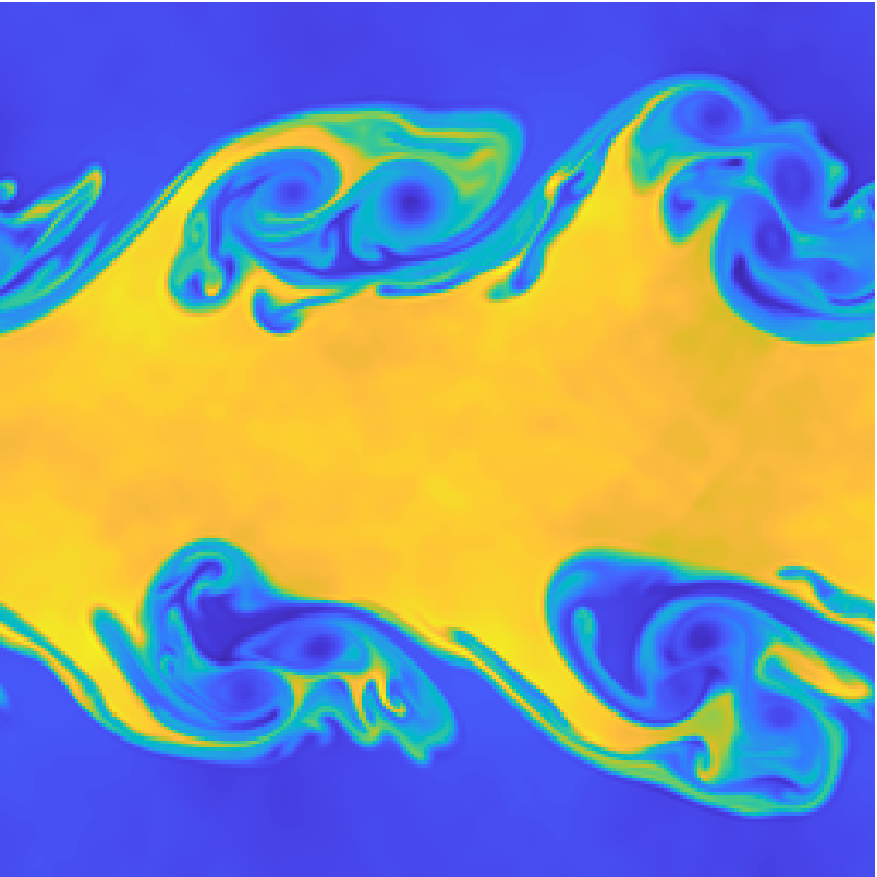} 
			\includegraphics[width=0.3\textwidth]{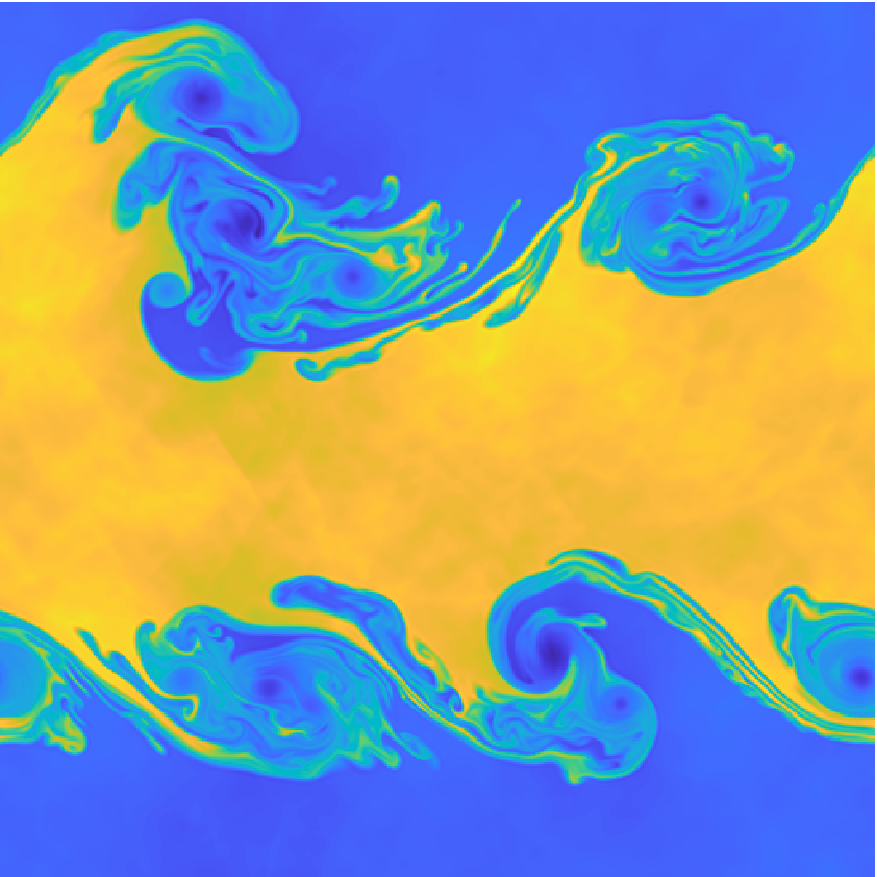}}\\[4pt]
		\makebox[\textwidth][c]{\includegraphics[width=0.3\textwidth]{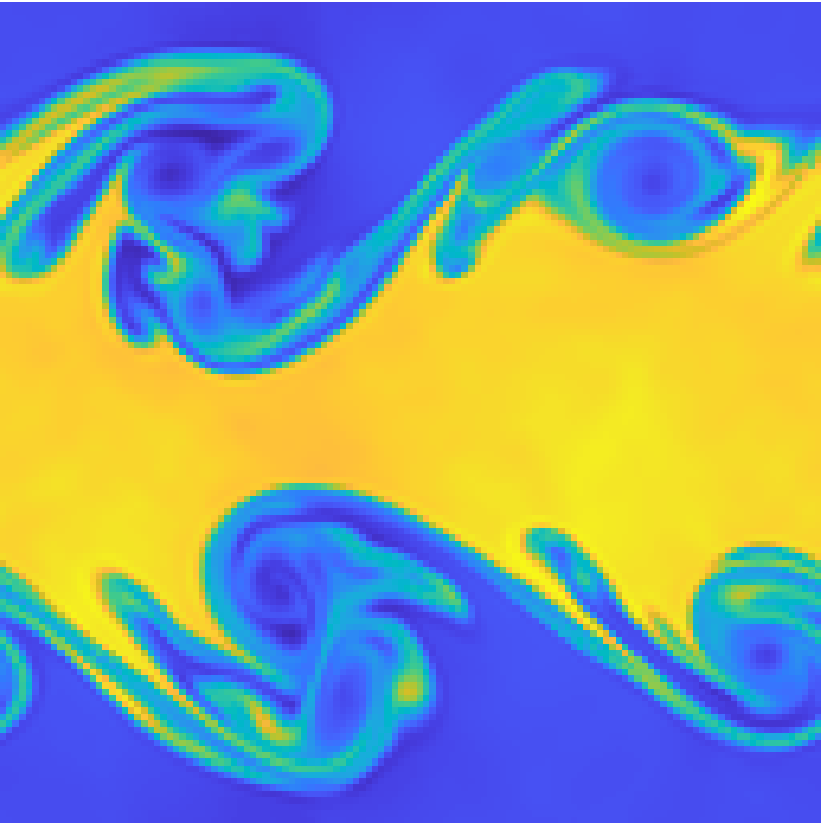}
			\includegraphics[width=0.3\textwidth]{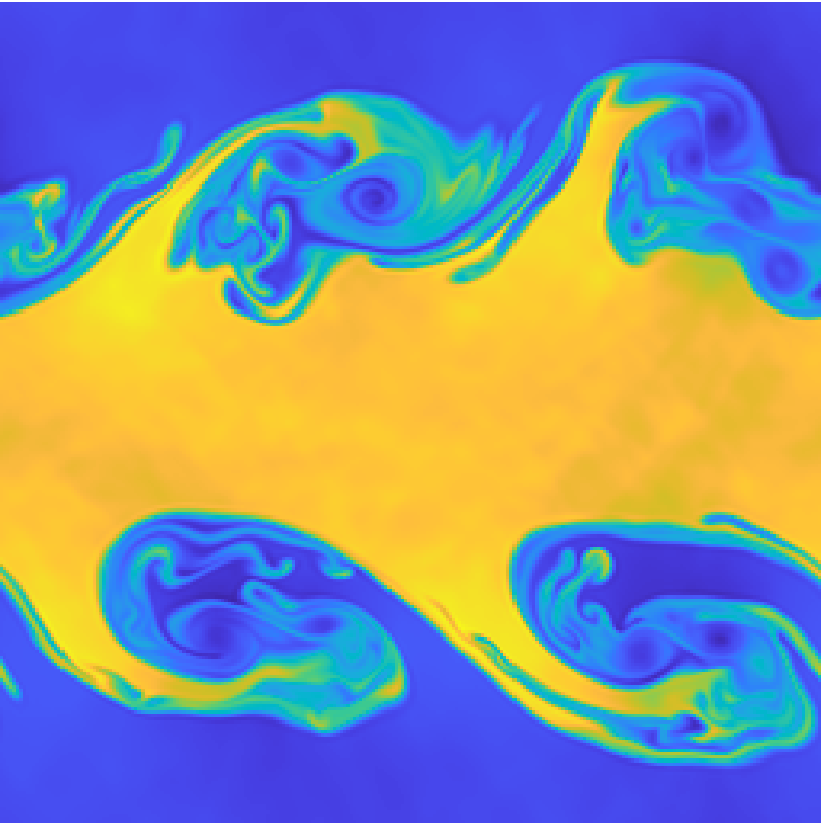}
			\includegraphics[width=0.3\textwidth]{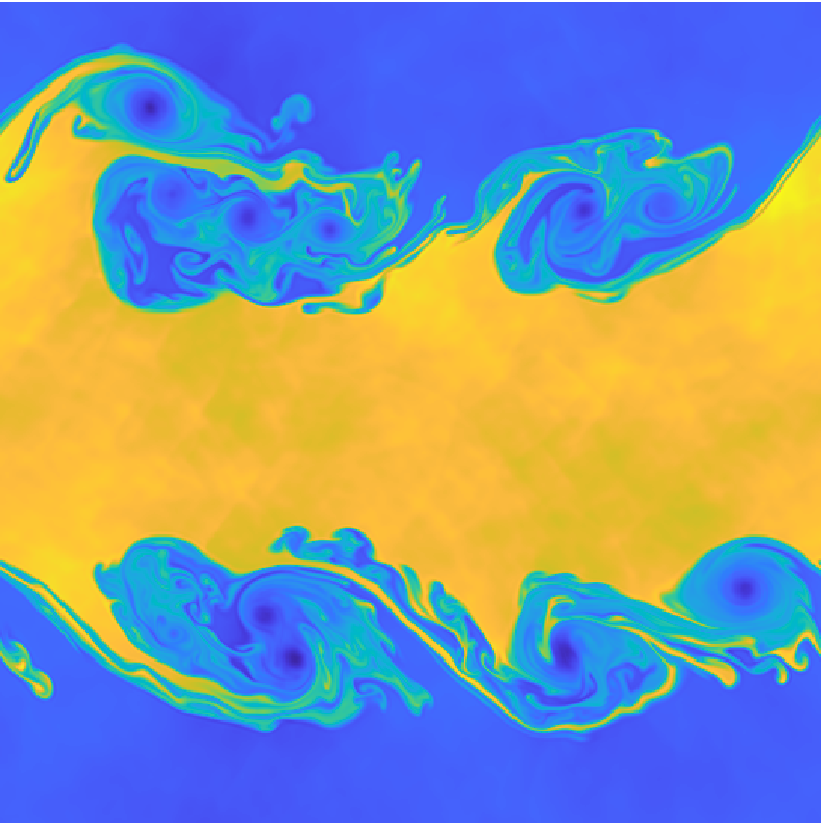}}\\[4pt]
		\makebox[\textwidth][c]{\includegraphics[width=0.3\textwidth]{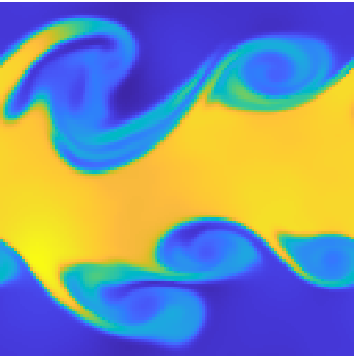}
			\includegraphics[width=0.3\textwidth]{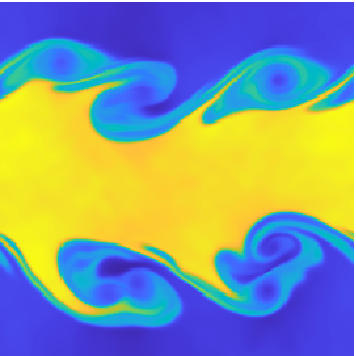}
			\includegraphics[width=0.3\textwidth]{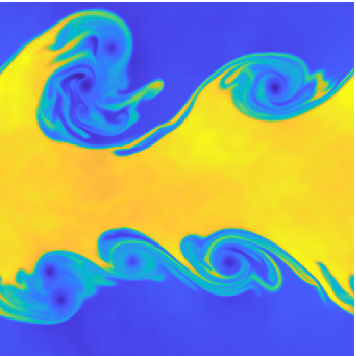}}
		\caption{\label{fig:KHrand2} Approximation of Example \ref{ex:kelvin_helmholtz} at $t=2$ on grids with $128^2$ (left), $256^2$ (middle), and $512^2$ (right) cells, computed using AFEG2 (top row, $\mathrm{CFL}\le 0.279$), AFEG2$_{1.0,1.0}^{0.4}$ (middle row, $\mathrm{CFL}\le 0.444$), and AFCW (bottom row, $\mathrm{CFL}\le 0.7$). No limiting was used.}
		\end{figure}
		
		In Figure \ref{fig:KHrand1} we show the solutions at $t=1$ calculated on different grids using AFEG2, AFEG2$_{1.0,1.0}^{0.4}$, and AFCW, while Figure \ref{fig:KHrand2} shows the corresponding solutions at $t=2$. Again, the time-step length was chosen close to the stability limit of each method. At the earlier time, only minor differences can be observed, whereas at the later time, the differences become more pronounced.

        To complement the visual comparison, we evaluate the entropy
and energy diagnostics considered in \cite{ChuEtAl2026}.
For each method, we average the numerical solutions on
$N^2$ grids with
$N = \{64,128,256,512\}$ after mapping them to
the common $512^2$ grid:
\[
\bigl(
\widetilde{\rho},\widetilde{\mathbf m},
\widetilde{S},\widetilde{E}
\bigr)
=
\frac{1}{4}
\sum_{N\in\mathcal N}
\mathcal P_N
\bigl(\rho_N,\mathbf m_N,S_N,E_N\bigr).
\]
Here, $\mathcal P_N$ denotes the mapping to the common grid,
$\mathbf m_N=\rho_N(u_N,v_N)$ is the momentum density,
$E_N$ is the total energy density, and
\[
S_N=c_v\rho_N
\log\left(\frac{p_N}{\rho_N^\gamma}\right),
\qquad c_v=\frac{1}{\gamma-1},
\]
is the entropy density.
The average over grid resolutions is computed at the output times $t=0,1,2$.

Writing the total energy density in terms of density,
momentum, and entropy as
\[
\mathcal E(\rho,\mathbf m,S)
=
\frac{|\mathbf m|^2}{2\rho}
+
\frac{\rho^\gamma}{\gamma-1}
\exp\left(\frac{(\gamma-1)S}{\rho}\right),
\]
we define
\begin{align*}
S(t)
&=\int_\Omega \widetilde S\,\mathrm{d}x\,\mathrm{d}y,
&
E_1(t)
&=\int_\Omega \widetilde E\,\mathrm{d}x\,\mathrm{d}y,
\\
E_2(t)
&=\int_\Omega
\mathcal E\bigl(
\widetilde\rho,\widetilde{\mathbf m},\widetilde S
\bigr)\,\mathrm{d}x\,\mathrm{d}y,
&
D_E(t)
&=\int_\Omega
\left|
\widetilde E-
\mathcal E\bigl(
\widetilde\rho,\widetilde{\mathbf m},\widetilde S
\bigr)
\right|\,\mathrm{d}x\,\mathrm{d}y.
\end{align*}
Thus, $E_1$ is the energy averaged over grid resolutions,
whereas $E_2$ is the energy of the averaged state.
The quantity $D_E$ measures the energy defect. 


Following the numerical diagnostics considered in \cite{ChuEtAl2026}, 
we also evaluate
\[
J(t)=E_2(t)-\overline{\theta} S(t),
\qquad
\overline{\theta}
=
\frac{\displaystyle
\int_\Omega E_0\,\mathrm{d}x\,\mathrm{d}y}
{\displaystyle
c_v\int_\Omega \rho_0\,\mathrm{d}x\,\mathrm{d}y},
\]
where $\rho_0$ and $E_0$ denote the initial density
and total energy density, respectively.
The quantity $J$ is motivated by the relative energy criterion
introduced in \cite{FL25}.
Minimizing $J$ favors the solution closest to the
maximum-entropy equilibrium as measured by the Bregman distance
(relative energy). In \cite{FL25}, minimizing $J$ or maximizing $S$ and $D_E$ are proposed to select a unique solution and to recover well-posedness of the compressible Euler equations; we also refer to \cite{FJL} for further related results.
 
Table~\ref{tab:kh_mean_diagnostics} reports the finite-time averages
\[
\langle X\rangle_T
=
\frac{1}{T}\int_0^T X(t)\,\mathrm{d}t,
\qquad
X\in\{S,E_1,E_2,D_E,J\},
\]
over the simulation interval $[0,T]$.
All of our methods capture the roll-up of the perturbed interfaces and reveal progressively finer vortical structures as the grid is refined. The differences are more visible at $T=2$.
The averaged energies $\mathcal{E}_1$ are nearly identical,
whereas the entropy and energy defect show more pronounced
differences.
AFCW yields the largest time-averaged entropy and the smallest
energy defect and is therefore favored by the maximum-entropy
and minimum-defect criteria.
In contrast, AFEG2$^{0.4}_{1.0,1.0}$ yields the largest energy
defect and is favored by the maximum-defect criterion.
This larger defect represents a greater energy contribution
from oscillations not captured by the averaged state.

AFEG2$^{0.4}_{1.0,1.0}$ also yields the smallest Bregman
distance; note, however, that the differences between the methods
are small.  For AFEG2$^{0.4}_{1.0,1.0}$, the lower mean-state energy
outweighs the effect of its lower entropy.
Thus, entropy maximization and the Bregman distance minimization favor
different methods in this comparison.

These observations concern the computed solutions at the
respective CFL numbers used for each method.
The differences therefore reflect the combined effects of
the spatial and temporal discretizations and cannot be
attributed solely to the reconstruction or evolution operator.
The resulting rankings describe preferences under the
specified selection criteria, rather than an overall
ranking of numerical accuracy.
	\begin{table}[htbp]
		\centering
		\caption{Kelvin--Helmholtz instability: Time averages of the entropy, averaged energy, mean-state energy, energy defect, and Bregman distance for AFEG2, AF$\mathrm{EG2}_{1.0,1.0}^{0.4}$ and AFCW.}
		\label{tab:kh_mean_diagnostics}
		\begin{tabular}{cccccc}
			\hline
			Operator 
			& \(\displaystyle \frac{1}{T}\int_0^T S(t)\,dt\)
			& \(\displaystyle \frac{1}{T}\int_0^T E_1(t)\,dt\)
			& \(\displaystyle \frac{1}{T}\int_0^T E_2(t)\,dt\)
			& \(\displaystyle \frac{1}{T}\int_0^T D_E(t)\,dt\)
			& \(\displaystyle \frac{1}{T}\int_0^T J(t)\,dt\) \\
			\hline
			EG2    
			& \(1.050382\) 
			& \(6.437890\) 
			& \(6.380600\) 
			& \(5.728970\times 10^{-2}\) 
			& \(4.577404\) \\
			$\mathrm{EG2}_{1.0,1.0}^{0.4}$
			& \(1.049953\) 
			& \(6.437898\) 
			& \(6.375874\) 
			& \(6.202381\times 10^{-2}\) 
			& \(4.573414\) \\
			AFCW
			& \(1.058007\)
			& \(6.437506\)
			& \(6.394706\)
			& \(4.293490\times 10^{-2}\)
			&\(4.578420\) \\
			\hline
		\end{tabular}
	\end{table}

	\ignore{
		The second example was proposed in \cite{article:LABH2024}. 
		\begin{example}\label{ex:KH2}
			We consider the two-dimensional Euler equations with initial values of
			the form
			\begin{align*}
				\rho _0(x,y) = &\gamma + R(1-2\eta(y))\\
				u _0(x,y) =& M(1-2\eta(y))\\
				v_0 (x,y) = & \delta M \sin(2\pi x)\\		
				p_0 (x,y)  = & 1,
			\end{align*}
			with 
			\begin{equation*}
				\eta (y)= 
				\left\{ \begin{array}{lcl}
					\frac{1}{2}(1+\sin (16 \pi (y+\frac{1}{4}))) & : & -\frac{9}{32}\le y<-\frac{7}{32},\\
					1 & : & -\frac{7}{32}\le y<\frac{7}{32},\\
					\frac{1}{2}(1+\sin (16 \pi (y+\frac{1}{4}))) & : & \frac{7}{32}\le y<\frac{9}{32},\\
					0 & : & else\\\end{array}\right.
			\end{equation*}
			and $R=10^{-3}$, $\delta = 0.1$ and $M = 0.01$, which is governing the Mach number of the flow. The domain is $[0,2]\times[-\frac{1}{2},\frac{1}{2}]$ with double periodic boundary conditions.  
		\end{example}
		In Figure \ref{fig:KH2} we show approximations at $t = 0.8/M = 80$ again calculated on different grids with the methods AFEG2, AFEG2$_{0.7}$ and  AFEG2$_{1.0,0,4}$. We do not observe any differences in the solution structure. The same test case was also considered in \cite{preprint:BarsukowAFEuler}. While spurious structures appear in that work, we do not observe such behavior for our method.
		\begin{figure}
			\centering
			
			\includegraphics[width=0.32\textwidth]{./figs/KH2/EG2_t80_64x32.pdf} 
			\includegraphics[width=0.32\textwidth]{./figs/KH2/EG2_t80_128x64.pdf} 
			\includegraphics[width=0.32\textwidth]{./figs/EG2_t80_256x128.pdf}  \\
			\vspace{4pt}
			\includegraphics[width=0.32\textwidth]{./figs/KH2/EG2_10_04_t80_64x32.pdf}
			\includegraphics[width=0.32\textwidth]{./figs/KH2/EG2_10_04_t80_128x64.pdf} 
			\includegraphics[width=0.32\textwidth]{./figs/KH2/EG2_10_04_t80_256x128.pdf}

			\caption{\label{fig:KH2} Approximation of Example \ref{ex:KH2} a at $t=80$ on a grid with $64 \times 32$ (left), $128 \times 64$ (middle) and $256 \times 128$ (right)  grid cells with AFEG2 (top row) and AFEG2$_{1.0,0.4}$ (bottom row).}
			
	\end{figure}}
\FloatBarrier
\section*{Conclusions}
We presented several new results on third-order accurate, fully discrete AF and closely related Evolution Galerkin methods for the Euler equations on two-dimensional Cartesian grids.  A crucial component of our methods is the evolution operator that updates point values. This evolution is based on approximations of two-dimensional linearized Euler equations in primitive variables using the method of bicharacteristics. In this paper, we introduced several new third-order accurate evolution operators, i.e., $\mathrm{EG}^{\text{quad}}$, which is exact for quadratic plane waves, as well as the family of EG2$_{\nu,\delta}^{\omega}$ operators, motivated by the previously known EG2 operator. We explored how the choice of evolution operators affects stability and showed that the stability of the classical AF method with globally continuous reconstruction can be improved from CFL$\le0.279$ to CFL$\le0.444$.

Furthermore, we studied the influence of the reconstruction on the stability and accuracy. To this end, we replaced the globally continuous, piecewise quadratic AF reconstruction with a piecewise quadratic CWENO reconstruction based on neighboring cell averages. The CWENO reconstruction enlarges the stencil, improving stability and allowing time steps up to CFL$\le0.7$. 
For smooth test problems, we observe third-order convergence rates for both methods but noticeably smaller error values when using the globally continuous AF reconstruction. For the AF method with globally continuous reconstruction, limiting of point and cell average values was introduced in \cite{article:CHP2026} and used here for all test problems. When using the AFCW method, limiting of point values was only needed for the Sod shock tube problem, the two-dimensional Riemann problem for Configuration 3 and the Mach 80 jet problem. For the other test problems, the CWENO reconstruction provided enough numerical viscosity to avoid spurious solution structures. While both methods produce comparable results, we observe a sharper approximation of contact waves when using the AF method with globally continuous reconstruction.

Finally, we observed an interesting difference for the Shu vortex problem. While the AFCW method provides enough numerical viscosity to accurately approximate this problem, the AF method with globally continuous reconstruction and EG2 evolution operator leads to an unstable approximation. The newly derived EG2$_{2.0}$ operator can be interpreted as an entropy-fix and leads to accurate results. 

In \cite{article:Roe2017}, Roe pointed out that in classical finite volume methods, a discontinuous reconstruction typically introduces one-dimensional physics. This is not the case in our AFCW method. Both of our methods are based on approximations of two-dimensional linearized Euler equations, making them truly multidimensional.  Our computations show that the more local stencil used in AF methods, where a globally continuous reconstruction is computed from point and cell average values, can yield more accurate approximations, particularly on coarse grids.  However, the lower numerical viscosity also makes the method less robust. This might be seen in stability constraints, which are quite sensitive to the precise choice of the evolution operator, the dependence of the method on the choice of the local linearization or difficulties in the approximation of the Shu vortex. The more viscous AFCW method was less sensitive in all respects.  
\section*{Acknowledgements}
This work was funded by DFG Projects 525800857 and 525853336 within the Priority Programme SPP 2410 ``Hyperbolic Balance Laws: Complexity, Scales and Randomness''. M.L.-M. gratefully acknowledges support of the Gutenberg Research College and of the Mainz Institute for Multiscale Modeling.

\bibliography{references}
\bibliographystyle{plain}

\end{document}